\input amstex
\documentstyle{amsppt}
\pageheight{18.3cm}
\magnification=1200
\nologo

\define\eps{\varepsilon}
\redefine\phi{\varphi}
\define\z{{\frak z}}
\define\k{{\frak k}}
\define\m{{\frak m}}
\define\LL{{\Cal L}}
\define\CC{{\Cal C}}
\define\J{{\Cal J}}
\define\vv{{\frak v}}
\define\g{{\frak g}}
\define\h{{\frak h}}
\define\zh{{\frak z(\frak h)}}
\define\phip{\phi'}
\define\psip{\psi'}
\define\uu{{\frak u}}
\define\prh{{\text{\rm{pr}}_{\frak h}}}
\define\fzb{\bar F_Z}
\define\HH{{\Cal H}}
\define\ph{{\varphi}}
\define\la{\lambda}
\define\lap{\lambda'}
\define\trans{{}^T}
\define\ab{_{a,b}}
\define\sig{{\sigma}}
\define\sab{{\frak S\ab}}
\define\subla{_\lambda}
\define\supla{^\lambda}
\define\sublap{_{\lambda'}}
\define\suplap{^{\lambda'}}
\define\om{\omega}
\define\subom{_\omega}
\define\supom{^\omega}
\define\omp{\omega'}
\define\subomp{_{\omega'}}
\define\suboii{_{\om,\phi,\psi}}
\define\suboiip{_{\omp,\phip,\psip}}
\define\sublii{_{\la,\phi,\psi}}
\define\subliip{_{\lap,\phip,\psip}}
\define\subpz{_{(p,z)}}
\define\scp{{\<\,.\,,.\,\>}}
\define\inv{^{-1}}
\define\Sym{\text{\rm{Sym}}}
\define\ourspace{{\Sym^2(\R^3)^*\otimes\R^3}}

\define\End{{\text{\rm{End}}}}
\define\Ad{{\text{\rm{Ad}}}}
\define\ad{{\text{\rm{ad}}}}
\define\Id{{\text{\rm{Id}}}}
\define\SU{{\text{\rm{SU}}}}
\define\su{{\text{$\frak s$}\text{$\frak u$}}}
\define\gl{{\text{$\frak g$}\text{$\frak l$}}}
\redefine\O{{\text{\rm{O}}}}
\define\SO{{\text{\rm{SO}}}}
\define\so{{\text{$\frak s$}\text{$\frak o$}}}
\define\Spin{{\text{\rm{Spin}}}}
\define\voll{{\text{\rm{vol}}}}
\define\dvol{{\text{{\it dvol}}}}
\define\spec{{\text{\rm{spec}}}}
\define\spann{{\text{\rm{span}}\,}}
\define\kernn{{\text{\rm{ker}}\,}}
\define\Imm{{\text{\rm{Im}}\,}}
\define\dimm{{\text{\rm{dim}}}}
\define\grad{{\text{\rm{grad}}}}
\redefine\exp{{\text{\rm{exp}}}}
\define\expst{{\text{\rm{exp}}}^*}
\define\tr{{\text{\rm{tr}}}}
\define\scal{{\text{\rm{scal}}}}
\define\Ric{{\text{\rm{Ric}}}}
\define\Isom{{\text{\rm{Isom}}}}
\define\const{\text{{\rm const}}}
\define\N{{\Bbb N}}
\define\R{{\Bbb R}}
\define\C{{\Bbb C}}

\define\<{{\langle}}
\define\>{{\rangle}}
\define\minzero{\setminus\{0\}}

\define\restr#1{{\lower0.4ex\hbox{$\vert$}\lower0.7ex
  \hbox{$\ssize{#1}$}}}
\define\drestr#1{{\lower0.6ex\hbox{$\vert$}\lower1ex
  \hbox{$\ssize{#1}$}}}
\define\tdiffzero#1{{\frac d{d{#1}}\lower0.2ex\hbox{\restr{{#1}=0}}\,}}
\define\diffzerosec#1{{\frac{d^2}{d{#1}^2}\lower1ex\hbox{$\big
  \vert_{\raise1\jot\hbox{${\tsize {#1}=0}$}}$}\,}}

\define\kapitel#1{\hskip 15pt\llap{{\bf#1}\hbox to6pt{}}}
\define\subkapitel#1{\hskip 34pt\llap{#1\hbox to6pt{}}}
\define\subsubkapitel#1{\hskip 61pt\llap{#1\hbox to 6pt{}}}

\font\gross=cmbx12

\topmatter

\title Isospectral manifolds with different local geometries\endtitle
\author Dorothee Schueth\endauthor
\address Mathematisches Institut, Universit\"at Bonn, Beringstr\.~1,
  D-53115 Bonn, Germany\endaddress
\email schueth\@math.uni-bonn.de\endemail
\thanks The author is partially supported by Sonderforschungsbereich 256,
  Bonn.\endthanks
\abstract
We construct several new classes of isospectral
manifolds with different local geometries.
After reviewing a theorem by Carolyn Gordon on isospectral torus
bundles and presenting certain useful specialized versions (Chapter~1)
we apply these tools to construct the first examples of isospectral
four-dimensional manifolds which are not locally isometric (Chapter~2).
Moreover, we construct the first examples of isospectral left invariant
metrics on compact Lie groups (Chapter~3). Thereby we also obtain the
first continuous isospectral families of globally homogeneous
manifolds and the first examples of isospectral manifolds which are
simply connected and irreducible.
Finally, we construct the first pairs of isospectral manifolds which
are conformally equivalent and not locally isometric (Chapter~4).
\endabstract

\toc\widestnumber\subsubhead{3.2.3}
\specialhead {}Introduction\endspecialhead
\head 1. Constructions of isospectral, locally non-isometric
  manifolds\endhead
\subhead 1.1 Isospectral torus bundles with totally geodesic fibers
  \endsubhead
\subhead 1.2 A first specialization and its application to
  products\endsubhead
\subhead 1.3 Review of some previously known examples\endsubhead
\head 2. Four-dimensional examples\endhead
\subhead 2.1 Isospectral, locally non-isometric metrics on
  $S^2\times T^2$\endsubhead
\subhead 2.2 Curvature properties\endsubhead
\head 3. Isospectral left invariant metrics on compact Lie groups\endhead
\subhead 3.1 Application of the torus bundle construction
  to compact Lie groups\endsubhead
\subhead 3.2 Examples\endsubhead
\subsubhead 3.2.1 Isospectral deformations on
  $\text{SO}(m)\times T^2$ ($m\ge5$),\newline
  $\text{Spin}(m)\times T^2$ ($m\ge5$),
  and $\text{SU}(m)\times T^2$ ($m\ge3$)\endsubsubhead
\subsubhead 3.2.2 Isospectral deformations on $\text{SO}(n)$
  ($n\ge9$), $\text{Spin}(n)$ ($n\ge9$),
  and $\text{SU}(n)$ ($n\ge6$)\endsubsubhead
\subsubhead 3.2.3 Isospectral deformations on $\text{SO}(8)$
  and $\text{Spin}(8)$\endsubsubhead
\subhead 3.3 Ricci curvature and $1$@-form heat invariants\endsubhead
\subhead 3.4 Infinitesimal spectral rigidity of bi-invariant
  metrics\endsubhead
\head 4. Conformally equivalent manifolds which are isospectral
   and not locally isometric\endhead
\subhead 4.1 Isospectral torus bundles whose fibers are not totally
  geodesic\endsubhead
\subhead 4.2 Conformally equivalent examples\endsubhead
\specialhead {}Bibliography\endspecialhead
\endtoc

\endtopmatter

\document

\noindent
{\gross Introduction} 

\bigskip

\noindent
In this work we construct several new classes of isospectral
manifolds with different local geometries. More precisely, we obtain
\roster
\item"{\bf--}" the first examples of four-dimensional isospectral manifolds
  which are not locally isometric (Chapter~2),
\item"{\bf--}" the first examples of isospectral left invariant metrics
  on compact Lie groups (Chapter~3) and the first examples of
  isospectral manifolds which are simply connected and irreducible,
\item"{\bf--}" the first examples of conformally equivalent manifolds
  which are isospectral and not locally isometric (Chapter~4).
\endroster

\medskip

The spectrum of a closed Riemannian manifold is the eigenvalue spectrum
of the associated Laplace operator acting on functions, counted with
multiplicities; two
manifolds are said to be isospectral if their spectra coincide.
Spectral geometry deals with the mutual influences between the spectrum
of a Riemannian manifold and its geometry (see the books \cite{1},
\cite{4}, \cite{11} for an introduction to spectral geometry).

Which geometric properties are determined by the spectrum?
Y.~Colin de Verdi\`ere showed that generically the Laplace spectrum
determines the spectrum of lengths of closed geodesics~\cite{12}.
Moreover, the spectrum determines a sequence of so-called heat invariants,
the first few of which are the dimension, the volume, and the total scalar
curvature (see, e.g., \cite{4}, \cite{18}).
A~few Riemannian manifolds are known to be completely characterized
by their spectra. For example, S.~Tanno proved this for the round spheres
in dimensions up to six~\cite{45}, using heat invariants;
for round spheres of arbitrary dimension it is only known that the
spectrum on functions together with the spectrum on $1$@-forms characterizes
them completely, as was shown by V.~Patodi~\cite{36}.
Moreover, there are several rigidity
and compactness results in special situations. S.~Tanno~\cite{46} showed
that every round sphere is infinitesimally spectrally rigid, that is,
one cannot continuously deform the round metric without changing the
spectrum.
C.~Croke and V.~Sharafutdinov proved infinitesimal spectral rigidity for
metrics of negative sectional curvature~\cite{14}.
B.~Osgood, R.~Phillips, and P.~Sarnak showed that
the set of metrics on a surface which are isospectral to a given metric
is always compact in the $C^\infty$@-topology~\cite{34}; the same
holds for bounded plane domains with respect to their Dirichlet
spectrum~\cite{35}. A~similar
result, although restricted to isospectral metrics within a fixed conformal
class, was shown for closed manifolds of dimension three
(\cite{10}, \cite{8}).
Moreover, R.~Brooks, P.~Perry, and P.~Petersen showed in~\cite{7}
that on a three-dimensional manifold every Riemannian metric which is
close enough to a metric of constant curvature has the property that
the corresponding set of isospectral metrics is compact.
H.~Pesce proved a compactness result in the case of a fixed Riemannian
covering: For any Riemannian manifold $(M,g)$ the set of discrete subgroups
$\Gamma<\Isom(M,g)$ for which the compact quotient manifold
$(\Gamma\backslash M,g)$ is isospectral to a given such manifold is compact
in the set of discrete subgroups of $\Isom(M,g)$ \cite{38}.

The general questions in these contexts, namely, whether ``most'' Riemannian
metrics are infinitesimally spectrally rigid, or whether each isospectral
set of metrics is compact in some appropriate topology, are still open.

\smallskip

On the other hand, many examples of isospectral manifolds have been
constructed, mainly during the last two decades. Note that the study
of such examples is the only possibility of finding geometric properties
which are {\it not\/} determined by the spectrum.
The first example of isospectral manifolds was given in 1964 by J.~Milnor:
a pair of flat tori in dimension sixteen~\cite{33} (by now there are also
examples of isospectral flat tori in dimension four~\cite{13}).
This was the first proof of the fact that the spectrum does not determine
the isometry class of a Riemannian manifold.
In 1980, M.-F.~Vign\'eras discovered examples of isospectral Riemann
surfaces and of isospectral hyperbolic manifolds in dimension
three, the latter showing that the fundamental group is not
spectrally determined (\cite{47}; see also P.~Buser's book~\cite{9} on
the spectral theory of Riemann surfaces).
The first examples of {\it continuous\/} families
of isospectral metrics were found by Carolyn Gordon and Edward Wilson
in 1984~\cite{25}; these were locally homogeneous metrics, induced by
left invariant ones, on compact quotients of nilpotent or solvable Lie
groups. In 1985 T.~Sunada established a general isospectrality
principle~\cite{43} which, either in its original or certain generalized
versions, came to be known as ``the Sunada method''. One generalized version,
established by Carolyn Gordon and Dennis DeTurck~\cite{16} in 1987,
is the following:
\roster
\item"\hbox{\qquad}"
  Let $(M,g)$ be a Riemannian manifold and $G<\Isom(M,g)$. Suppose
  $\Gamma_1\,,\Gamma_2$ are two discrete cocompact subgroups of~$G$ which
  act freely and properly discontinuously on~$M$ with compact quotients.
  If the quasi-regular representations of~$G$ on $L^2(\Gamma_1\backslash G)$
  and $L^2(\Gamma_2\backslash G)$ are unitarily equivalent, then
  $\Gamma_1\backslash M$ and $\Gamma_2\backslash M$, each endowed with
  the metric induced by~$g$, are isospectral.
\endroster
This theorem not only covered most of the isospectral examples known at
that time, but also led (via another generalization by
P.~B\'erard (\cite{2},\cite{3})
to the case of orbifolds) to the famous first examples of bounded plane
domains with the same Dirichlet (and Neumann) spectrum; these were found
in 1991 by C.~Gordon, D.~Webb, and S.~Wolpert~\cite{24}. Thereby,
M.~Kac's question of 1966, ``Can one hear the shape of a drum?''~\cite{32},
was finally answered negatively. Note, however, that these domains have
nonsmooth boundaries; the answer to Kac's question in the smoothly bounded
case is still open.

There is a generic converse to the Sunada theorem in the case where the
covering manifold~$M$ is compact:
H.~Pesce proved that there exists an open and dense set of metrics~$g$
on each closed manifold~$M$ such that all possible pairs of isospectral
quotients manifolds of the form $(\Gamma_1\backslash M,g)$,
$(\Gamma_2\backslash M,g)$ must necessarily arise from the Sunada
construction~\cite{39}.

By the very prinicple of the Sunada method described above, the isospectral
manifolds which arise from it always have a common Riemannian covering. In
particular, they are always locally isometric. Their geometries can
be distinguished only by global properties; for example, by the continuously
changing mass
of certain homology classes~\cite{17} or the changing distance of certain
geometrically distinguished families of geodesic loops~\cite{41}.

\smallskip

We now come to the history of isospectral manifolds which are not locally
isometric. In 1991, Zoltan Szab\'o discovered the first pairs of such
manifolds (see~\cite{44}, published much later); these were manifolds
with boundary, diffeomorphic to the product of an eight-dimensional ball
and a three-dimensional torus, arising as domains in quotients of certain
harmonic manifolds. Motivated by Szab\'o's examples, and related to them,
were the first pairs of isospectral manifolds without boundary which
Carolyn Gordon gave in 1992 (\cite{20}, \cite{21}); these were pairs
of two-step nilmanifolds with different underlying group structures.
Her isospectrality proof for these examples revealed another general
principle which is quite different from Sunada's and
does not imply local isometry of the resulting isospectral manifolds:
\roster
\item"\hbox{\qquad}"
  If a torus acts on two Riemannian manifolds freely and isometrically with
  totally geodesic fibers, and if the quotients of the manifolds by any
  subtorus of codimension at most one are isospectral when endowed with
  the submersion metric, then the original two manifolds are isospectral.
\endroster
Using this principle, C.~Gordon and E.~Wilson \cite{27}
generalized Z.~Szab\'o's
examples and obtained continuous multiparameter families of isospectral,
locally non-isometric metrics on products of $(m\ge5)$@-dimensional balls
with $(r\ge2)$@-dimension\-al tori.
These arise as domains in certain Riemannian nilmanifolds whose Ricci
tensors have in general different eigenvalues.
Next, it was observed during a workshop in Grenoble in 1997 that the
boundaries of these manifolds are again isospectral (by the same principle)
and not locally isometric~\cite{22}. Among the isospectral families
discovered in this way are some where the maximal scalar curvature changes
during the deformation, which shows that the range of the scalar curvature
is not spectrally determined (although the total scalar curvature is).
Independently, Z.~Szab\'o showed, in a more special setting, that the
boundaries of his original examples are isospectral (see again~\cite{44}).
Among his examples is a pair of isospectral metrics one of which is
homogeneous while the other is not even locally homogeneous; thus (local)
homogeneity is not encoded in the spectrum. By embedding the torus
factor of the manifolds given in \cite{22} into a compact Lie group
and extending the metrics in such a way that the above isospectrality
principle for torus bundles still applied, the author constructed, also
in 1997, the first examples of simply connected, closed isospectral
manifolds~\cite{42}; note that non-simple connectivity was another
feature always present in the Sunada type examples of isospectral closed
manifolds.

In~\cite{42} the author also used the new examples to show
that the individual terms in the linear combination $5\int\scal^2-2\int
\|\Ric\|^2+2\int\|R\|^2$, which is one of the heat invariants, are not
spectrally
determined; more precisely, the corresponding heat invariant for the Laplace
operator acting on $1$@-forms, which is a different linear combination of
the same terms, changes during the isospectral deformations given
in~\cite{42}. In particular, these manifolds are not isospectral on
$1$@-forms.
Before that, examples of manifolds which are isospectral on
functions but not on $1$@-forms had been given by A.~Ikeda~\cite{31}
(lens spaces, see also \cite{30}, Carolyn Gordon~\cite{19}
(Heisenberg manifolds, see also~\cite{26}), and
Ruth Gornet (\cite{28}, \cite{29}) who also constructed the first
continuous families with this property. Note that none of those examples
arose from the Sunada method, but used special constructions. In fact,
the Sunada setting 
--- except for a certain further generalization of it established by H.~Pesce
\cite{37} which also explains Ikeda's examples ---
always implies isospectrality not only on functions,
but also on all $p$@-forms~\cite{16}. Still, in the above examples by
Ikeda, Gordon, and Gornet, the isospectral manifolds do have a common
Riemannian covering.
It is not hard to see that integrals of functions which are induced
on manifolds of the same volume by isometry invariant functions
on a common covering manifold (such as $\scal^2$, $\|\Ric\|^2$, etc.)~must
always be the same. Thus isospectral manifolds with a common Riemannian
covering always share the same heat invariants, also for the Laplace operator
on $p$@-forms, even if they are not isospectral on $p$@-forms.
So it was the first time in~\cite{42} that heat invariants alone
were used to prove non-isospectrality on $1$@-forms.

Recently Carolyn Gordon and Zoltan Szab\'o gave a version of the isospectral
torus bundle theorem for the case where the fibers are not necessarily
totally geodesic, imposing certain other restrictions instead. In particular,
they obtained by this approach the first continuous families of negatively
curved isospectral manifolds with boundary~\cite{23}, contrasting with
the above rigidity result by Croke and Sharafutdinov~\cite{14} for the
case of closed manifolds.

\medskip

Thereby we finish our account of the previously known examples of isospectral
manifolds, and turn now to the description of the contents of the present
work.

\smallskip

In the preliminary Chapter~1 we first review Carolyn Gordon's above
principle of isospectral torus bundles with totally geodesic fibers
(Theorem~1.3) and give a slightly more special version (Theorem~1.6) in which
we assume the submersion quotients to be not only isospectral, but isometric,
and formulate this condition in terms of bundle connection forms with
which the metrics are associated. The formulation becomes
quite simple in the case of trivial bundles; i.e., products of the base
manifold and a torus (Proposition~1.8). We then review the examples
given in \cite{22} and \cite{42} and interpret them as applications
of Proposition~1.8 and Theorem~1.6, respectively.

\smallskip

Recall that the isospectral manifolds from \cite{22} were diffeomorphic
to $S^{m-1}\times T^2$ with $m\ge5$, and thus at least six-dimensional.
Examples of isospectral, locally non-isometric manifolds in lower dimensions
were not known until now. In Chapter~2 we use the point of view developed
in Chapter~1 on these previous examples in order to drop an unnecessary
property of their metrics and apply Proposition~1.8 in a systematic way
to construct isospectral, locally non-isometric metrics in dimension four,
namely, on $S^2\times T^2$ (see Example~2.6). We show that in all isospectral
pairs arising from Proposition~1.8 in which the base manifold is --- as here
--- two-dimensional, the associated scalar curvature functions share
the same range and the same integrals of each of their powers
(Theorem~2.11), which contrasts with the properties of the higher-dimensional
examples. Nevertheless, in one of our pairs of isospectral metrics on
$S^2\times T^2$ the preimages of the maximum of the associated scalar
curvature functions have different dimensions (Proposition 2.10), which
shows that the manifolds are not locally isometric. In many examples the
metrics can also be distinguished by the integral of $(\Delta\,\scal)^2$
(see Remark~2.12(iii)).

\smallskip

In Chapter~3 we apply Theorem~1.6\,/\,Propositon~1.8 to the case of
compact Lie groups with left invariant metrics
(Proposition~3.1\,/\,Corollary~3.2). The isospectrality result formulated
in Proposition~3.1 can be shown not only ``geometrically'' by deducing
it from Theorem~1.6 or~1.3, but also by purely algebraic methods involving
the expression of the Laplace operator associated to a left invariant
metric in terms of the right-regular representation of the Lie group.
We use Proposition~3.1 and Corollary~3.2 to construct continuous isospectral
families of left invariant metrics on $\SO(m\ge5)\times T^2$,
$\Spin(m\ge5)\times T^2$, $\SU(m\ge3)\times T^2$, $\SO(n\ge8)$,
$\Spin(n\ge8)$, and $\SU(n\ge6)$ (see the examples~3.3 and 3.7--3.10
in Section~3.2). These are not only the first examples of isospectral
left invariant metrics on compact Lie groups in general, but among them
are also the first examples of {\it irreducible\/} simply connected
isospectral manifolds.
Moreover, they are the first examples of continuous families of globally
homogeneous isospectral manifolds. We also obtain the first examples
of continuous isospectral families of manifolds of positive Ricci
curvature. In Section~3.3 we prove that the left
invariant metrics in our isospectral examples are not locally isometric; more
precisely, the norm of the associated Ricci tensors, $\|\Ric\|^2$ (which
is a constant function on each of the manifolds) changes during the
deformations (see Proposition~3.15 and Theorem~3.14). In particular,
a certain heat invariant for the Laplace operator on $1$@-forms also
changes during the deformations (Corollary~3.17).
Continuous isospectral families of left invariant metrics of the type
constructed in Section~3.2 can occur arbitrarily close to bi-invariant
metrics while still enjoying the non-isometry properties established
in Section~3.3 (see Remark 3.12(i)). In contrast to this, we prove in
Section~3.4 a rigidity result for bi-invariant metrics: Any continuous
isospectral family of left invariant metrics which contains a bi-invariant
metric must be trivial (Theorem~3.19).

\smallskip

In Chapter~4 we formulate a canonical generalization of Theorem~1.6
which in turn can be viewed as a special form of C.~Gordon's and
Z.~Szab\'o's above-mentioned isospectrality theorem for principal
torus bundles whose fibers are not assumed to be totally geodesic
(see Theorem~4.3 and Remark~4.4). Here again, we also give a simpler
formulation for the case of trivial bundles (Proposition~4.5).
We then use Proposition~4.5 to construct the first examples of conformally
equivalent isospectral manifolds with different local geometries
(Example~4.6). The non-isometry proof consists of showing that in our
pairs of conformally equivalent isospectral manifolds the respective
preimages of the maximal scalar curvature constitute a pair of isospectral,
globally homogeneous submanifolds of the type studied in Chapter~3,
whose Ricci tensors were already shown to have different norms
(Proposition~4.7).
Note that the first examples of isospectral, conformally
equivalent manifolds were constructed by Robert Brooks and Carolyn
Gordon~\cite{6} in 1990 using the Sunada method as formulated in~\cite{16};
in particular, those manifolds were locally isometric.
Our new examples are the first ones which show that even within a
fixed conformal class the local geometry is not determined by the spectrum.

\medskip
\noindent
{\bf Acknowledgments.}
First of all, it is a pleasure to thank Carolyn Gordon for many interesting
conversations over the past few years.
Also, there are many other members of the ``isospectral family'' to whom
I am grateful for stimulating meetings and discussions,
especially Pierre B\'erard, Robert Brooks, Peter Buser,
Ruth Gornet, Toshikazu Sunada, David Webb, Edward Wilson, and our
deceased friend and colleague Hubert Pesce.
Moreover, I thank Werner Ballmann, Ursula Hamenst\"adt, and Hermann
Karcher for their continuing interest and encouragement, and I thank the
Sonderforschungsbereich~256 at Bonn for partially supporting my research.
\hbox{Finally,} I am very grateful to my family for their love and support.

\bigskip

\bigskip
\noindent
\gross 1.\hbox to6pt{}Constructions of isospectral, locally
  non-isometric
\newline
\hbox to17pt{}manifolds\rm

\bigskip

\noindent
In this chapter we present the tools which we will use in Chapters~2 and~3
to construct new examples of isospectral manifolds which are not locally
isometric.

The starting point is a general theorem by Carolyn Gordon concerning
torus bundles with totally geodesic fibers; see Theorem~1.3 in Section~1.1.
We formulate somewhat more special versions
(Theorem~1.6 and Proposition~1.8 in Section~1.2) which account for almost
all previous applications of Theorem~1.3 and also for most of the
new examples of isospectral manifolds which
we present in this work. In Section~1.3 we review some previously
known examples which we interpret as applications of Theorem~1.6 or
Proposition~1.8.

\bigskip

\smallskip

\subhead 1.1\ \ Isospectral torus bundles with totally geodesic fibers
\endsubhead
\bigskip

\definition{Definition 1.1}
Let $(M,g)$ be a closed Riemannian manifold, and let~$\Delta_g$ be the
Laplacian acting on functions by
$$
(\Delta_g f)(p):=
-\sum_{i=1}^n \diffzerosec{t} f(c_i(t))\qquad\text{for $p\in M$,}
$$
where the~$c_i$ are geodesics starting in~$p$ such that $\{\dot c_1(0),
\dots,\dot c_n(0)\}$ is an orthonormal basis for~$T_pM$. The discrete
sequence $0=\la_0<\la_1\le\la_2\le\,\ldots\,\to\infty$ of the
eigenvalues of~$\Delta_g$\,, counted with the corresponding
multiplicities, is called the {\it spectrum} of $(M,g)$; we will denote
it by $\spec(M,g)$ or $\spec(\Delta_g)$. If $\HH\subseteq L^2(M,g)$ is
a subspace of functions such that $C^\infty(M)\cap\HH$ is invariant
under $\Delta_g$\,, we will
denote the corresponding spectrum of eigenvalues by $\spec(\HH)$.
Two closed Riemannian manifolds are called
{\it isospectral\/} if their spectra coincide.\enddefinition

\bigskip

All previously known examples of closed isospectral manifolds which are
isospectral and not locally isometric, except for some recent examples
by Gordon and Szab\'o \cite{23} (see Remark~4.4 in Section~4.1), arise from
the following theorem by Carolyn Gordon \cite{21} which we present below.

\bigskip

\subheading{Notation 1.2}
By a {\it torus\/}, we always mean
a compact connected abelian Lie group. If a torus $H$ acts
smoothly and freely by isometries on a closed Riemannian
manifold $(M,g)$, then there is
a unique Riemannian metric, denoted~$g^H$, on the quotient
manifold $M/H$ such that the canonical projection
$\pi_H:(M,g)\to(M/H,g^H)$ becomes a Riemannian submersion.

\bigskip

\proclaim{Theorem 1.3 ([Go3])}
Let $H$ be a torus, and let $(M,g)$ and $(M',g')$ be two
principal $H$@-bundles
such that the Riemannian metrics $g,g'$ are invariant under the action
of~$H$. Assume:
\roster
\item"(i)" The fibers of the action of~$H$ are totally geodesic
  submanifolds of $(M,g)$, resp\. of $(M',g')$.
\item"(ii)" For any closed subgroup $W$ of $H$ which is either $H$ itself
or a subtorus of codimension~$1$ in~$H$, the manifolds
$(M/W,g^W)$ and $(M/W,g^{\prime\,W})$ are isospectral.
\endroster
Then $(M,g)$ and $(M',g')$ are isospectral.\endproclaim

\smallskip

\demo{Proof}
We consider the unitary representation of~$H$ on the Hilbert space
$\HH:=L^2(M,g)$, defined by
$(zf)(x)=f(zx)$ for all $f\in\HH$, $z\in H$, $x\in M$.
Write $H=\h/\LL$, where $\h$ is isomorphic to some~$\R^r$,
and let $\LL^*$ be the dual lattice.
Since~$H$ is abelian, $\HH$ decomposes as the orthogonal sum
$\bigoplus_{\mu\in\LL^*}\HH_\mu$ with
$\HH_\mu=\{f\in\HH\mid zf=e^{2\pi i\mu(Z)}f\quad\text{for
all }z\in H\}$,
where $Z$ denotes any representative for~$z$ in~$\h$.
In particular, this implies the following coarser decomposition:
$$\HH=\HH_0\oplus{\tsize\bigoplus}_W(\HH_W\ominus\HH_0),\tag1
$$
where $W$ runs though the set of all closed connected subgroups of
codimension~$1$ in~$H$,
and $\HH_W$ is the sum of all $\HH_\mu$ such that $\mu\in\LL^*$
and $T_eW\subseteq\kernn\mu$.
Note that~$\HH_W$ is just the space of $W$@-invariant functions in~$\HH$.
Let $\CC_W:=C^\infty(M)\cap\HH_W$ and $\CC_0:=C^\infty(M)\cap\HH_0$\,.
Since the action of~$H$ is by isometries and therefore commutes
with $\Delta_g$\,, the spaces~$\CC_W$ and $\CC_0$ are invariant
under~$\Delta_g$\,.
Note that~$\pi_H^*$ is a linear bijection from
$C^\infty(M/H)$ to~$\CC_0$\,. Since $\pi_H$ is a Riemannian submersion with
totally geodesic fibers by~(i), $\pi_H^*$~intertwines the corresponding
Laplacians. Thus $\spec(\HH_0)=
\spec(M/H,g^H)$. Assumption~(ii) for $W=H$ implies, with the obvious
analogous notations for $(M',g')$, that $\spec(\HH_0) = \spec(\HH'_0)$.

Now let $W$ be a closed connected subgroup of codimension~$1$ in~$H$.
Then $\pi_W^*$ is a linear bijection
from $C^\infty(M/W)$ to $\CC_W$\,.
Note that assumption~(i) implies that also the $W$@-orbits are
totally geodesic, since the induced metric on the $H$@-orbits is
invariant and thus flat. This, together with assumption (ii) for~$W$
implies, by the same argument as before, that $\spec(\HH_W)=\spec(\HH'_W)$.
{}From~\thetag{1} we now conclude that $\spec(M,g)= \spec(M',g')$.
\qed\enddemo

\smallskip

\subheading{Remarks 1.4}

(i) Note that Theorem~1.3 is trivial if the torus~$H$ is one-dimensional:
In that case, condition (ii) on subtori of codimension~$1$ is equivalent
to the assertion of the theorem.

(ii) We remark that {\it all continuous\/} isospectral families
$(M,g_t)$ of principal torus bundles with totally geodesic fibers
must necessarily arise from Theorem~1.3.
In fact, if $H$ is the fiber of~$M$, then for each closed subgroup~$W$ of~$H$
the fibers of the $W$@-action are again totally geodesic in~$(M,g_t)$,
and thus the spectrum of $(M/W,g_t^W)$ is contained in the discrete set
$\spec(M,g_t)$ which is independent of~$t$ by assumption. Hence
$\spec(M/W,g_t^W)$ must itself be independent of~$t$. Therefore
also condition~(ii) of Theorem~1.3 is satisfied, which means that
the isospectral family $(M,g_t)$ arises from this theorem.

(iii) Although the isospectral manifolds arising from Theorem~1.3 are in
general not locally isometric, certain well-known families of isospectral,
locally isometric manifolds too can be viewed as applications of this
theorem.
For example, there are isospectral families $(\Phi_t(\Gamma)\backslash G,g)$
of two-step nilmanifolds~\cite{15} (where $G$ is a simply
connected two-step nilpotent Lie group, $g$~is a left invariant metric,
$\Gamma$~is a cocompact discrete sugroup of~$G$, and~$\Phi_t$ a family
of so-called almost inner automorphisms of~$G$). These arose originally
in the Sunada type context (see the Introduction and~\cite{25}, \cite{16}),
but can as well be viewed as arising from Theorem~1.3 (or the more special
Theorem~1.6 below). This is not only clear from~(ii) above, but in fact
the isospectrality proof given in~\cite{15} for these particular families
was already similar in vein to the proof of Theorem~1.3.

(iv) The observation in~(ii) illustrates that Theorem~1.3 is in fact a quite 
general result. Its ``broadness'' does not make it very obvious how
to find explicit applications --- except for the ones which motivated its
discovery in the first place (i.e., the two-step nilmanifolds from
\cite{21}) and their closely related companions from
\cite{22} and~\cite{42} (see Examples~1.11 and~1.14).

Therefore, we formulate in this work several specialized
versions of Theorem~1.3 (namely, Theorem~1.6 and Proposition~1.8 below,
as well as Proposition~3.1 and Corollary~3.2 in Chapter~3).
Using these we will be able to give interesting new applications.

(v) In almost all known applications of Theorem~1.3, the pairs of quotient
manifolds $(M/W,g^W)$ and $(M'/W, g^{\prime\,W})$ are not only
isospectral, but isometric. The only exceptions of this are certain pairs
of two-step nilmanifolds constructed in \cite{21} and \cite{27}.
The isospectrality of the (non-isometric) quotient manifolds follows
there from explicit knowledge of their spectra. We will not review
these examples in the present work.
All other known applications of Theorem~1.3 are actually applications of
Theorem~1.6 in the following section which is somewhat more special and
does imply isometric quotient manifolds, as we will see in the proof.

\bigskip

\smallskip

\subhead 1.2\ \ A first specialization and its application to products
  \endsubhead
\nopagebreak
\bigskip
\nopagebreak
\subheading{Notation and Remarks 1.5}
Let $H$ be a torus and $\h:=T_eH$.
When we call a metric on a torus {\it invariant\/} we will always
mean left invariant (which is here the same as bi-invariant).
Now let $H$ be equipped with a fixed
invariant metric. Let $M$ be a principal $H$@-bundle over a
Riemannian manifold $(N,h)$. Each fiber canonically inherits an
$H$@-invariant metric from the given metric on~$H$.
\roster
\item"(i)" For $Z\in\h$ we denote the corresponding vector field
$p\mapsto\tdiffzero{t}p\cdot\exp(tZ)$ on~$M$ by~$Z$ again.
Note that $Z$ is $H$@-invariant since $H$ is abelian.
\item"(ii)" A connection form $\om$ on~$M$ is an $\h$@-valued, $H$@-invariant
  $1$@-form on~$M$ such that $\om(Z)=Z$ for all $Z\in\h$.
  For any connection form $\om$ on~$M$ and any
  $Z\in\h$, we define the $1$@-form $\om_Z:=\<\om(\,.\,),Z\>$ on~$M$, where
  $\scp$ denotes the scalar product induced on~$\h$ by the metric on~$H$.
\item"(iii)" For any connection form $\om$ on~$M$,
  we denote by $g\subom$ the unique $H$@-invariant Riemannian metric on~$M$
  which satisfies:
\newline
\qquad\qquad 1.) $g\subom$ induces the given invariant metric on each fiber.
\newline
\qquad\qquad 2.) The projection $\pi_H:(M,g\subom)\to(N,h)$ is a
  Riemannian submersion.
\newline
\qquad\qquad 3.) The $\om$@-horizontal distribution $\kernn\om$
  is $g\subom$@-orthogonal to the fibers.
\newline
  Note that in particular, $g\subom(X,Z)=\<\om(X),Z\>=\om_Z(X)$ for
  all $X\in TM$ and $Z\in\h$.
\item"(iv)" If $M$ is the trivial bundle $N\times H$ and $\la$ is an
  $\h$@-valued $1$@-form on~$N$ then we write $g\subla:=g\subom$\,,
  where~$\om$ is the unique connection form extending $\la\restr{TN\times
  \{0\}}$\,; that is, $\om(X,Z)=\la(X)+Z$ for all $(X,Z)\in
  T_{(p,z)}(N\times H)\cong T_pN\times\h$.
  In particular, $g\subla$ has the properties 1.),~2.)~from (iii)~above, and
  the vector $(X,-\la(X))$ is horizontal for all $X\in TN$.
  For each $Z\in\h$ we define the $1$@-form $\la_Z:=\<\la(\,.\,),Z\>$ on~$N$.
\item"(v)" If $F:M\to M$ is a gauge transformation, that is, a bundle
  automorphism which induces the identity on~$N$, then it is obvious from
  the definitions that $F$ is an isometry from $(M,g_{F^*\om})$ to
  $(M,g\subom)$ for any connection form~$\om$ on~$M$.
  If $\alpha$ and~$\om$ are two connection forms such that
  $\alpha=\om+df$ for some $f\in C^\infty(M,\h)$, then the gauge
  transformation $F:p\mapsto p\cdot\exp(f(p))$ satisfies $\alpha=F^*\om$;
  thus $(M,g_\alpha)$ and $(M,g\subom)$ are isometric.
  Analogously, if two $\h$@-valued
  $1$@-forms $\la,\mu$ on~$N$ differ by~$df$ for some $f\in C^\infty(N,\h)$,
  then the associated metrics $g\subla$ and~$g_\mu$ on $N\times H$ are
  isometric.
\endroster

\bigskip

\proclaim{Theorem 1.6}
Let $(N,h)$ be a closed Riemannian manifold and $H$ be a torus equipped
with an invariant metric. Let $M$ be a principal $H$@-bundle over
$(N,h)$, and let $\om,\omp$ be two connection forms on~$M$. Assume:
\roster
\item"($*1$)" For every $Z\in\h$ there exists a bundle automorphism
  $F_Z:M\to M$ which induces an isometry on the base manifold $(N,h)$
  and satisfies $\omp_Z=F_Z^*\om_Z$\,.
\endroster
Then $(M,g\subom)$ and $(M,g\subomp)$ are isospectral.
\endproclaim

\smallskip

\demo{Proof}
We show that $(M,g\subom)$ and $(M,g\subomp)$ satisfy the hypotheses
of Theorem~1.3. Denote by $\nabla\supom$ the Levi-Civit\`a connection
of~$g\subom$\,.
Let $X\in TM$ be arbitrary and extend it to an $H$@-invariant vector field
on~$M$. Then $\om(X):M\to\h$ is $H$@-invariant; moreover, for each $Z\in\h$
the flows of~$X$ and~$Z$ commute, whence $[X,Z]=0$. Thus
$$g\subom(\nabla\supom_Z Z,X)=Z(g\subom(X,Z))-g\subom(Z,\nabla\supom_Z X)
=Z\<\om(X),Z\>-g\subom(Z,\nabla\supom_X Z)=0-0=0.
$$
Since $X\in TM$ was arbitrary, we conclude $\nabla\supom_Z Z=0$ for each
$Z\in\h$. Therefore the
$H$@-orbits are totally geodesic in~$(M,g\subom)$. The first condition
of Theorem~1.3 is thus satisfied for~$g\subom$\,, and similarly
for~$g\subomp$\,.

It remains to check condition~(ii) of Theorem~1.3.
For $W=H$, there is nothing to show since $(M/H,g\subom^H)$ and
$(M/H,g\subomp^H)$ both equal $(N,h)$ by the choice of $g\subom$
and~$g\subomp$\,.
Let $W$ be a closed subgroup of codimension~$1$ in~$H$. Choose $Z\perp T_eW$
in $\h\minzero$, and let $F_Z:M\to M$ be a bundle
automorphism as in~($*1$). We claim that $F_Z$ induces an isometry
from $(M/W,g\subomp^W)$ to $(M/W,g\subom^W)$, where
$g\subom^W$ and $g\subomp^W$ are the submersion metrics induced
by~$g\subom$ and~$g\subomp$\,.
Since $F_Z$ commutes with the $H$@-action and induces an isometry
on the base manifold $(N,h)$, we only need to check that for any
$\omp$@-horizontal vector~$X$, the vector $F_{Z*}X$ is $\om$@-horizontal
up to an error tangent to the $W$@-orbits; in other words, $F_{Z*}X$
is $g\subom$@-orthogonal to~$Z$. But $g\subom(F_{Z*}X,Z)
=\<\om(F_{Z*}X),Z\> = \om_Z(F_{Z*}X)=\omp_Z(X)=\<\omp(X),Z\>=0$
since $\omp_Z=F_Z^*\om_Z$ and $\omp(X)=0$.
\qed\enddemo

\smallskip

\subheading{Remarks 1.7}

(i) For any closed subgroup~$W$ of codimension~$1$ in~$H$, let $\CC_W$
denote the space of smooth functions on~$M$ which are invariant under~$W$.
Recall that $\pi_W^*:C^\infty(M/W)\to
\CC_W$ intertwines the Laplacians associated with $g\subom$ and $g\subom^W$
(or $g\subomp$ and $g\subomp^W$) because, as we have seen
in the proof of Theorem~1.6, the fibers of 
$\pi_W$ are totally geodesic for $g\subom$ and $g\subomp$\,.
We also saw in the proof of Theorem~1.6 that~$F_Z$ induces
an isometry from $(M/W,g\subomp^W)$ to $(M/W,g\subom^W)$, where
$Z\in\h\minzero$ is orthogonal to $T_eW$, and $F_Z$ is chosen as
in~\thetag{$*1$}. The pullback of this
isometry intertwines the corresponding Laplacians on $M/W$.
Combining these intertwining maps, we conclude that
$${\Delta_{g\subomp}\vphantom)}\restr{\CC_W} = (F_Z^*\circ\Delta_{g\subom}
  \circ{F_Z^*}\inv)\restr{\CC_W}\,.
$$
This last fact can of course also be derived directly from the assumptions
on~$F_Z$\,, without even introducing the quotient manifolds $M/W$.
We do not present this alternative argument here because we will
do so later in a more general situation, namely, in the proof of
Theorem~4.3 in Chapter~4 (see also Remark~4.4(i)) which is a
generalization of Theorem~1.6.

(ii) If in the context of Theorem~1.6 there exists a bundle automorphism
$F:M\to M$ which satisfies~\thetag{$*1$} for each $Z\in\h$, then we have
$\omp=F^*\om$; hence $F:(M,g\subomp)\to(M,g\subom)$ is an isometry.
In order to obtain nontrivial pairs of isospectral manifolds from
Theorem~1.6 it is thus crucial that \thetag{$*1$} be satisfied without there
being a choice of the~$F_Z$ independent of~$Z$.
Such examples do exist; see Section~1.3 below and Chapters~2 and~3.

\bigskip

In the following proposition we specialize Theorem~1.6 to the case of
products. We use Notation~1.5(iv).

\bigskip

\proclaim{Proposition 1.8}
Let $(N,h)$ be a closed Riemannian manifold and $H$ be a torus with
Lie algebra $\h:=T_eH$, equipped with an invariant metric.
Let $\la,\lap$ be two $\h$@-valued 1@-forms on~$N$ which satisfy:
$$
\lap_Z\in\Isom(N,h)^*(\la_Z)\text{ for each }Z\in\h.\tag{$*2$}
$$
Then $(N\times H,g\subla)$ and $(N\times H,g\sublap)$ are isospectral.
\endproclaim

\smallskip

\demo{Proof}
Let $\om,\omp$ be the connection forms on the trivial $H$@-bundle
$N\times H$ which extend $\la$ and~$\lap$, respectively; thus
$g\subla=g\subom$ and $g\sublap=g\subomp$\,. We check that~$\om$ and~$\omp$
satisfy condition~($*1$) of Theorem~1.6. Let $Z\in\h$. By~($*2$),
there exists $f_Z\in\Isom(N,h)$ such that $\lap_Z=f_Z^*(\la_Z)$.
Define $F_Z:=(f_Z\,,\Id):N\times H\to N\times H$. Then $F_Z$ is obviously
a bundle isomorphism which induces the isometry~$f_Z$ on $(N,h)$ and
by the definition of $\om,\omp$ satisfies
$\omp_Z=F_Z^*(\om_Z)$.\qed\enddemo

\smallskip

\subheading{Remark 1.9}
In analogy with Remark~1.7(ii) we observe that if $\lap\in\Isom(N,h)^*
\la$, then $(N\times H,g\subla)$ and $(N\times H,g\sublap)$ are isometric.
But there do exist examples where this is not the case although \thetag{$*2$}
is satisfied; see Example~1.11 in the following 
section, and various new examples in the Chapters~2 and~3.

\bigskip

\smallskip

\subhead 1.3\ \ Review of some previously known examples\endsubhead
\bigskip

\noindent
In this section we will explain some previously known examples of
closed isospectral, locally non-isometric manifolds from the point of
view of Theorem~1.6 and Proposition~1.8.

The first class of examples (Example~1.11) concerns manifolds diffeomorphic
to $S^{m-1}\times T^r$ with $m\ge5$ and $r\ge2$.
Continuous isospectral families of locally non-isometric metrics on such
manifolds were constructed in \cite{22}; independently, Z.~Szab\'o
\cite{44} found pairs of isospectral, locally non-isometric metrics on
$S^{4k-1}\times T^3$ with $k\ge2$.
These manifolds arise as the boundaries of certain isospectral, locally
non-isometric manifolds with boundary which were constructed
in~\cite{27} and~\cite{44}, respectively.

The second class of examples (Example~1.14) concerns manifolds diffeomorphic
to $S^{m-1}\times S$, where $m\ge5$ and $S$ is a compact Lie group
of rank at least two. The author constructed in~\cite{42} continuous
families of isospectral, locally non-isometric metrics on these manifolds,
thereby providing, in particular, the first examples of simply
connected isospectral manifolds. Otherwise, the assumption that $S$ is simply
connected which was made in~\cite{42} plays no role in the construction of
the isospectral metrics and in the non-isometry arguments.
Thus the first class of examples, mentioned above, can actually be viewed
as a subclass of the second one.

As we will see, both classes of examples arise from Theorem~1.6;
the first class arises even from the more special Proposition~1.8.

We recall that together with the pairs of two-step nilmanifolds
constructed in \cite{21} and \cite{27} (see Remark~1.4(v)), and
some recent examples by C.~Gordon and Z.~Szab\'o \cite{23} (see Remark~4.4
in Chapter~4), these manifolds provide all previously known examples
of isospectral, locally non-isometric, closed manifolds.

We start by a definition introduced by C.~Gordon and E.~Wilson.

\bigskip

\definition{Definition 1.10 \cite{27}}
Two linear maps $j,j':\R^r\to\so(m)$ are called {\it isospectral},
denoted $j\sim j'$, if for every $Z\in\R^r$ there exists $A_Z\in\O(m)$
such that $j'_Z=A_Z j_Z A_Z\inv$.
\enddefinition

\smallskip

\subheading{Example 1.11:\ \ Products of spheres with tori
  (\cite{22}, \cite{44})}

\noindent
Let $S^{m-1}\subset\R^m$ be the
$(m-1)$@-dimensional unit sphere.
Let $H$ be a torus with a fixed invariant metric and Lie algebra
$\h\cong\R^r$.

For each linear map $j:\h\to\so(m)$ we define an $\h$@-valued $1$@-form~$\la$
on~$S^{m-1}$ by requiring that
$$\la_Z(X)=-\tsize\frac12 \<j_Z p,X\>\tag2
$$
for each $X\in T_pS^{m-1}$ and $Z\in\h$, where $\scp$ denotes the standard
scalar product on~$\R^m$.
Here, $j_Z\in\so(m)$ acts on $p\in S^{m-1}\subset\R^m$ by usual
multiplication.
(Why the factor $-1/2$ is convenient will become clear at the end of the
following Remark~1.12.)
         
Now let $j,j':\h\to\so(m)$ be such that $j\sim j'$.
Let $\la,\lap$ be the associated $\h$@-valued $1$@-forms on~$S^{m-1}$.
For each $Z\in\h$ choose $A_Z\in\O(m)$ such that $j'_Z=A_Z j_Z A_Z\inv$.
Then it follows immediately from~\thetag{2} that $\lap_Z=A_Z^{-1*}\la_Z$\,.
Let $h$ be the round standard metric on~$N:=S^{m-1}$.
Since $A_Z\inv$ is an isometry of $(N,h)$, condition~\thetag{$*2$}
from Proposition~1.8 is satisfied for $\la$ and $\lap$.
We conclude that $(S^{m-1}\times H,g\subla)$
and $(S^{m-1}\times H,g\sublap)$ are isospectral, where the metrics
$g\subla$ and $g\sublap$ are associated with $\la$ and~$h$, resp\. 
with $\lap$ and~$h$, as in Notation~1.5(iv).

\bigskip

\subheading{Remark 1.12}

\noindent
We explain why the isospectral manifolds from Example~1.11
are exactly those constructed in \cite{22}.
In that paper the approach was as follows (up to minor changes of notation).

Let $\vv:=\R^m$ and $\h:=\R^r$ be endowed with the euclidean standard
metrics, and let $\LL\subset\h$ be a lattice of full rank.
For each linear map $j:\h\to\so(\vv)$ consider the two-step nilpotent
Lie algebra $\g_j:=\vv\oplus\h$ whose Lie bracket is defined by requiring
that $\h$ be central, $[\g_j\,,\g_j]\subseteq\h$, and
$\<[X,Y],Z\>=\<j_Z X,Y\>$ for all $X,Y\in\vv$ and $Z\in\h$.
Let $G_j$ be the associated simply connected Lie group,
and let $g_j$ be the left invariant metric on~$G_j$
which corresponds to the standard scalar product on $\g_j=\vv\oplus\h$.
The group exponential map $\exp:\g_j\to G_j$ is a diffeomorphism
which restricts to a linear isomorphism between $\h$ and
$\exp\,\h\subset G_j$\,.
Denote by $\bar G_j$ the quotient of~$G_j$ by the discrete central subgroup
$\exp\,\LL\subset\exp\,\h$, and denote by~$\bar g_j$ the left invariant
metric on~$\bar G_j$ induced by~$g_j$\,.
The group exponential map $\overline\exp:\g_j\to\bar G_j$ induces a
diffeomorphism from $\vv\times(\h/\LL)$ to~$\bar G_j$\,.
Define $M_j\subset\bar G_j$ as the image of $S^1(\vv)\times(\h/\LL)$
under this diffeomorphism,
and denote the induced metric on~$M_j$ by~$\bar g_j$ again.

In \cite{22} it was proven that for $j\sim j'$, the Riemannian
manifolds $(M_j\,,\bar g_j)$ and $(M_{j'}\,,\bar g_{j'})$ (denoted there by
$N(j)$ and $N(j')$) are isospectral. This is exactly the same as what
we stated in Example~1.11 above: As we are going to see now,
$(M_j\,,\bar g_j)$ is isometric to our above $(S^{m-1}\times H,g\subla)$,
where we let $H:=\h/\LL$ and where $\la$ is associated with~$j$ as
in~\thetag{2}. More precisely, under the identification of $S^{m-1}$ with
$S^1(\vv)$, we claim that
$\overline\exp:\vv\oplus\h\to\bar G_j$ induces an isometry
from $(S^1(\vv)\times(\h/\LL),g_\la)$ to $(M_j\,\bar g_j)$.

We first consider the metric $\expst g_j$ on $\vv\oplus\h$.
Extend~$\la$ to an $\h$@-valued $1$@-form on~$\vv$ by letting
$\<\la(X),Z\>=-\frac12\<j_ZV,X\>$ for each $V\in\vv$ and
$X\in T_V\vv$.
By the Campbell-Baker-Hausdorff formula and the definition of $\la$ and
$[\,,\,]$ (both associated with~$j$) we have
$$\exp_*\restr{V+W}(X+Z)\;=\;L_{\exp(V+W)*}(X+Z-\tsize\frac12[V,X])
  \;=\;L_{\exp(V+W)*}(X+\la(X)+Z)
$$
for all $V,X\in\vv$ and $W,Z\in\h$. Thus
$$(\expst g_j)\restr{V+W}(X_1+Z_1\,,\,X_2+Z_2)
  \;=\;\<X_1+\la(X_1)+Z_1\,,\,X_2+\la(X_2)+Z_2\>
$$
for all $V,X_1\,,X_2\in\vv$ and $W,Z_1\,,Z_2\in\h$.
This shows that $\expst g_j$ induces the given euclidean metric on the
$\h$@-fibers, that $X-\la(X)$ is orthogonal to $\h\subset
T_{V+W}(\vv\oplus\h)$,
and that $(\expst g_j)(X_1-\la(X_1),X_2-\la(X_2))=\<X_1\,,X_2\>$
for all $X,X_1\,,X_2\in\vv\subset T_{V+W}(\vv\oplus\h)$\,.

Note that $\overline{\text{$\exp$}}^*\bar g_j=\expst g_j$\,. Thus the above
properties hold also for $\overline{\text{$\exp$}}^*\bar g_j$\,, and in
particular for its restriction to $S^1(\vv)\times\h$.
By the definition of~$g_\la$ this implies that $\overline\exp$
induces indeed an isometry from $(S^1(\vv)\times(\h/\LL),g_\la)$
to $(\overline\exp(S^1(\vv)\times\h),\bar g_j)=(M_j\,,\bar g_j)$.

\bigskip

\subheading{Remarks 1.13}

(i) Given $H=\h/\LL$ as above, we say that two linear maps~$j,j':\h\to
\so(m)$ are {\it trivially isospectral\/} if there
exist $A\in\O(m)$ and $C\in\O(\h)$ such that $C(\LL)=\LL$ and
$j'_{CZ}=Aj_ZA\inv$ for all $Z\in\h$. (Note that in contrast to the
isospectrality condition in Definition~1.10, the map~$A$ is assumed
to be independent of~$Z\in\h$.)
If this is the case, then the corresponding metrics~$g\subla$
and~$g\sublap$ on $S^{m-1}\times H$ are obviously isometric;
an isometry from~$g\subla$ to~$g\sublap$ is induced by $(A,C)$.

For $m\le4$ and $\dimm\,\h\le2$, and also for $m\le3$ and arbitrary
dimension of~$\h$, elementary arguments show that isospectrality
of two linear maps $j,j':\h\to\so(m)$ always implies triviality
in the above sense. Nontrivial isospectral manifolds
of the type described in Example~1.11 can therefore occur only
in the case $m+\dimm\,\h\ge7$, that is, if $\dimm(S^{m-1}\times H)\ge6$.
Dimension six is indeed attained here since for $m=5$ and $\dimm\,\h=2$
there do exist nontrivial families of isospectral metrics; see~(ii) below
and also Proposition~1.16.

(ii) It was proven in \cite{22}, using a result from \cite{27}, that
for each $m\ge5$ the above method provides continuous
$d$@-parameter families of isospectral metrics on $S^{m-1}\times T^2$,
where $d$ is of order at least~$O(m^2)$. The proof that the manifolds in
these multiparameter-families have pairwise different local geometries is
rather abstract in the sense that it does not distinguish the metrics in
terms of straightforwardly formulated geometrical quantities.

However, it was also shown in \cite{22} that for some of the isospectral
families the maximum of the scalar curvature varies during the deformation.
(In particular, the manifolds are not locally isometric.)
For a more detailed statement concerning the critical values of the scalar
curvature see Proposition~1.16 below
which refers to the following Example~1.14.
Note that the above Example~1.11 can be considered as a special case
of Example~1.14; thus all statements of Proposition~1.16 hold also for the
isospectral manifolds of Example~1.11.

\bigskip

\subheading{Example 1.14:\ \ Products of spheres with compact Lie groups
  (\cite{42})}

\noindent
Let $S$ be a compact Lie group containing a torus~$H$ of dimension at
least two, and let $\h:=T_eH\subseteq T_eS$.
Let $k$ be a bi-invariant metric on~$S$, and
let $H$ be equipped with the metric induced by~$k$.
Let $S^{m-1}\subset\R^m$ be the $(m-1)$@-dimensional unit sphere.
View $S^{m-1}\times S$ as a principal $H$@-bundle with respect to the
left action of~$H$ on the second factor.

For each linear map $j:\h\to\so(m)$ we define a connection form~$\om$
on the $H$@-bundle $S^{m-1}\times S$ by requiring (using Notation~1.5(ii))
that
$$\om_Z(X,U)=-\tsize\frac12\<j_Zp,X\>+k(U,R_{s*}Z)\tag3
$$
for all $(X,U)\in T_{(p,s)}(S^{m-1}\times S)$ and $Z\in\h$, where
$R_s$ denotes right multiplication by~$s$ and $\scp$ denotes the
standard scalar product on~$\R^m$.
In order to check that $\om$ is indeed a connection form, note that
the definition implies (using Notation~1.5(i)) that
$\om(Z_{(p,s)})=\om(0,R_{s*}Z)=Z$ for all $Z\in\h$; furthermore, $\om$~is
$H$@-(left) invariant because for all $(X,U)\in T_{(p,s)}(S^{m-1}\times S)$
and all $Z\in\h$ and $z\in H$ we have
$$\align\om_Z(X,L_{z*}U)&=-\tsize\frac12\<j_Zp,X\>+k(L_{z*}U,R_{zs*}Z)
  =-\tsize\frac12\<j_Zp,X\>+k(U,R_{s*}\Ad_z\inv Z)\\
  &=\,\om_Z(X,U)\endalign
$$
by the bi-invariance of~$k$ and the commutativity of~$H$.

Now let $j,j':\h\to\so(m)$ be such that $j\sim j'$, and let
$\om,\omp$ be the associated connection forms on $S^{m-1}\times S$.
For each $Z\in\h$ choose $A_Z\in\O(m)$ such that
$j'_Z=A_Z j_Z A_Z\inv$. Define $F_Z:=(A_Z\inv,\Id):S^{m-1}\times S
\to S^{m-1}\times S$.
Obviously $F_Z$ is a bundle automorphism. Moreover, it follows
immediately from~\thetag{3} that $\omp_Z=F_Z^*\om_Z$.
Finally, define the metric $h$ on the base manifold
$N:=S^{m-1}\times(H\backslash S)$
as the product of the round standard metric on~$S^{m-1}$ and the
submersion metric~$\bar k$ on $H\backslash S$ induced by~$k$.
Then $F_Z$ induces the $h$@-isometry $\bar F_Z:=(A_Z\inv,\Id)$ on~$N$.
Thus $\om$ and~$\omp$ satisfy
condition~\thetag{$*1$} of Theorem~1.6. We conclude
that $(S^{m-1}\times S,g\subom)$ and $(S^{m-1}\times S,g\subomp)$
are isospectral, where the metrics $g\subom$ and~$g\subomp$
are associated with $\om$ and~$h$, resp\. with $\omp$ and~$h$,
as in Notation~1.5(iii).

\bigskip

\subheading{Remark 1.15}

\noindent
We explain why the isospectral manifolds from Example~1.14
are exactly those constructed in~\cite{42}.
Actually, $S$~was  assumed to be simply connected there,
but this did not play a role in any of the arguments, except
for the fact that it caused the manifolds $S^{m-1}\times S$ to be simply
connected.

Let $S$, $H$, $\h$, and $k$ be as in Example~1.14.
Let $r:=\dimm\,\h$, and let $\vv:=\R^m$ be equipped with the euclidean
standard metric.
In~\cite{42} we associated with each linear map $j:\h\to\so(\vv)$
a metric~$g_j$ on $\vv\times S$. Instead of reviewing its original
definition here, we recall a characterization of~$g_j$ resulting from
Lemma~1.7(i) in~\cite{42}.
Namely, a $g_j$@-orthogonal basis at $(V,s)\in\vv\times S$ is given by
$$\{(X_a\,,R_{s*}(\tsize\frac12[V,X_a]))\mid a=1,\,\ldots,m\}
  \;\cup\;\{(0,R_{s*}U_i)\mid i=1,\,\ldots,\dimm\,S\},\tag4
$$
where $\{X_1\,,\,\ldots,X_m\}$ is an orthonormal basis of~$\vv$,
$\;\{U_1\,,\,\ldots,U_{\dimm\,S}\}$ is a $k$@-ortho\-normal basis
of $T_eS$, and $[\,\,,\,]$ is the Lie bracket on $\g_j:=\vv\oplus\h$
associated with~$j$ as in Remark~1.12;
that is, $\h$ is central, $[\g_j\,,\g_j]\subseteq\h$,
and $\<[X,Y],Z\>=\<j_ZX,Y\>$ for all $X,Y\in\vv$ and $Z\in\h$.
Denote the restriction of~$g_j$ to $S^1(\vv)\times S\subset\vv\times S$
by~$g_j$ again.

In~\cite{42} we showed that for $j\sim j'$ the Riemannian manifolds
$(S^1(\vv)\times S,g_j)$ and $(S^1(\vv)\times S,g_{j'})$
are isospectral.
This is exactly the same as what was stated in Example~1.14 above:
As we are going to see now, the metric~$g_j$ on $S^1(\vv)\times S
=S^{m-1}\times S$ is equal to our above metric~$g\subom$ on $S^{m-1}
\times S$, where $\om$ is associated with~$j$ as in~\thetag{3}.

In fact, it follows from the description of~$g_j$ on $\vv\times S$ by the
orthonormal bases given in~\thetag{4} that
\roster
\item"1.)" The metric induced by~$g_j$ on the left $H$@-orbits in
  $S^{m-1}\times S$ is the one inherited from the metric $k\restr H$
  on~$H$.
\item"2.)" The projection from $(S^{m-1}\times S,g_j)$ to
  $(S^{m-1}\times(H\backslash S),h)$ is a Riemannian submersion,
  where $h$ is as in Example~1.14.
\item"3.)" The vector $(X,U)\in T_{(p,s)}(S^{m-1}\times S)$ is
  $g_j$@-orthogonal to the fiber $\{p\}\times Hs$ if and only if
  for each $Z\in\h$ we have
$$\align
\qquad 0&=g_j\bigl((X,U),(0,R_{s*}Z)\bigr)\\
   &=g_j\bigl((X,R_{s*}(\tsize\frac12[p,X])),(0,R_{s*}Z)\bigr)
    -k(R_{s*}(\tsize\frac12[p,X]),R_{s*}Z)
    +k(U,R_{s*}Z)\\
   &=0-\tsize\frac12\<j_Zp,X\>+k(U,R_{s*}Z)
   =\om_Z(X,U).\endalign
$$
\endroster
{}From these properties and the definition of~$g_\om$ (recall
Notation~1.5(iii)) it follows that $g_j=g_\om$\,.

\bigskip

Concerning non-isometry criteria for the above manifolds
$(S^{m-1}\times S,g_j)$, we showed in~\cite{42} the following result,
parts of which we will need again in Chapter~3.

\bigskip

\proclaim{Proposition 1.16 (\cite{42})}

\noindent
In the context of Example~{\rm1.14} and Remark~{\rm1.15}
suppose $m\ge5$ and $\dimm\,H=2$, and let $\{Z_1\,,Z_2\}$ be an
orthonormal basis of~$\h\subseteq T_eS$.
\roster
\item"(i)"
Let $j,j'$ be two linear maps from~$\h$ to~$\so(m)$ such that
$j\sim j'$.\newline
\qquad\qquad{\rm(a)} If $j_{Z_1}^{\phantom\prime2}+j_{Z_2}^{\phantom\prime2}$
and $j_{Z_1}^{\prime\,2}+j_{Z_2}^{\prime\,2}$ have different sets
of eigenvalues, then the scalar curvature of~$g_j$ and the scalar curvature
of~$g_{j'}$ on $S^{m-1}\times S$ have different sets of critical values
(\cite{42}, Proposition~3.5).\newline
\qquad\qquad{\rm(b)} If $\|j_{Z_1}^{\phantom\prime2}
+j_{Z_2}^{\phantom\prime2}\|^2\ne\|j_{Z_1}^{\prime\,2}
+j_{Z_2}^{\prime\,2}\|^2$ {\rm(}where $\|\,.\,\|$ denotes the standard
euclidean norm on real $m\times m\,$@-matrices{\rm)}, then the heat
invariants for the Laplacian on $1$@-forms are not equal for $g_j$
and~$g_{j'}$ (\cite{42}, Proposition~5.1).
\item"(ii)"
There exists a Zariski open subset~$\Cal U$ of the space~$\J$
of all linear maps from~$\h$ to~$\so(m)$ with the property that every
$j\in\Cal U$ is contained in a smooth isospectral family
$j(t)$, defined in some open interval around zero, such that $j(0)=j$ and
such that $\|j_{Z_1}(t)^2
+j_{Z_2}(t)^2\|^2$ has nonzero derivative at $t=0$
(\cite{42}, Proposition~4.1).
\item"(iii)"
{}From {\rm(i)} it follows that the isospectral manifolds
$(S^{m-1}\times S,g_{j(t)})$ from~{\rm(ii)} are not isospectral
on $1$@-forms and have different sets of critical values for the scalar
curvature. In particular, the manifolds are not locally isometric.
\item"(iv)"
For $m=5$, an explicit example of an isospectral family~$j(t)$ with the
property that
$\|j_{Z_1}(t)^2+j_{Z_2}(t)^2\|^2\ne
\text{\rm const}$ is given by
$$j_{Z_1}(t):=\left(\smallmatrix 0&0&-t&0&0\\0&0&0&t-1&0\\t&0&0&0&-\ph(t)\\
0&1-t&0&0&-\psi(t)\\0&0&\ph(t)&\psi(t)&0\endsmallmatrix\right),\qquad
j_{Z_2}(t):=\left(\smallmatrix 0&1&0&0&0\\-1&0&0&0&0\\0&0&0&1&0\\0&0&-1&0&0\\
0&0&0&0&0\endsmallmatrix\right),
$$
where $\ph(t)=\bigl((t^4-3t^2+1)/(1-2t)\bigr)^{1/2}$ and
$\psi(t)=\bigl((-t^4+4t^3-3t^2-2t+1)/(1-2t)\bigr)^{1/2}$
{\rm(}see \cite{42}, Remark~4.3(ii){\rm)}.
This family is defined for $t\in\bigl[\frac12(1-\sqrt5),
\frac12(3-\sqrt5)\bigr]$. The $j(t)$ are pairwise isospectral since
$\det\bigl(\la\Id-(sj_{Z_1}(t)+uj_{Z_2}(t))\bigr)=\la^5+(3s^2+2u^2)\la^3+
(s^2+u^2)^2\la$
is independent of~$t$; however,
$\|j_{Z_1}(t)^2+j_{Z_2}(t)^2\|^2
=4t^2-4t+26$ is nonconstant in~$t$.
\endroster\endproclaim

\bigskip

\bigskip
\noindent
{\gross 2.\hbox to6pt{}Four-dimensional examples} 

\bigskip

\noindent
Recall from Remark~1.13(i) that nontrivial isospectral pairs of linear
maps $j:\R^r\to\so(m)$ can occur only if $m+r\ge7$; therefore, all
nontrivial examples of isospectral manifolds arising from Example~1.11 are of
dimension at least six. The manifolds from Example~1.14 have at least as
many dimensions as those from Example~1.11. Finally, the locally
non-isometric,
two-step nilmanifolds from \cite{21} and \cite{27}, too, are based on
pairs of isospectral maps~$j$ as above; in particular, they are of dimension
at least seven. Thus the smallest dimension for which examples of
isospectral, locally non-isometric manifolds have been known is six.

In this chapter we will present pairs of isospectral, locally non-isometric
metrics in dimension four; more precisely, on $S^2\times T^2$.
The idea is the following: Let~$h$ be the standard metric on~$S^2$.
Try to construct nontrivial pairs of $\R^2$@-valued $1$@-forms $\la,\lap$
on~$S^2$ such that $\la_Z'\in\Isom(S^2,h)^*\la_Z=\O(3)^*\la_Z$ for all
$Z\in\R^2$. ``Nontrivial'' means that $\lap\notin\O(3)^*\la$ (or more
strongly, in view of Remark~2.1(i) below, $\lap\notin\O(3)^*(\la+df)
$ for all $f\in C^\infty(S^2,\R^2)$\,).
Then apply Proposition~1.8 to obtain isospectral metrics on $S^2\times T^2$.

Since we know that this would not be possible here by the approach of
Example~1.11\,/\,Remark~1.12, we must drop an unnecessary assumption
which was made there.
Remember that in Example~1.11 each $\la_Z$ was of the form $X\mapsto
-\frac12\<j_Zp,X\>$ for $X\in T_pM$; in particular, $\la$ depended {\it
linearly on the coordinates of the basepoint\/}~$p\in S^{m-1}\subset\R^m$
of~$X$.
This plays no role for the applicability of Proposition~1.8. The reason
why this linearity was present throughout the earlier examples (1.11 and
1.14) lies in the way they were originally discovered, namely, as companions
to certain two-step nilmanifolds (recall Remark~1.12; for a different
point of view which naturally includes the mentioned linearity see
Remark~3.4).

The key of our construction of examples in dimension four is to replace the
linear dependence of~$\la$ on the basepoint by quadratic dependence.
This is what we do in Section~2.1, obtaining explicit nontrivial examples.
In Section~2.2 we derive a formula for the associated scalar curvature
and show
for a specific isospectral pair $g\subla\,,g\sublap$ given in Section~2.1
that the preimages in $S^2\times T^2$ of the maximal scalar curvature
of~$g\subla$ and~$g\sublap$ have different dimensions. In particular,
$(S^2\times T^2,g\subla)$ and $(S^2\times T^2,g\sublap)$ are not locally
isometric.
Finally, we present a more general result concerning the scalar curvature:
Let $(N\times H,g\subla)$, $(N\times H,g\sublap)$ be any pair of isospectral
manifolds arising from Proposition~1.8, and assume that the base manifold~$N$
is two-dimensional. In this situation we prove that
the associated scalar curvature functions have the same range,
and for each~$k\in\N$ the integrals of their $k$@-th
powers coincide (Theorem~2.11).
If in addition the torus~$H$ is also two-dimensional,
then the same holds for the $L^2$@-norms of the Ricci and curvature tensors.
This contrasts with the
higher-dimensional examples $(S^{m-1}\times T^2,g\subla)$ from
Example~1.11, where the isospectral manifolds did in general neither share
the range of the scalar curvature nor $\int\scal^2$ or $\int\|\Ric\|^2$
(compare Remark~1.13(ii), Proposition~1.16, and \cite{42}, Lemma~5.4).

\bigskip

\smallskip

\subhead 2.1\ \ Isospectral, locally non-isometric metrics on $S^2\times T^2$
  \endsubhead
\bigskip

\noindent
Before starting the search for suitable pairs of $\R^2$@-valued
$1$@-forms~$\la,\lap$ on~$S^2$ in the ``quadratic'' category (i.e.,
with $\la,\lap$ depending quadratically on the basepoint),
we make an observation which implies that it would be useless
to instead replace 
the skew-symmetric linear maps~$j_Z$ of Example~1.11 by more general
linear maps or by constant ones;
moreover, improving the dimension of the manifold to three instead of four
is not possible by our methods.

\bigskip

\subheading{Remarks~2.1}

(i) Let $H$ be a torus with an invariant metric and Lie algebra $\h:=T_eH$,
and let $\{Z_1\,,\,\ldots,Z_r\}$ be an orthonormal basis of~$\h$.
Let $h$ be the standard metric on~$S^{m-1}$, and let $\mu,\nu$ be two
$\h$@-valued $1$@-forms on~$S^{m-1}$.
For $\la:=\mu-\nu$ suppose that each $\la_Z:=\<\la(\,.\,),Z\>$ is of
the type $X\mapsto\<w_Z(p),X\>$ for $X\in T_pS^{m-1}$, where $w_Z$ is
a vector field on~$\R^m$ which has a potential function~$\phi_Z$\,.
For example, this is the case if $w_Z(p)$ is constant, or if $w_Z$ is
linear and symmetric.
Then $\la=df$ with $f(p):=\sum_{i=1}^r\phi_{Z_i}(p)Z_i\in\h$.
Thus $\mu=\nu+df$, and by~1.5(v) the associated
metrics~$g_\mu$ and~$g_\nu$ on $S^{m-1}\times H$, defined as in
Notation~1.5(iv), are isometric.

(ii) If $\la,\lap$ are two $\h$@-valued $1$@-forms on~$S^1$ then each of
them is the sum of an $S^1$@-invariant form and an exact one. By~1.5(v)
we can thus assume, without changing the isometry classes of $g\subla$ and
$g\sublap$ on $S^1\times H$, that~$\la$ and~$\lap$ are $S^1$@-invariant.
But if two such forms satisfy condition~\thetag{$*2$} of Proposition~1.8,
then it is easy to see that $\la=\pm\lap$ and thus the associated metrics
on $S^1\times H$ are isometric. This shows that in order to get nontrivial
pairs of isospectral metrics on $S^{m-1}\times H$ by using Proposition~1.8
we cannot lower the dimension of $S^{m-1}$ to less than two. Since the
dimension of~$H$ must also be at least two (recall Remark~1.4(i)), examples
in dimension four (which we obtain in this chapter) is the best we can
hope for.

(iii) More generally, if in the context of Theorem~1.3 the base manifold
$M/H$ is one-dimensional, then any $H$@-invariant metric~$g$ for which
the $H$@-orbits are totally geodesic must be flat. Thus any isospectral
pair of metrics~$g,g'$ arising from Theorem~1.3 in this situation
would be locally isometric anyway.

\bigskip

\subheading{Notation and Remarks 2.2}
\roster
\item"(i)"
Let $\scp$ denote the standard scalar product on~$\R^3$,
and let $\ourspace$ be endowed with the canonically induced scalar product;
that is, $\<(X_1Y_1)^*\otimes V_1\,,(X_2Y_2)^*\otimes V_2\>
=\frac12(\<X_1\,,X_2\>\<Y_1\,,Y_2\>+\<X_1\,,Y_2\>\<X_2\,,Y_1\>)
\<V_1\,,V_2\>$.
The group $\O(3)$ acts orthogonally on $\ourspace$ in the canonical way;
that is, $A\in\O(3)$ maps $(XY)^*V$ to $(AX\cdot AY)^*AV$.
\item"(ii)"
Denote by $P:\ourspace\to\ourspace$ the linear map defined by
$P:(XY)^*\otimes V\mapsto\frac13((XY)^*\otimes V
+(XV)^*\otimes Y+(YV)^*\otimes X)$.
Obviously $P$ is an orthogonal projection onto its image $\Imm P$
which is isomorphic, as a representation space of~$\O(3)$,
to the ten dimensional space $\Sym^3\R^3$.
The projection~$P^\perp$ onto the orthogonal complement $\kernn P$ of
$\Imm P$ is given by $P^\perp:(XY)^*\otimes V\mapsto\frac13
(2(XY)^*\otimes V-(XV)^*\otimes Y-(YV)^*\otimes X)$.
Moreover, $P$ is $\O(3)$@-equivariant; thus $\Imm P$ and
$\kernn P$ are invariant under the $\O(3)$@-action.
\item"(iii)"
Denote by $\End_0(\R^3)$ the space of traceless endomorphisms of~$\R^3$.
Define a linear map $\Phi:\End_0(\R^3)\to\ourspace$ by letting
$$\Phi(b):\;X\cdot X\;\mapsto\;bX\times X
$$
for all $X\in\R^3$,
where~$\times$~denotes the vector product in~$\R^3$ and we interpret
$\Phi(b)\in\ourspace$ as a linear map from $\Sym^2(\R^3)$ to~$\R^3$
(obtained from the above formula by polarization and linear extension).
Obviously $\Phi$ maps $\End_0(\R^3)$ to $\kernn P$ since
$3\<P(\Phi(b))(X\cdot X),V\>=\<bX\times X,V\>+\<bX\times V,X\>+\<bV\times
X,X\>=0$ for all $X,V\in\R^3$. Moreover, $\Phi$ is $\SO(3)$@-equivariant,
where $\SO(3)$ acts on $\End_0(\R^3)$ by conjugation.
\endroster

\bigskip

\proclaim{Lemma 2.3}
$\Phi:\End_{\,0}(\R^3)\to\kernn P=\Imm P^\perp$ is an $\SO(3)$@-invariant
isomorphism.
In particular, $\ourspace=\Imm P\oplus\kernn P$
is isomorphic, as a representation space of $\SO(3)$, to the sum of
$\Sym^3(\R^3)$ and $\End_{\,0}(\R^3)$.
\endproclaim

\smallskip

\demo{Proof}
If $\Phi(b)=0$ for some $b\in\End_0(\R^3)$ then $bX\times X=0$ and thus
$bX\parallel X$ for all $X\in\R^3$. Hence $b$ is a multiple of the identity;
from $\tr(b)=0$ it follows that $b=0$. Therefore $\Phi$ is injective.
Since the dimension of $\kernn P$
equals $\dimm(\ourspace)-\dimm(\Imm P)=18-10=8=\dimm(\End_{\,0}(\R^3))$,
it follows that $\Phi$ is a vector space isomorphism.
\qed\enddemo

\bigskip

\proclaim{Proposition 2.4}
Denote by $S_0(\R^3)\subset\End_0(\R^3)$ the space of symmetric
traceless endomorphisms of~$\R^3$.
\roster
\item"(i)"
There exist pairs of linear maps $c,c':\R^2\to S_0(\R^3)$ such that the
following conditions are satisfied:
\newline
\quad\quad{\rm1.)}
For each $Z\in\R^2$ the elements $c_Z$ and $c'_Z$ of $S_0(\R^3)$
are conjugate by an element of $\SO(3)$.\newline
\quad\quad{\rm 2.)}
There is no $A\in\O(3)$ such that either $c'=I_A\circ c$ or $-c'=I_A\circ c$,
where~$I_A$ denotes conjugation by~$A$.
\item"(ii)"
Let $c,c'$ be as in~{\rm(i)}, and let $q:=\Phi\circ c$, $q':=\Phi\circ c':
\R^2\to\kernn P\subset\ourspace$. Then $q,q'$ have the following properties:
\newline
\quad\quad{\rm1.)}
For each $Z\in\R^2$ the elements $q_Z$ and $q'_Z$ of $\ourspace$
belong to the same $\SO(3)$@-orbit.\newline
\quad\quad{\rm2.)}
There is no $A\in\O(3)$ such that $q'=A\circ q$.
\item"(iii)"
An example of a pair $c,c':\R^2\to S_0(\R^3)$ satisfying the
conditions in~{\rm(i)} is given by
$$
\hskip29pt c_{Z_1}=c'_{Z_1}=\pmatrix -1&0&0\\0&0&0\\0&0&1\endpmatrix,\ \
c_{Z_2}=\pmatrix 0&1&0\\1&0&1\\0&1&0\endpmatrix,\ \
c'_{Z_2}=\pmatrix 0&0&\sqrt2\\0&0&0\\ \sqrt2&0&0\endpmatrix,
$$
where $\{Z_1\,,Z_2\}$ is the standard basis of~$\R^2$, and we identify
elements of $S_0(\R^3)$ with their matrix representation with respect
to the standard basis of~$\R^3$.
\endroster\endproclaim

\smallskip

\demo{Proof}
$\Phi:\End_0(\R^3)\to\kernn P$ is an isomorphism of
$\SO(3)$@-re\-presentation spaces by Lemma~2.3. The element $-\Id$ of $\O(3)$
acts as multiplication by~$-1$ on $\kernn P$ and as multiplication by~$1$
on $\End_0(\R^3)$. Therefore~(i) implies~(ii). We now show~(i).

Equivalence of two elements of $S_0(\R^3)$ modulo conjugation by orthogonal
endomorphisms can easily be checked, namely, by deciding whether the 
characteristic polynomials are equal. In the following we identify
elements of $S_0(\R^3)$ with their matrix representation with respect
to the standard basis.

Let $a:=\left(\smallmatrix -1&0&0\\0&0&0\\0&0&1\endsmallmatrix\right)$,
and consider all matrices~$b$ of the form
$\left(\smallmatrix 0&b_{12}&b_{13}
\\b_{12}&0&b_{23}\\b_{13}&b_{23}&0\endsmallmatrix\right)$.
Define $c^b:\R^2\to S_0(\R^3)$ by $c^b_{Z_1}:=a$ and $c^b_{Z_2}:=b$,
where $\{Z_1\,,Z_2\}$ is the standard basis of~$\R^2$.
We claim that there are pairs $b,b'$ such that
$c^b$ and~$c^{b'}$ satisfy the conditions~1.)~and~2.)~above.

Condition~1.)~is equivalent to the characteristic polynomials of
$sa+tb$ and $sa+tb'$ being equal for all $(s,t)\in\R^2$.
We have
$$\det\bigl(\la\Id-(sa+tb)\bigr)
=\la^3-\bigl(t^2(b_{12}^2+b_{13}^2+b_{23}^2)+s^2\bigr)\la
-st^2(b_{23}^2-b_{12}^2)-2\,t^3b_{12}b_{13}b_{23}\,.
$$
Thus condition~1.)~for $c^b$ and~$c^{b'}$
is equivalent to the following three equations being satisfied:
$$\aligned b_{12}^2+b_{13}^2+b_{23}^2&=b_{12}^{\prime\,2}+b_{13}^{\prime\,2}
  +b_{23}^{\prime\,2}\,,\\b_{23}^2-b_{12}^2
  &=b_{23}^{\prime\,2}-b_{12}^{\prime\,2}\,,\\
  b_{12}b_{13}b_{23}&=b'_{12}b'_{13}b'_{23}\,.\endaligned\tag5 
$$
Concerning condition~2.), note that if there does exist $A\in\O(3)$ such
that either $c^{b'}=I_A\circ c^b$ or $-c^{b'}=I_A\circ c^b$ then
$\pm a=AaA\inv$ and $\pm b'=AbA\inv$. In particular, we then have
$\tr(a^2b^2)=\tr(a^2b^{\prime2})$; that is,
$$b_{12}^2+2b_{13}^2+b_{23}^2=b_{12}^{\prime\,2}+2b_{13}^{\prime\,2}
+b_{23}^{\prime\,2}\,.
$$
Hence for all pairs $b,b'$ which satisfy \thetag{5}, but not the latter
equation, the corresponding pairs $c^b, c^{b'}$ satisfies 1.)~and 2.)~of~(i).
It is easy to see that many such pairs exist. One example is given by
$$b:=\left(\smallmatrix 0&1&0\\1&0&1\\0&1&0\endsmallmatrix\right)
  \quad\text{and}\quad
  b':=\left(\smallmatrix 0&0&\sqrt2\\0&0&0\\ \sqrt2&0&0\endsmallmatrix\right)
$$
The corresponding pair $c^b,c^{b'}$ is just the
specific example given in~(iii).
\qed\enddemo

\smallskip

\subheading{Notation and Remarks 2.5}

(i) Let $\R^2$ and $\R^3$ be equipped with the standard scalar products.
With each linear map $q:\R^2\to\ourspace$ we associate an
$\R^2$@-valued $1$@-form~$\la$ on~$S^2$ by letting
$\la_Z(X)=\<q_Z(p\cdot p),X\>$ for $X\in T_pS^2$.
Here, $\la_Z$~denotes $\<\la(\,.\,),Z\>$ as usual;
moreover, we have interpreted $q_Z\in\ourspace$ as a linear
map from $\Sym^2\R^3$ to~$\R^3$.
In other words, each $q_Z\in\ourspace$ defines a quadratic vector field
$v_Z:p\mapsto q_Z(p\cdot p)$ on~$\R^3$, and
$\la_Z$ is the pullback (by the inclusion $S^2\hookrightarrow\R^3$)
of the $1$@-form $\<q_Z(p\cdot p),\,.\,\>$ which is dual to~$v_Z$\,.

(ii) Let $q$ and $\hat q$ be two linear maps from~$\R^2$ to $\ourspace$.
For each $Z\in\R^2$ consider the vector field
$v_Z:p\mapsto(\hat q_Z-q_Z)(p\cdot p)$ on~$\R^3$.
If the image of $\hat q-q$ is contained in~$\Imm P$
then $\<dv_Z\restr p X,Y\>=2\<(\hat q_Z-q_Z)(p\cdot X),Y\>$ is symmetric in
$X,Y\in\R^3$ for each $p\in\R^3$.
This means that $v_Z$ satisfies the integrability conditions and thus
has a potential function on~$\R^3$.
Similarly, if the image of $\hat q-q$ is contained in
$\Phi(\text{Skew}(\R^3))$, where $\text{Skew}(\R^3)\subset\End_0(\R^3)$
denotes the skew-symmetric endomorphisms of~$\R^3$, then for each
$Z\in\R^2$ there exists $a_Z\in\R^3$ such that~$v_Z$ is of the form
$p\mapsto (p\times a_Z)\times p$, and thus the component of~$v_Z$ tangent
to~$S^2$ is the same as that of the constant vector field~$a_Z$\,.
In either of the two cases, it follows from Remark~2.1(i) that
the $\R^2$@-valued $1$@-forms $\la,\hat \la$
which are associated with~$q,\hat q$ as in~(i) differ by~$df$ for some
$f\in C^\infty(S^2,\R^2)$.

Thus if $\la$ is associated with any linear map $q:\R^2\to\ourspace$,
then there exists $f\in C^\infty(S^2,\R^2)$ such that $\hat\la:=\la+df$
is associated with a linear map $\hat q:\R^2\to\Phi(S_0(\R^3))$.
Recall from~1.5(v) (or Remark~2.1(i)) that the associated
metrics~$g\subla$ and~$g_{\widehat\la}$ on $S^2\times T^2$ are isometric.
This is the reason why we are interested only in
pairs~$q,q'$ with image in $\Phi(S_0(\R^3))$, as in Proposition~2.4.

\bigskip

\subheading{Example 2.6}

(i) Let $h$ be the standard round metric on~$S^2$, and let
$T^2:=\R^2/\LL$, where $\LL$ is any uniform lattice in~$\R^2$.
Let $T^2$ be equipped with the metric induced from the euclidean metric
on~$\R^2$.
For any pair of linear maps $q,q':\R^2\to\ourspace$ let $\la,\lap$ be the
corresponding $\R^2$@-valued $1$@-forms on~$S^2$ as in~2.5(i),
and let $g\subla$ and~$g\sublap$ be the associated metrics on
$S^2\times T^2$.

If $q$ and $q'$ satisfy condition~1.)~of
Proposition~2.4(ii), then $(S^2\times T^2,g\subla)$ and
$(S^2\times T^2,g\sublap)$ are isospectral by Proposition~1.8.
In fact, the condition implies that for each $Z\in\R^2$
there exists $A_Z\in\O(3)$ (even $A_Z\in\SO(3)$) such that
$$\la'_Z\drestr p=\<q'_Z(p\cdot p),\,.\,\>=\<A_Zq_Z(A_Z\inv p\cdot A_Z\inv
p),\,.\,\>=\la_Z\drestr{A_Z\inv p}\circ A_Z\inv=(A_Z^{-1\,*}\la_Z)\drestr p
$$
for all $p\in S^2$. Hence $\la'_Z=(A_Z^{-1})^*\la_Z\in\Isom(S^2,h)^*\la_Z$\,.
Thus condition~\thetag{$*2$} of Proposition~1.8 is satisfied, and
consequently the two manifolds are isospectral.

(ii) In the context of~(i), consider now the pair $q:=\Phi\circ c$ and
$q':=\Phi\circ c'$,  where $c,c'$ is the explicit pair of linear maps given
in Proposition~2.4(iii).
By~(i), the associated Riemannian manifolds
$(S^2\times T^2,g\subla)$ and $(S^2\times T^2,g\sublap)$ are isospectral.
In the next section we will show that they are not locally isometric.
More precisely, the preimage of the maximal scalar curvature in the first
manifold has greater dimension than the preimage of the maximal scalar
curvature in the second manifold (Proposition~2.10).

\bigskip

\subheading{Remark 2.7}
For amusement of the reader, we illustrate the specific pair of metrics
$g\subla\,,g\sublap$ on $S^2\times T^2$ from Example~2.6(ii) by explicitly
writing down the horizontal distributions. Let $\{Z_1\,,Z_2\}$ be the
standard basis of~$\R^2$, and denote the corresponding vector fields
on $S^2\times T^2$ by $Z_1\,,Z_2$ again. Then for each $p\in S^2$, the
$g\subla$@-orthogonal complement in $T_pS^2$ to the $T^2$@-fiber through~$p$
is given by the vectors
$$\allowdisplaybreaks
  \align X&-\left<\pmatrix -p_2p_3\\2p_1p_3\\-p_1p_2\endpmatrix,\,X\right>Z_1
  -\left<\pmatrix p_1p_3-p_2^2+p_3^2\\p_1p_2-p_2p_3\\-p_1p_3-p_1^2+p_2^2
  \endpmatrix,\,X\right>Z_2\,,\\
\intertext{where $X\in T_pS^2$ and $\scp$ is the standard scalar product
in~$\R^3$.
Similarly, the $g\sublap\,$@-horizontal distribution is given by the vectors}
X&-\left<\pmatrix -p_2p_3\\2p_1p_3\\-p_1p_2\endpmatrix,\,X\right>Z_1
   -\left<\pmatrix -\sqrt2\,p_1p_2\\ \sqrt2\,(p_1^2-p_3^2)\\ \sqrt2\,p_2p_3
   \endpmatrix,\,X\right>Z_2\endalign
$$
with $X\in T_pS^2$, $p\in S^2$.
This, together with the fact that for both metrics the $T^2$@-fibers are
isometrically embedded and that the projection to the round sphere
$(S^2,h)$ is a Riemannian submersion, describes $g\subla$ and~$g\sublap$
completely.

\bigskip

\smallskip

\subhead 2.2\ \ Curvature properties
  \endsubhead
\bigskip

\noindent
In this section we first compute the scalar curvature of the manifolds
$(N\times H,g\subla)$ from Notation~1.5(iv) in terms of~$\la$ and the scalar
curvature of $(N,h)$ (Proposition~2.8).
Applying this formula to the specific kind of metrics from the previous 
section we conclude, in particular, that the manifolds from Example~2.6(ii)
are not locally isometric because the preimages of their maximal
scalar curvatures have different dimensions (Proposition~2.10).
Finally, we make some general observations about the curvature
properties of pairs of isospectral manifolds $(N\times H,g\subla)$,
$(N\times H,g\sublap)$ arising from Proposition~1.8 in the special case
that $N$ is of dimension two (Theorem~2.11 and Remarks~2.12).

\bigskip

In the following we always assume that $H$ is a torus, $\h=T_eH$, and
$H$ is equipped with an invariant metric whose restriction to~$\h$
we denote by~$\scp$.
We fix an orthonormal basis $\{Z_1\,,\,\ldots,Z_r\}$ of~$\h$.

\bigskip

\proclaim{Proposition 2.8}
Let $(N,h)$ be a closed Riemannian manifold and $\la$ be an $\h$@-valued
$1$@-form on~$N$. Let $g\subla$ be the associated metric on $N\times H$
as in Notation~{\rm1.5(iv)}. Then we have
$$\scal^{g\subla}_{(p,z)}=\scal^h_p-\frac14\|d\la\restr p\|_h^2
$$
for all $(p,z)\in N\times H$, where $\scal^{g\subla}$ and~$\scal^h$
denote the scalar curvature functions associated with $g\subla$ and~$h$,
respectively, and $\|\,.\,\restr p\|^2_h$ denotes the euclidean norm on
tensors associated with the scalar product~$h\restr{T_pN}$\,.
\endproclaim

\smallskip

Note that our choice of the norm $\|\,.\,\|_h$ on tensors means, for example,
that the canonical volume form~$v$ associated to the standard metric~$h$
on~$S^{m-1}$ satisfies $\|v\restr p\|_h^2=(m-1)!$ for all $p\in S^{m-1}$.

\bigskip

\demo{Proof}
Let $Z\in\h$ and $X,Y$ be vector fields on~$N$.
We denote by the same names also the corresponding $H$@-invariant
vector fields on $N\times H$. By $\tilde X,\tilde Y$ we denote the associated
horizontal vector fields, that is, $\tilde X\subpz=(X_p\,,-\la_p(X))\in
T\subpz(N\times H)$. Note that $\la(X),\la(Y)$ commute since they are
$H$@-invariant and tangent to the fibers. Thus
$$\align
  g\subla([\tilde X,\tilde Y],Z)&=-g\subla([X,\la(Y)],Z)
  +g\subla([Y,\la(X)],Z)+g\subla([X,Y],Z)\\
  &=-X(\la_Z(Y))+Y(\la_Z(X))+\la_Z([X,Y])=-d\la_Z(X,Y).
\endalign
$$
Moreover, $[\tilde X,Z]=0$ since $\tilde X$ is $H$@-invariant.
By the Koszul formula $\nabla_ZZ=0$, $\nabla_Z\tilde X$ is horizontal,
and $g\subla(\nabla_Z\tilde X,\tilde Y)=-\frac12g\subla
([\tilde X,\tilde Y],Z)=\frac12d\la_Z(X,Y)$.
Since the projection to $(N,h)$ is a Riemannian submersion,
$\nabla_{\tilde X}\tilde X$ is horizontal and thus $0=Z(g\subla
(\nabla_{\tilde X}\tilde X,Z))=g\subla(\nabla_Z\nabla_{\tilde X}\tilde X,Z)$.
Hence $g\subla(R(Z,\tilde X)\tilde X,Z)=-g\subla(\nabla_{\tilde X}\nabla_Z
\tilde X,Z)\allowbreak=\|\nabla_Z
\tilde X\|^2_{g\subla}=\frac14\|d\la_Z(X,\,.\,)\|^2_h$\,.
By the flatness of the fibers this implies
$$\Ric^{g\subla}(Z,Z)=\tfrac14\|d\la_Z\|_h^2\,.
$$
Moreover, from O'Neill's formula and the above formula for the vertical part
of $[\tilde X,\tilde Y]$ it follows that
$$\Ric^{g\subla}(\tilde X,\tilde X)=\bigl(\Ric^h(X,X)-\tfrac34\tsize\sum_i
\|d\la_{Z_i}(X,\,.\,)\|_h^2\bigr)+\tfrac14\tsize\sum_i\|d\la_{Z_i}(X,\,.\,)
\|_h^2\,.
$$
Consequently $\scal^{g\subla}=\scal^h-\frac12\sum_i\|d\la_{Z_i}\|_h^2
+\frac14\sum_i\|d\la_{Z_i}\|_h^2$\,, which gives the desired formula.
\qed\enddemo

\bigskip

Next, we focus on the case $N=S^2$, $H=T^2$, and $\la$ depending
quadratically on the basepoint as in Section~2.1.
Recall from 2.2(iii) and~2.5(ii) that the relevant $1$@-forms~$\la$
in this category are of the type $\la_Z\restr p=\<c_Zp\times p,\,.\,\>$,
where~$c$ is a linear map from~$\R^2$ to the space $S_0(\R^3)$
of symmetric traceless endomorphisms of~$\R^3$.

\bigskip

\proclaim{Lemma 2.9}
Let $c:\R^2\to S_0(\R^3)$ be a linear map and $\la$ be the associated
$\R^2$@-valued $1$@-form on~$S^2$ as given by {\rm2.2(iii)} and~{\rm2.5(i)};
that is,
$\la_Z\restr p(X)=\<c_Zp\times p,X\>$ for all $p\in S^2$, $X\in T_pS^2$,
and $Z\in\R^2$, where $\scp$ denotes the standard scalar product on~$\R^3$.
Then for the associated metric~$g\subla$ on $S^2\times T^2$ we have
$$\scal\subpz^{g\subla}=2-\frac92\<c_{Z_1}p,p\>^2-\frac92\<c_{Z_2}p,p\>^2,
$$
where $\{Z_1\,,Z_2\}$ is the standard basis of~$\R^2$.
\endproclaim

\smallskip

\demo{Proof}
Let $p\in S^2$ and $\{X,Y\}$ be an orthonormal basis of $T_p S^2$ with
respect to the standard round metric~$h$ on~$S^2$.
Then $\|d\la\restr p\|^2=2\|d\la\restr p(X,Y)\|^2$.
By Proposition~2.8 and the fact that the scalar curvature of $(S^2,h)$
equals~$2$, the lemma will follow if we show that $\|d\la_Z\restr p(X,Y)\|^2=
9\<c_Zp,p\>^2$ for all $Z\in\R^2$. We interpret~$\la_Z$ as the restriction
of a $1$@-form on~$\R^3$ (defined by the same formula as~$\la_Z$) and
extend $X,Y$ to constant vector fields on~$\R^3$. Then
$$\align d\la_Z\restr p(X,Y)&=\<c_ZX\times p+c_Zp\times X,Y\>
  -\<c_ZY\times p+c_Zp\times Y,X\>\\
  &=\<c_ZX\times p,Y\>+2\<c_Zp\times X,Y\>-\<c_ZY\times p,X\>\\
  &=\<c_ZX\times p,Y\>+2\<c_Zp\times X,Y\>+\<Y\times c_Zp,X\>+\<Y\times p,
  c_ZX\>\\
  &=3\<c_Zp\times X,Y\>.\endalign
$$
Note that in the third equation we have used the fact that $\tr(c_Z)=0$
implies $\det(c_ZY,p,X)+\det(Y,c_Zp,X)+\det(Y,p,c_ZX)=0$.
Thus $\|d\la_Z\restr p(X,Y)\|^2=9\det(c_Zp,X,Y)^2=9\<c_Zp,p\>^2$,
as claimed.
\qed\enddemo

\bigskip

\proclaim{Proposition~2.10}
Let $c,c':\R^2\to S_0(\R^3)$ be the specific pair of linear maps given
in Proposition~{\rm2.4(iii)}, and $g\subla\,,g\sublap$ be the associated
metrics on $S^2\times T^2$ as in Example~{\rm2.6(ii)}.
Then the maximal scalar curvature equals~$2$ for both metrics.
Its preimage under $\scal^{g\subla}$ is
$$\bigl(\{p\in S^2\mid p_1=-p_3\}\cup\{(\pm1/\sqrt2,0,\pm1/\sqrt2)\}\bigr)
  \times T^2
$$
and thus contains a submanifold of codimension one,
whereas its preimage under $\scal^{g\sublap}$ is
$$\{(0,\pm1,0)\}\times T^2
$$
and thus is of codimension two.
\endproclaim

\smallskip

\demo{Proof}
By Lemma~2.9 we have that $\scal^{g\subla}$, respectively $\scal^{g\sublap}$,
attains its maximum precisely in those points where the function
$$\align
\<c_{Z_1}p,p\>^2+\<c_{Z_2}p,p\>^2&=(p_1^2-p_3^2)^2+(2p_2(p_1+p_3))^2,
\text{ respectively}\\
\<c'_{Z_1}p,p\>^2+\<c'_{Z_2}p,p\>^2&=(p_1^2-p_3^2)^2+(2\sqrt2p_1p_3)^2,
\endalign
$$
attains its minimum. Both minima are obviously zero and are attained
precisely in the sets given in the statement.
\qed\enddemo

\bigskip

We conclude this chapter by some general observations about the behaviour
of the scalar curvature functions associated with two isospectral manifolds
$(N\times H,g\subla)$, $(N\times H,g\sublap)$ arising from Proposition~1.8
in the case that $N$ is of dimension two.

\bigskip

\proclaim{Theorem 2.11}
Let $(N,h)$ be a closed Riemannian manifold and $H$ be a torus with Lie
algebra $\h=T_eH$, equipped with an invariant metric.
Let $\la,\lap$ be two $\h$@-valued $1$@-forms on~$N$ which satisfy
condition~\thetag{$*2$} of Proposition~{\rm1.8}, and
let $g\subla\,,g\sublap$ be the associated
isospectral metrics on $M=N\times H$. If $N$ is two-dimen\-sion\-al, then
the following holds:
\roster
\item"(i)" $\dsize\int_M(\scal^{g\subla})^k\dvol_{g\subla}
  =\dsize\int_M(\scal^{g\sublap})^k\dvol_{g\sublap}\;$ for all $k\in\N$.
\item"(ii)" {\vphantom{$\dsize\sum$}}The functions
  $\scal^{g\subla}$ and $\,\scal^{g\sublap}$ have the same range.
\endroster\endproclaim

\smallskip

\demo{Proof}

(i) We first consider the case that $N$ is orientable. We choose an
orientation and let~$v$ be the associated volume form on $(N,h)$,
that is, $v\restr p(X,Y)=1$ for a positively oriented $h$@-orthonormal
basis $\{X,Y\}$ of $T_pN$. Then there exist linear~(!)~maps
$\phi,\psi:\h\to C^\infty(N)$ such that $d\la_Z=\phi_Zv$,
$d\lap_Z=\psi_Zv$ for all $Z\in\h$. Since $A^*v=\pm v$ for $A\in\Isom(N,h)$,
and $\lap_Z\in\Isom(N,h)^*\la_Z$ by the condition of Proposition~1.8,
it follows that $\psi_Z^2\in\Isom(N,h)^*\phi_Z^2$ for all $Z\in\h$.
Thereby, for all $s\in\R^r$, $m\in\N$ and any $\Isom(N,h)$@-invariant
function~$f$ on~$N$ we have
$$\int_Nf\cdot(\phi_{s_1Z_1+\ldots+s_rZ_r}^2)^m\dvol_h
  =\int_Nf\cdot(\psi_{s_1Z_1+\ldots+s_rZ_r}^2)^m\dvol_h\,.
$$
Expanding this into monomials in the~$s_i$\,, using the linearity of
~$\phi$ and~$\psi$, we get
$$\int_Nf\tsize\prod_{i=1}^r\phi_{Z_i}^{n_i}\,\dvol_h
  =\int_Nf\tsize\prod_{i=1}^r\psi_{Z_i}^{n_i}\,\dvol_h\tag6
$$
for all $\Isom(N,h)$@-invariant functions~$f$ and all $n_1\,,\,\ldots,n_r
\in\N$ such that $\sum_in_i$ is even.
By Proposition~2.8 and the fact that $\|v\restr p\|_h^2\equiv2$,
we have
$$\int_M(\scal^{g\subla})^k\dvol_{g\subla}=\voll(H)\cdot
  \int_N(\scal^h-\tfrac12\tsize\sum_{i=1}^r\phi_{Z_i}^2)^k\dvol_h\,,
$$
and similarly
$$\int_M(\scal^{g\sublap})^k\dvol_{g\sublap}=\voll(H)\cdot
  \int_N(\scal^h-\tfrac12\tsize\sum_{i=1}^r\psi_{Z_i}^2)^k\dvol_h\,.
$$
Since $\scal^h$ and its powers are $\Isom(N,h)$@-invariant, these two
integrals decompose into two sums whose summands are of the type given
in \thetag{6} and thus match up pairwise.
This proves~(i) in the case that $N$ is orientable.

If $N$ is not orientable, consider its orientable double covering
$\pi:\bar N\to N$ and let $\bar h=\pi^*h$, $\bar\la=\pi^*\la$,
$\bar\lap=\pi^*\lap$. Clearly, $\bar\la$ and $\bar\lap$ again satisfy
the condition of Proposition~1.8 since every isometry of $(N,h)$ lifts
to an isometry of $(\bar N,\bar h)$. With respect to the associated
Riemannian metrics $g_{\bar\la}\,,g_{\bar\lap}$ on $\bar N\times H$,
the projections to $N\times H$ are Riemannian coverings. Since the
assertion for $g_{\bar\la}\,,g_{\bar\lap}$ holds by the above arguments,
it also follows for $g\subla\,,g\sublap$ by the fact that the integrals in
question equal just one half of the corresponding integrals over~$\bar M=
\bar N\times H$.

(ii) From~\thetag{6} we can also conclude, using the same argument
as in the proof of~(i), that
$$\int_M(a+b\,\scal^{g\subla})^k\dvol_{g\subla}=
  \int_M(a+b\,\scal^{g\sublap})^k\dvol_{g\sublap}
$$
for all $k\in\N$ and all $a,b\in\R$.
If we choose $a\ge-\min\{\min\scal^{g\subla},\min\scal^{g\sublap}\}$
and $b=1$ then $a+b\,\scal^{g\subla}$ and $a+b\,\scal^{g\sublap}$
induce nonnegative functions~$\Phi,\Psi$ on $(N,h)$ with the property
that for each $k\in\N$ their $L^k$@-norms coincide.
This implies that $\Phi$ and~$\Psi$ --- and consequently $\scal^{g\subla}$
and $\scal^{g\sublap}$ --- have the same maximum. In fact, we could
otherwise rescale $\Phi$ and~$\Psi$ simultaneously such that $\max\Phi<1$
and $\max\Psi>1$ (or vice versa). But then it would follow that
$\int_N\Phi^k\dvol_h\to0$ for $k\to\infty$, while $\int_N\Psi^k\dvol_h\to
\infty$ for $k\to\infty$, which is a contradiction.
A similar argument, using $a\ge\max\scal^{g\subla}=\max\scal^{g\sublap}$
and $b=-1$ shows that $\scal^{g\subla}$ and $\scal^{g\sublap}$ have the
same minimum.
\qed\enddemo

\smallskip

\subheading{Remarks 2.12}

(i) If in Theorem~2.11 the manifold $M=N\times H$ is four-dimensional
(that is, if $\dimm\,H=2$), then we can also conclude
$\int_M\|\Ric^{g\subla}\|^2\dvol_{g\subla}=\int_M\|\Ric^{g\sublap}\|^2
\dvol_{g\sublap}$ and $\int_M\|R^{g\subla}\|^2\dvol_{g\subla}=
\int_M\|R^{g\sublap}\|^2\dvol_{g\sublap}$\,. In fact,
$5\int\scal^2-2\int\|\Ric\|^2+2\int\|R\|^2$ is a heat invariant
(see \cite{18}, Theorem~4.8.18) and thus is the same for $g\subla$
and~$g\sublap$ because of their isospectrality; moreover,
$\int\scal^2-4\int\|\Ric\|^2+\int\|R\|^2$
is a topological invariant in dimension four because of the
Gauss-Bonnet-Chern formula (see, e.g., \cite{40}, p.~291).
Since the vectors $(5,-2,2)$, $(1,-4,1)$, and $(1,0,0)$ are linearly
independent, equality of $\int(\scal^{g\subla})^2$ and $\int(\scal^{
g\sublap})^2$ implies equality of the other two pairs of integrals too.

(ii) The results of Theorem~2.11 contrast with the properties of the
scalar curvature in isospectral examples with higher dimensional base
manifold~$N$. In \cite{22} examples were given of isospectral metrics on
$S^{m-1}\times T^2$ (interpretable as arising from Proposition~1.8;
see Example~1.11\,/\,Remark~1.12)
such that the associated scalar curvature functions have different maxima.
In \cite{42} it was shown that in these examples (and in a certain
generalization of them) the isospectral metrics have in general different
total squared curvatures and different total squared norms of the Ricci
tensor. (See Example~1.11, Remark~1.13(ii), Proposition~1.16, and
\cite{42}, Lemma~5.4.)

(iii) We finally return to the special case where $(N,h)$ is the standard
two-dimensional sphere and $\la,\lap$ are quadratic $\R^2$@-valued
$1$@-forms on~$S^2$ which are
associated with linear maps $c,c':\R^2\to S_0(\R^3)$, as in Lemma~2.9.
Assume that $c,c'$ satisfy the isospectrality condition~1.)~from
Proposition~2.4(i).
For the associated metrics $g\subla\,,g\sublap$
on $M=S^2\times T^2$ it is then possible to show:
$$\int_M(\Delta_{g\subla}\scal^{g\subla})^2\dvol_{g\subla}\ne
  \int_M(\Delta_{g\sublap}\scal^{g\sublap})^2\dvol_{g\sublap}
  \;\Longleftrightarrow\;
  \tr(c_{Z_1}^2c_{Z_2}^2)\ne\tr(c^{\prime\,2}_{Z_1}
  c^{\prime\,2}_{Z_2}).
$$
Recall from the proof of Proposition~2.4(i) that there exist many examples
where the latter is the case (one of them being the specific pair
$c,c'$ given in~2.4(iii) which was used for Example~2.6(ii)). The proof of
the above equivalence statement is quite elementary but somewhat tedious,
and we do not present it here.

\bigskip

\bigskip
\noindent
{\gross 3.\hbox to6pt{}Isospectral left invariant metrics on compact Lie
  groups}

\bigskip

\noindent
In the first section of this chapter, we will formulate versions of
Theorem~1.6 and
Proposition~1.8 for the special case of left invariant metrics on compact
Lie groups (see Proposition~3.1 and Corollary~3.2, respectively).
In Section~3.2 we will then give explicit applications.

We obtain the first examples of left invariant isospectral metrics on
compact Lie groups, even continuous families of such metrics.
This provides us with the first examples of continuous families of globally
homogeneous isospectral metrics.
(Note that there have previously been {\it pairs\/}
of globally homogeneous isospectral metrics, namely, pairs of isospectral
flat tori \cite{33}, \cite{13}. There also have been continuous families
of {\it locally\/} homogeneous isopectral manifolds, namely, isospectral
families of nil- and solvmanifolds; see, e.g., \cite{25}, \cite{41},
\cite{29}.)

In particular, we will obtain continuous isospectral families of left
invariant metrics on $\SO(m\ge5)\times T^2$, $\Spin(m\ge5)\times T^2$,
$\SU(m\ge3)\times T^2$, $\SO(n\ge8)$, $\Spin(n\ge8)$, and
$\SU(n\ge6)$.
Among these are the first examples of simply connected {\it irreducible\/}
isospectral manifolds. (The first examples of simply connected isospectral
manifolds, given in~\cite{42}, were products; recall Example~1.14.)
We also obtain the first continuous families of isospectral manifolds
of positive Ricci curvature (see Remark~3.12(ii)).

In Section~3.3 we prove that the isospectral homogeneous manifolds
constructed in Section~3.2 are not locally isometric
by computing the norm of their Ricci
tensors, which turns out to be in general nonconstant during the isospectral
deformations. From this we will also conclude, using heat invariants, 
that the manifolds are not
isospectral for the Laplace operator acting on $1$@-forms.

We will finish this chapter by proving in Section~3.4 that, although
the deformations from Section~3.2 can occur arbitrarily close to bi-invariant
metrics, a bi-invariant metric itself can never be
contained in a nontrivial continuous isospectral family of left invariant
metrics on a compact Lie group; in other words, bi-invariant metrics are
infinitesimally spectrally rigid within the class of left invariant metrics.

\bigskip

\smallskip

\subhead 3.1\ \ Application of the torus bundle construction to compact
  Lie groups\endsubhead
\bigskip

\proclaim{Proposition 3.1}
Let $G$ be a compact Lie group with Lie algebra~$\g=T_eG$, and let $g_0$
be a bi-invariant metric on $G$\,.
Let $H\subset G$ be a torus in~$G$ with Lie algebra
$\h\subset\g$.
Denote by~$\uu$ the $g_0$@-orthogonal complement of the centralizer~$\zh$
of~$\h$ in~$\g$.
Let $\la,\lap:\g\to\h$ be two linear maps with
$\la\restr{\h\oplus\uu}=\lap\restr{\h\oplus\uu}=0$ which satisfy:
\roster
\item"($*3$)" For every $Z\in\h$ there exists $a_Z\in G$ such that
  $a_Z$ commutes with~$H$ and $\la'_Z=\Ad_{a_Z}^*\la_Z$\,,
  where $\la_Z:=g_0(\la(\,.\,),Z)$ and $\lap_Z:=g_0(\lap(\,.\,),Z)$.
\endroster
Denote by $g\subla$ and $g\sublap$ the left invariant metrics on~$G$
which correspond to the scalar products $(\Id+\la)^*g_0$ and
$(\Id+\lap)^*g_0$ on~$\g$.
Then $(G,g\subla)$ and $(G,g\sublap)$ are isospectral.
\endproclaim

\smallskip

Note that $\Id+\la$, $\Id+\lap$ are indeed invertible maps since
$\la^2=\la^{\prime\,2}=0$. Thus the definition of $g\subla$
and $g\sublap$ makes sense.

We give two proofs of Proposition~3.1.
Although these proofs are related to each other, the
``geometric'' one uses Theorem~1.6, whereas the ``algebraic'' one is
self-contained and uses only the expression for the Laplacian on unimodular
Lie groups with a left invariant metric.

\bigskip

\subheading{Geometric proof of Proposition~3.1}

\noindent
We want to apply Theorem~1.6 to the torus~$H$, equipped with the restriction
of~$g_0$\,, and to $M:=G$ which we interpret as a principal $H$@-bundle with
respect to the right action of~$H$ on~$G$.
Denote by $\prh$ the $g_0$@-orthogonal projection from~$\g$ to~$\h$.
Let $\om:=\la+\prh$\,, $\omp:=\lap+\prh:\g\to\h$, and extend
$\om,\omp$ to left invariant, $\h$@-valued $1$@-forms on~$G$.
We claim that~$\om$ and~$\omp$ are invariant under the right action of~$H$.
First note that for every $Z\in\h$, the map $\ad_Z$ annihilates $\zh$
and thus, being $g_0$@-skew symmetric, maps $\g$ to~$\uu$.
Therefore, if $X\in\g$ and $z\in H$ we have indeed
$$\om(R_{z*}X)-\om(X)=\om(\Ad_z\inv X-X)\in\om([\g,\h])\subseteq
  \om(\uu)=0$$
by $\la\restr{\uu}=0$ and the definition of~$\om$; analogously for~$\omp$.
Moreover, $\la\restr\h=0$ implies $\om\restr\h=\omp\restr\h=\Id_\h$\,.
Hence we can view $\om,\omp$ as connection forms on the $H$@-bundle~$G$.

Let $N:=G/H$ and $h$ be the submersion metric on~$N$ induced by~$g_0$\,.
The definitions imply immediately that the metrics $g\subom$ and
$g\subomp$ (as defined in Notation~1.5(iii)) then equal $g\subla$ and
$g\sublap$\,, respectively.
In order to prove Proposition~3.1 it now suffices to check condition~($*1$)
of Theorem~1.6 for $\om$ and~$\omp$.

Let $Z\in\h$ and choose $a_Z\in G$ such that
$a_Z$ commutes with~$H$ and $\lap_Z=\Ad_{a_Z}^*\la_Z$\,.
Then $F_Z:=L_{a_Z}\circ R_{a_Z}\inv:G\to G$ is a bundle
automorphism satisfying $\omp_Z=F_Z^*\om_Z$\,. Moreover, since $g_0$
is bi-invariant, $F_Z$ is an isometry of $(G,g_0)$ and therefore
induces an isometry on $(N,h)=(G/H,g_0^H)$.
\qed

\bigskip

\subheading{Algebraic proof of Proposition~3.1}

\noindent
The group~$G$ acts on the Hilbert space $L^2(G)$ by
$(\rho_xf)=R_x^*f$ for all $f\in L^2(G)$
and $x\in G$. By unimodularity of~$G$ this action is unitary.
In particular, $H$ acts unitarily on $L^2(G)$ by the restriction
of~$\rho$ to~$H$.
Let $\LL:=\exp\inv(e)\cap\h\subset\h$, and denote by $\LL^*$ the dual
lattice in~$\h^*$. Then $L^2(G)=\oplus_{\mu\in\LL^*}\HH_\mu$\,,
where
$$\HH_\mu:=\{f\in L^2(G)\mid\rho_{\exp(Z)}f=e^{2\pi i\mu(Z)}f\text{ for all }
  Z\in\h\}.$$
We claim that $\Delta_{g\subla}$ and $\Delta_{g\sublap}$ leave each
$\CC_\mu:=C^\infty(G)\cap\HH_\mu$
invariant, and that their spectra on $\HH_\mu$ coincide.
Isospectrality of $(G,g\subla)$ and $(G,g\sublap)$ will then follow.

To prove our claim we first note that
if $g$ is any left invariant metric on~$G$ and $\{X_1\,,\,\ldots,X_d\}$
is a $g$@-orthonormal basis of~$\g$, then~$\Delta_g$ is given by
$-\sum_{i=1}^d X_i^2-\sum_{i=1}^d\nabla_{X_i}X_i$\,. The second
term is zero here because for each $Y\in\g$ we have
$\<\sum_{i=1}^d \nabla_{X_i}X_i\,,Y\>=-\sum_{i=1}^d\<X_i\,,[X_i\,,Y]\>
=\tr(\ad_Y)=0$ by unimodularity of the compact Lie group~$G$.
Thus
$$\Delta_g\;=\;-\sum_{i=1}^d X_i^2\;=\;-\sum_{i=1}^d
  (\rho_*X_i)^2.\tag7
$$
Now let $\{E_1\,,\,\ldots,E_d\}$ be a $g_0$@-orthonormal basis of~$\g$\,.
Then $\{E_1-\la(E_1),\,\ldots,E_d-\la(E_d)\}$ is a left invariant
$g\subla$@-orthonormal frame.
We assume $\{E_1\,,\,\ldots,E_d\}$ to be adapted to the $g_0$@-orthogonal
decomposition $\g=\zh\oplus\uu$; then in particular~$E_i$ and~$\la(E_i)$
commute for each~$i$ since $\la\restr\uu=0$. Consequently,
$$\align\Delta_{g\subla}\restr{\CC_\mu}&\;=\;-\sum_{i=1}^d
  \big(E_i^2-2\,E_i\circ\la(E_i)+\la(E_i)^2\big)
  \restr{\CC_\mu}\\
  &=\;-\sum_{i=1}^d\big(E_i^2-4\pi i\mu(\la(E_i))E_i-4\pi^2
    \mu(\la(E_i))^2\Id
    \big)\restr{\CC_\mu}\\
  &=\;\big(\Delta_{g_0}+4\pi iY_\mu(\la)+4\pi^2\|Y_\mu(\la)\|_{g_0}^2\Id\big)
    \restr{\CC_\mu}\,,\tag8\endalign$$
where $Y_\mu(\la):=\sum_{i=1}^d \mu(\la(E_i))E_i$\,.
Note that for each $x\in G$ the map $R_x$ is an isometry with respect to
the bi-invariant metric~$g_0$\,. Therefore $\Delta_{g_0}$ commutes with~$R_x$
and hence with~$\rho_x$ for each $x\in G$. In particular, $\CC_\mu$
is invariant under $\Delta_{g_0}$\,.
We claim that $\CC_\mu$ is also invariant under $Y_\mu(\la)$.
Let $Z_\mu$ be the dual vector to~$\mu$ with respect to $g_0\restr\h$\,,
and denote by $\trans\la:\h\to\g$ the transpose of~$\la$ with respect
to~$g_0$\,. Then
$$Y_\mu(\la)\;=\;\sum_{i=1}^d g_0(Z_\mu\,,\la(E_i))E_i
  =\trans\la(Z_\mu)\subseteq\uu^\perp=\zh$$
since $\uu\subset\kernn\la$. But $Y_\mu(\la)\in\zh$ implies that
$\CC_\mu$ is indeed invariant under $Y_\mu(\la)$. Thus $\CC_\mu$ is
invariant under $\Delta_{g\subla}$ by~\thetag{8}, and analogously under
$\Delta_{g\sublap}$\,.

Moreover, if we let $a_\mu:=a_{Z_\mu}$ be as in~\thetag{$*3$},
then~$\rho_{a_\mu}$~too leaves~$\CC_\mu$ invariant since~$a_\mu$ commutes
with~$H$.
By $\lap_{Z_\mu}=\Ad_{a_\mu}^*\la_{Z_\mu}$ and the bi-invariance of~$g_0$
we have
$$Y_\mu(\lap)=\trans\lap(Z_\mu)=\trans\Ad_{a_\mu}(\trans\la(Z_\mu))
  =\Ad_{a_\mu}\inv(Y_\mu(\la)).
$$
This implies $\|Y_\mu(\la)\|_{g_0}^2=\|Y_\mu(\lap)\|_{g_0}^2$\,,
and by~\thetag{8}:
$$\Delta_{g\sublap}\restr{\CC_\mu}\;=\;
  (\rho_{a_\mu}\inv\circ\Delta_{g\subla}\circ\rho_{a_\mu})
  \restr{\CC_\mu}\,.
$$
Therefore the spectra of $\Delta_{g\subla}$ and $\Delta_{g\sublap}$
on $\HH_\mu$ coincide, as claimed.
\qed

\bigskip

We conclude this section with a corollary of Proposition~3.1
which follows as well from Proposition~1.8 and can be regarded
as the intersection of both.

\bigskip

\proclaim{Corollary 3.2}
Let $K$ be a compact Lie group with Lie algebra $\k=T_eK$, and let
$h$ be a bi-invariant metric on~$K$.  Let $H$ be a torus with
Lie algebra $\h:=T_eH$, equipped with an invariant metric.
Let $\la,\lap:\k\to\h$ be two linear maps which satisfy:
$$\lap_Z\in\Ad_K^*(\la_Z)\text{ for each }Z\in\h.\tag{$*4$}
$$
Then $(K\times H,g_\la)$ and $(K\times H,g_\lap)$ are isospectral,
where $\la,\lap$ are interpreted as left invariant $\h$@-valued $1$@-forms
on~$K$, and $g\subla$\,, $g\sublap$ are the associated metrics on
$K\times H$ as in Notation~{\rm1.5(iv)}.
\endproclaim

\smallskip

\demo{Proof}
The corollary follows immediately from either Proposition~1.8 or
Proposition~3.1.
In the context of Proposition~1.8, the manifold $(K,h)$ plays the
role of $(N,h)$, and it suffices to note that the inner automorphisms
of~$K$ are isometries with respect to~$h$.

In the context of Proposition~3.1, the group $K\times H$, equipped with
the product metric, plays the role of~$(G,g_0)$; we extend $\la,\lap$
to linear maps from $\k\oplus\h$ to~$\h$
by letting $\la\restr{\h}=\lap\restr{\h}=0$.
It suffices to note that now~$H$ is central in~$G$ and thus $\uu=0$; hence
$\la\restr\uu=\lap\restr\uu=0$ is trivially satisfied.
Condition~\thetag{$*3$} is implied
by~\thetag{$*4$} and the fact that each $a_Z\in K$ commutes with~$H$.
\qed\enddemo

\bigskip

\smallskip

\subhead 3.2\ \ Examples\endsubhead
\bigskip

\noindent
We will now exploit Proposition~3.1 to construct the first examples of
isospectral left invariant metrics on compact Lie groups, as announced
at the beginning of this chapter.

In Subsection~3.2.1, we apply Proposition~3.1 via the more special
Corollary~3.2 to obtain continuous isospectral families of left invariant
metrics on $\SO(m)\times T^2$ and $\Spin(m)\times T^2$ for $m\ge5$,
and also on $\SU(m)\times T^2$ for $m\ge3$.

In the subsections~3.2.2\,/\,3.2.3 we will then use Proposition~3.1 in its
general form to find continuous isospectral families of left invariant
metrics on the irreducible groups
$\SO(n)$ and $\Spin(n)$ for $n\ge8$ and on $\SU(n)$ for $n\ge6$.

\bigskip

\bigskip

\noindent
{\bf 3.2.1\ \ Isospectral deformations on $\SO(m)\times T^2$
  ($m\ge5$),\newline
\phantom{3.2.1\ \ }$\Spin(m)\times T^2$ ($m\ge5$), 
  and $\SU(m)\times T^2$ ($m\ge3$).}\nopagebreak
\bigskip
\nopagebreak
\noindent
One application of Corollary~3.2 has already been waiting in a barely
disguised form, as we are going to see now. The main tool is the 
existence of nontrivial
isospectral families of linear maps $j(t):\R^2\to\so(m)$ with $m\ge5$
which is guaranteed by Proposition~1.16(ii) from the last chapter
and which were the key tool in Example~1.11 and~1.14.

\bigskip

\subheading{Example~3.3}

\noindent
Let $K=\SO(m)$ and $\k:=\so(m)=T_eK$; assume $m\ge5$.
Consider a bi-invariant metric~$h$ on~$K$ (unique up to scaling)
and the associated $\Ad_K$@-invariant scalar product on~$\k$.
Let $H$ be a two-dimensional torus with Lie algebra $\h=T_eH$,
equipped with some invariant metric.
Let $\{Z_1\,,Z_2\}$ be an orthonormal basis of~$\h$.

Recall from Proposition~1.16(ii) that there exists a Zariski open
subset~$\Cal U$
of the space~$\J$ of linear maps $j:\h\to\so(m)=\k$ such that for each
$j\in\Cal U$ there is a continuous family~$j(t)$ in~$\J$, defined in some
open neighbourhood of~$t=0$, such that $j(0)=j$ and such that:
\roster
\item"1.)" The maps~$j(t)$ are pairwise isospectral in the sense of
  Definition~1.10.
\item"2.)" The function $t\mapsto\|j_{Z_1}(t)^2+j_{Z_2}(t)^2\|^2=
  \tr\bigl((j_{Z_1}(t)^2+j_{Z_2}(t)^2)^2\bigr)$ is nonconstant in~$t$
  in every interval around zero.
\endroster
An explicit example for a family $j(t)$ satisfying 1.)~and~2.)~in case
$m=5$ was given in Proposition~1.16(iv).

Now let $\{j(t)\}_{t\in(-\eps,\eps)}$ be any continuous family in~$\J$ which
satisfies 1.)~and~2.). Recall that condition~1.)~means that for each
$Z\in\h$ the path $t\mapsto j_Z(t)$ is contained in the
$\Ad_{\O(m)}$@-orbit of $j_Z(0)$. By continuity it follows that
$j_Z(t)$ must even be contained in the $\Ad_{\SO(m)}$@-orbit of $j_Z(0)$.
Define linear maps $\la(t):\k\to\h$ by letting
$$\la_Z(t):=\<\,.\,,j_Z(t)\>
$$
for each $Z\in\h$, where $\la_Z(t)$ means $\<(\la(t))(\,.\,),Z\>$ as
usual. In other words, $\la(t):\k\to\h$ is the transpose of
$j(t):\h\to\k$ with respect to the given metrics.
By the bi-invariance of~$h$ and the fact that for any fixed $Z\in\h$
we have $j_Z(t)\in\Ad_K(j_Z(0))$ for all~$t$,
it follows that $\la_Z(t)\in\Ad_K^*(\la_Z(0))$ for all~$t$.
But this just means that the maps~$\la(t)$ pairwise satisfy
condition~\thetag{$*4$} of Corollary~3.2.
We conclude that the Riemannian manifolds
$$(\SO(m)\times H,g_{\la(t)})
$$
are isospectral, where $g_{\la(t)}$ is the left invariant metric
associated with $\la(t)$ and~$h$ as in Corollary~3.2.

Note that instead of $K=\SO(m)$ we may as well consider its universal
covering $\tilde K:=\Spin(m)$ because $\Ad_{\tilde K}$@-orbits
in~$\k$ are the same as $\Ad_K$@-orbits.
Hence our above families $\la(t):\k\to\h$ satisfy condition~\thetag{$*4$}
of Corollary~3.2 also with respect to~$\tilde K$.
Thus for each $m\ge5$ we also get isospectral families
$$(\Spin(m)\times H,g_{\la(t)}),
$$
where $H\cong T^2$ is as above and $g_{\la(t)}$ is the left invariant metric
which is associated, as in Corollary~3.2, with~$\la(t)$
and the bi-invariant metric~$\tilde h$ on $\Spin(m)$ which is the pullback
of the above metric~$h$ on $\SO(m)$.
With these notations the projection from $(\Spin(m)\times H,
g_{\la(t)})$ to $(\SO(m)\times H,g_{\la(t)})$ is a Riemannian covering;
thus isospectrality of the manifolds in the latter family is actually
implied by continuity and by the isospectrality of the covering manifolds.

For each~$t$, the norm of the associated Ricci tensor $\Ric^{g_{\la(t)}}$
is a constant function on $\SO(m)\times H$ (resp\. $\Spin(m)\times H$)
since $g_{\la(t)}$ is left invariant.
We claim that from the above conditions 1.),~2.)~on $j(t)={}^T\la(t)$
it follows that:
\roster
\item"(i)"
$\|\Ric^{g_{\la(t)}}\|^2$ is nonconstant in~$t$. In particular,
the manifolds $(\SO(m)\times H,g_{\la(t)})$ are not pairwise locally
isometric.
\item"(ii)"
The second heat invariant for the Laplace operator on $1$@-forms associated
with $(\SO(m)\times H,g_{\la(t)})$ depends nontrivially on~$t$. 
In particular, the manifolds are not pairwise isospectral on $1$@-forms.
\endroster
The analogous statements hold for $(\Spin(m)\times H,g_{\la(t)})$.
We postpone the proof of these facts to Section~3.3
(see Theorem~3.14, Proposition~3.15, and Corollary~3.17).

\bigskip

\subheading{Remark 3.4}
The relation between the isospectral metrics on $\SO(m)\times T^2$ from
Example~3.3 on the one hand and those on $S^{m-1}\times T^2$ from
Example~1.11 on the other hand can be explained by the following general
principle.

Let $K$ be a compact Lie group with Lie algebra~$\k$, equipped with a
bi-invariant metric, and suppose that $K$~acts by isometries on a closed
Riemannian manifold $(N,h)$. Each $j_Z\in\k$ is canonically identified
with the Killing vector field $p\mapsto\tdiffzero t\exp(tj_Z)p$ on~$N$;
by taking duals on both sides, each $\la_Z\in\k^*$ is canonically
identified with a certain $1$@-form on~$N$. If $\la_Z\,,\lap_Z$
belong to the same coadjoint orbit then the associated $1$@-forms on~$N$
are related by an element of $K\subseteq\Isom(N,h)$. In other words,
two linear maps $\la,\lap:\k\to\h\cong\R^r$ satisfying the isospectrality
condition~\thetag{$*4$} from Corollary~3.2 produce an associated pair
of $\h$@-valued $1$@-forms on~$N$ which satisfies the isospectrality
condition~\thetag{$*2$} from Proposition~1.8, and this for each
Riemannian manifold $(N,h)$ on which~$K$ acts by isometries.

In our above context, $K=\SO(m)$ and $(N,h)$ is the round standard
sphere~$S^{m-1}$. Obviously the above point of view opens prospects
for examples of isospectral metrics not only on $\SO(m)\times T^2$
or $S^{m-1}\times T^2$, but also on $N\times T^2$ where~$N$ is, for
example, a Grassmannian, or any other manifold admitting a metric~$h$
with respect to which $K=\SO(m\ge5)$ acts by isometries.
In view of Example~3.7 below, the same considerations are also
valid for $K=\SU(m\ge3)$.
We will not pursue these ideas in the present work, but concentrate
entirely on the Lie groups themselves.

\bigskip

Our next aim is to construct isospectral metrics on $\SU(m)\times T^2$
by analogous methods as those used above for $\SO(m)\times T^2$.
As we will see, this is indeed possible for $m\ge3$.
First we need a result analogous to the one cited in Proposition~1.16(ii)
(that is, to~\cite{42}, Proposition~4.1).

\bigskip

\definition{Definition 3.5}
Two linear maps $j,j':\R^r\to\su(m)$ are called {\it isospectral},
denoted $j\sim j'$, if for every $Z\in\R^r$ there exists $A_Z\in\SU(m)$
such that $j'_Z=A_Z j_Z A_Z\inv$.
\enddefinition

\smallskip

\proclaim{Proposition 3.6}
Let $m\ge3$ and $\{Z_1\,,Z_2\}$ be the standard basis of~$\R^2$.
\roster
\item"(i)"
There exists a Zariski open subset~$\Cal U$
of the space~$\J$ of linear maps $j:\R^2\to\su(m)$ such that for each
$j\in\Cal U$ there is a continuous family~$j(t)$ in~$\J$, defined in some
open neighbourhood of~$t=0$, such that $j(0)=j$ and:\newline
\qquad\qquad{\rm1.)} The maps~$j(t)$ are pairwise isospectral in the sense of
  Definition~{\rm3.5}.\newline
\qquad\qquad{\rm2.)} The function $t\mapsto\|j_{Z_1}(t)^2+j_{Z_2}(t)^2\|^2=
  \tr\bigl((j_{Z_1}(t)^2+j_{Z_2}(t)^2)^2\bigr)$ is not constant in~$t$
  in any interval around zero.
\item"(ii)"
For $m=3$, an explicit example of an isospectral family~$j(t):\R^2\to
\su(3)$ with
$\|j_{Z_1}(t)^2+j_{Z_2}(t)^2\|^2\ne
\text{\rm const}$ is given by
$$j_{Z_1}(t):=\left(\smallmatrix -i&0&0\\0&0&0\\0&0&i\endsmallmatrix\right),
\qquad j_{Z_2}(t):=\left(\smallmatrix 0&t&\sqrt{1-2t^2}\\-t&0&t\\-\sqrt
{1-2t^2}&-t&0\endsmallmatrix\right).
$$
This family is defined for $t\in[-1/\sqrt2,1/\sqrt2]$.
The $j(t)$ are pairwise isospectral since
$\det\bigl(\la\Id-(sj_{Z_1}(t)+uj_{Z_2}(t))\bigr)=\la^3+(s^2+u^2)\la$
is independent of~$t$. However, $\|j_{Z_1}(t)^2+j_{Z_2}(t)^2\|^2=8-4t^2$
is nonconstant in~$t$.
\endroster\endproclaim

\smallskip

\demo{Proof}
Part (ii) can be checked by straightforward computation.
As to part~(i), the proof of \cite{42}, Proposition~4.1, which asserted
the analogous statement for $\so(m\ge5)$ (see
Proposition~1.16(ii)) instead of $\su(m\ge3)$ carries
over almost verbatim.

First of all, elementary arguments show that two elements of~$\su(m)$
are conjugate by an element of $\SU(m)$ if and only if they have the
same characteristic polynomials. Therefore the condition $j\sim j'$
is equivalent to $sj_1+uj_2$ having the same characteristic
polynomial as $sj'_1+uj'_2$ for all~$s,u$, where we write
$j_1:=j_{Z_1}$\,, $j_2:=j_{Z_2}$\,. This is in turn
equivalent to $\tr\bigl((sj_1+uj_2)^k\bigr)=\tr\bigl((sj'_1+
uj'_2)^k\bigr)$ for all $k=1,\,\ldots,m$, or equivalently, for all $k\in\N$.
By expanding into monomials in~$s,u$ we get that
$$j\sim j'\;\;\Longleftrightarrow\;\;p\ab(j)=p\ab(j')\text{ for all
  $a,b\in\N_0$ with $a+b>0$},
$$
where
$$p\ab(j):=\sum_{\sig\in\sab}\tr(j_{\sig(1)}\ldots j_{\sig(a+b)})
$$
and $\sab$ denotes the set of all maps $\sig:\{1,\ldots,a+b\}\to\{1,2\}$
which satisfy $\#\sig^{-1}(1)=a$, $\#\sig^{-1}(2)=b$.

The algebraic vector field~$Y$ on~$\J$, given by
$Y(j)=(j_1^5\,j_2-j_2\,j_1^5\,,0)$ satisfies $dp\ab\restr jY=0$ for all 
$j\in\J$. The proof of this fact is purely combinatoric and reads exactly
as the proof of \cite{42}, Lemma~4.3, except that there we used the
exponent~3 instead of~5.
This implies that the (locally defined) flow lines of~$Y$ consist of pairwise
isospectral maps.

For $j\sim j'$, we obviously have $\tr(j_1^4)=\tr(j_1^{\prime\,4})$ and
$\tr(j_2^4)=\tr(j_2^{\prime\,4})$; hence in this case, the condition
$\tr\bigl((j_1^2+j_2^2)^2\bigr)=\tr\bigl((j_1^{\prime\,2}
+j_2^{\prime\,2})^2\bigr)$ is equivalent to $q(j)=q(j')$,
where $q(j):=\tr(j_1^2j_2^2)$.

We have $dq\restr jY=\tr(j_1^5j_2j_1j_2^2-j_1j_2j_1^5j_2^2)$.
This polynomial does not vanish identically on~$\J$ if $m\ge3$;
e.g., for 
$$j_1=\left(\smallmatrix i&0&0\\0&2i&0\\0&0&-3i\endsmallmatrix\right)
\text{ and }j_2=\left(\smallmatrix 0&1&1\\-1&0&1\\-1&-1&0
\endsmallmatrix\right)
$$
it equals $240\ne0$.
Therefore, the Zariski open subset $\Cal U:=\{j\in\J\mid dq\restr jY\ne0\}$
has the required properties.
\qed\enddemo

\bigskip

We now proceed in analogy with Example~3.3.

\bigskip

\subheading{Example 3.7}

\noindent
Let $K:=\SU(m)$ and $\k:=\su(m)=T_eK$; assume $m\ge3$.
Consider a bi-invariant metric~$h$ on~$K$ (unique up to scaling)
and the associated $\Ad_K$@-invariant scalar product on~$\k$.
Let $H$ be a two-dimensional torus with Lie algebra $\h=T_eH$,
equipped with some invariant metric. Let $\{Z_1\,,Z_2\}$ be an
orthonormal basis of~$\h$.

Let $\{j(t)\}_{t\in(-\eps,\eps)}$ be any continuous family of linear 
maps from~$\h$ to $\su(m)$ which
satisfies conditions 1.)~and~2.)~of Proposition~3.6(i).
Define linear maps $\la(t):\k\to\h$ by letting
$$\la_Z(t):=\<\,.\,,j_Z(t)\>
$$
for each $Z\in\h$; that is, $\la(t)={}^Tj(t)$ with respect to the
given metrics on $\k$ and~$\h$. From condition~1.)~it
follows that $\la_Z(t)\in\Ad_K^*(\la_Z(0))$ for all~$t$.
Thus the maps~$\la(t)$ pairwise satisfy
condition~\thetag{$*4$} of Corollary~3.2.
Hence the Riemannian manifolds
$$(\SU(m)\times H,g_{\la(t)})
$$
are isospectral, where $g_{\la(t)}$ is the left invariant metric
associated with $\la(t)$ and~$h$ as in Corollary~3.2.

In Section~3.3 we will prove that the conditions 1.)~and~2.)~imply the same
properties as in Example~3.3; that is:
\roster
\item"(i)"
$\|\Ric^{g_{\la(t)}}\|^2$ is nonconstant in~$t$. In particular,
the manifolds $(\SU(m)\times H,g_{\la(t)})$ are not pairwise locally
isometric (see Theorem~3.14 and Proposition~3.15).
\item"(ii)"
The second heat invariant for the Laplace operator on $1$@-forms associated
with $(\SU(m)\times H,g_{\la(t)})$ depends nontrivially on~$t$. 
In particular, the manifolds are not pairwise isospectral on $1$@-forms
(Corollary~3.17).
\endroster

\bigskip

\smallskip

\noindent
{\bf 3.2.2\ \ Isospectral deformations on $\SO(n)$ ($n\ge9$),
  $\Spin(n)$ ($n\ge9$),\newline
\phantom{3.2.2\ \ }and $\SU(n)$ ($n\ge6$).}
\bigskip

\noindent
We will now use the ideas from the previous subsection to construct
isospectral left invariant metrics on {\it irreducible\/} compact
Lie groups.
More precisely, we embed the above products $\SO(m)\times T^2$ etc\.
into certain irreducible groups and use Proposition~3.1 to obtain
isospectral metrics on these.

\bigskip

\subheading{Example 3.8}

\noindent
Let $G:=\SO(m+4)$ and $\g:=\so(m+4)=T_eG$; assume $m\ge5$.
Let $g_0$ be a bi-invariant metric on~$G$, and denote the corresponding
$\Ad_G$@-invariant scalar product on~$\g$ by~$g_0$ again.
Let $K_1:=\SO(m)$ and $K_2:=\SO(4)$. Since $\SO(m)\times\SO(4)$
is canonically embedded in~$\SO(m+4)$, we will from now on consider
$K_1$ and~$K_2$ as commuting subgroups of~$G$ which are orthogonal
with respect to the Killing metric on~$G$, and consequently with respect
to~$g_0$\,.
Let $H$ be a maximal torus in~$K_2$\,, endowed with the invariant metric
induced by~$g_0$\,. We denote by $\g$, $\h$, $\k_1$\,, $\k_2$ the Lie
algebras of $G$, $H$, $K_1$\,, and $K_2$\,, respectively.
Note that
$$[\k_1\,,\h]=0\quad{\text{ and }}\quad\k_1\perp_{g_0}\h,\tag{9}
$$
where $\perp_{g_0}$ denotes orthogonality with respect to~$g_0$\,.
Since $H$ is two-dimensional, there exist continuous families of linear
maps $j(t):\h\to\k_1$ satisfying the conditions~1\.) and~2\.)
from Example~3.3.
(Recall that there even exists a Zariski open subset~$\Cal U$ of the
space of linear maps from~$\h$ to~$\k_1$ such that each element of~$\Cal U$
is contained in a continuous family satisfying 1.)~and~2.).)

Let $\{j(t)\}_{t\in(-\eps,\eps)}$ be such a family. As in Example~3.3
we conclude from condition~1.)~that for each $Z\in\h$ we have
$j_Z(t)\in\Ad_{K_1}(j_Z(0))$ for all~$t$.
We now interpret $j(t):\h\to\k_1\subset\g$ as a linear map from
$\h$ to~$\g$ and define $\la(t):={}^Tj(t):\g\to\h$
as the transpose of~$j(t)$ with respect to~$g_0$\,.

We claim that the maps~$\la(t)$ pairwise satisfy the conditions of
Proposition~3.1.
In fact, by~\thetag{9} we have that $\k_1$ is $g_0$@-orthogonal both
to~$\h$ and~$\uu$, where $\uu$ is the $g_0$@-orthogonal complement of the
centralizer~$\z(\h)$ of~$\h$ in~$\g$.
Thus $\h\oplus\uu\perp\k_1\supseteq\Imm j(t)$, which implies
that $\h\oplus\uu\subseteq\kernn\la(t)$ for all~$t$.
Hence the first condition of Proposition~3.1 is satisfied.
Moreover, for each $Z\in\h$ we have $\la_Z(t)\in\Ad_{K_1}^*(\la_Z(0))$
for all~$t$ by the analogous property of the~$j_Z(t)$. Since $K_1$ commutes
with~$H$, condition~\thetag{$*3$} of Proposition~3.1 is satisfied, too.
We thus get isospectral families
$$(\SO(n),g_{\la(t)})
$$
for each $n=m+4\ge 5+4=9$, where $g_{\la(t)}$ is the left invariant
metric associated with $\la(t)$ and~$g_0$ as in the proposition.

Instead of $G,K_1\,,K_2$ we may as well consider their universal coverings
$\tilde G:=\Spin(m+4)$ and $\tilde K_1\times\tilde K_2:=\Spin(m)
\times\Spin(4)\subset\tilde G$,
endowed with a bi-invariant metric~$\tilde g_0$\,. Note that
$\Ad_{\tilde K_1}$@-orbits in~$\k_1$ are the same as
$\Ad_{K_1}$@-orbits, and that $\h$ is the Lie algebra of some
two-dimensional torus~$\tilde H$ in~$\tilde K_2$ which commutes with~$\tilde
K_1$ and is $\tilde g_0$@-orthogonal to~$\tilde K_1$\,.
Hence our above arguments go through to show that the family of linear maps
$\la(t):\g\to\h$ satisfies the conditions of Proposition~3.1 also with
respect to $\tilde G$, $\tilde H$, and~$\tilde g_0$\,.
Thus we also obtain isospectral families
$$(\Spin(n),g_{\la(t)})
$$
for each $n\ge9$, where $g_{\la(t)}$ is the left invariant metric
associated with $\la(t)$ and~$\tilde g_0$ as in Proposition~3.1.

Concerning local non-isometry, we have by Proposition~3.15 of Section~3.3
below that $\|\Ric_{g_{\la(t)}}\|^2$ is nonconstant in~$t$, and by
Corollary~3.17 the manifolds are not pairwise isospectral on $1$@-forms.

\bigskip

\subheading{Example 3.9}

\noindent
We replace the groups $G, K_1\,, K_2$ appearing in Example~3.8
by $G:=\SU(m+3)$ with $m\ge3$, $K_1:=\SU(m)$, and $K_2:=\SU(3)$.
Again, we consider $K_1$ and $K_2$ as commuting subgroups of~$G$ which
are orthogonal with respect to the Killing metric.
We choose a bi-invariant metric~$g_0$ on~$G$ (in particular, $g_0$ is
proportional to the Killing metric) and a maximal, hence two-dimensional,
torus~$H$ in~$K_2$\,.

Using the isospectral families $j(t):\h\to\k_1\subset\g$
from Proposition~3.6(i) this time, we obtain continuous families of linear
maps $\la(t):={}^Tj(t):\g\to\h$ which pairwise satisfy the conditions
of Proposition~3.1; the arguments read exactly as in Example~3.8.
We thus obtain isospectral families
$$(\SU(n),g_{\la(t)})
$$
for each $n=m+3\ge3+3=6$, where $g_{\la(t)}$ is the left invariant
metric associated with~$\la(t)$ and~$g_0$ as in Proposition~3.1.

Concerning non-isometry and non-isospectrality on $1$@-forms,
see again Proposition~3.15\,/\,Corollary~3.17 below.

\bigskip

\smallskip

\noindent
\subhead 3.2.3\ \ Isospectral deformations on $\SO(8)$ and $\Spin(8)$
  \endsubhead
\bigskip

\noindent
Recall that in the previous subsection, we modified Example~3.3
(respectively~3.7) to obtain isospectral families of metrics on the
irreducible Lie groups of Example~3.8 (respectively~3.9)
by considering the canonical embeddings
$K_1\times H=\SO(m)\times T^2\hookrightarrow\SO(m)\times\SO(4)\hookrightarrow
\SO(m+4)$
and $K_1\times H=\SU(m)\times T^2\hookrightarrow\SU(m)\times\SU(3)
\hookrightarrow\SU(m+3)$, respectively.
The lowest dimensions of irreducible examples obtained in this way
are $35=\dimm(\SU(6))$ and $36=\dimm(\SO(9))$.
As a final application of Proposition~3.1, we now construct isospectral
families of left invariant metrics on the $28$@-dimensional irreducible
Lie groups $\SO(8)$ and~$\Spin(8)$ by arranging~$K_1$ and~$H$ in a
more economical way.

\bigskip

\subheading{Example 3.10}

\noindent
For each~$m\in\N$ let $\Psi:\C^m\to\R^{2m}$ be the isomorphism
of real vector spaces which sends the standard basis vector $e_k\in\C^m$
to $e_k\in\R^{2m}$ and $ie_k\in\C^m$ to $e_{m+k}\in\R^{2m}$ for each
$k=1,\,\ldots,m$.
Consider the injective Lie algebra homomorphism 
$\phi:\su(m)\ni X\mapsto\Psi X\Psi\inv\in\so(2m)$
and the associated homomorphic embeddings
$\Phi:\SU(m)\to\SO(2m)$ and $\tilde\Phi:\SU(m)\to\Spin(2m)$.
Let $K_1:=\Imm\Phi\subset\SO(2m)=:G_1$ and $\tilde K_1:=\Imm\tilde\Phi\subset
\Spin(2m)=:\tilde G_1$\,.
We define $G:=\SO(2m+2)$, $\tilde G:=\Spin(2m+2)$, $G_2:=\SO(2)$,
$\tilde G_2:=\Spin(2)$, and consider the canonical embeddings
$G_1\times G_2\hookrightarrow G$ and $\tilde G_1\times
\tilde G_2\hookrightarrow
\tilde G$ which allow us to consider $G_1$ and~$G_2$ (respectively $\tilde
G_1\times\tilde G_2$) as commuting subgroups of~$G$ (respectively~$\tilde G$)
which are orthogonal with respect to the Killing metric.
Denote by $\k_1\,,\g_1\,,\g$ the Lie algebras of $K_1\,,G_1$\,, and~$G$,
respectively.

Let $J:=\left(\smallmatrix 0&-I\\I&0\endsmallmatrix\right)\in\so(2m)=\g_1$\,,
where $I$ denotes the $m$@-dimensional unit matrix. Note that $J$
commutes with $\Imm\phi=\k_1$\,.
We define two-dimensional tori~$H,\tilde H$ in~$G,\tilde G$ by
$$\align H&:=\exp(\R J)\times G_2\subset G_1\times G_2\subset G,\\
    \tilde H&:=\widetilde\exp(\R J)\times\tilde G_2\subset\tilde G_1\times
    \tilde G_2\subset\tilde G.\endalign
$$
For the Lie algebra~$\h$ of~$H$ (resp.~$\tilde H$) we then have
$$[\k_1\,,\h]=0\quad\text{ and }\quad\k_1\perp_{g_0}\h,\tag{10}
$$
where $g_0$ is any bi-invariant metric on~$G$ (resp.~$\tilde G$)
and $\perp_{g_0}$ denotes orthogonality with respect to~$g_0$\,.

Now assume $m\ge3$, and let $\{\widehat j(t)\}_{t\in(-\eps,\eps)}$
be a family of linear maps from~$\h$ to~$\su(m)$ satisfying conditions
1.)~and~2.)~from Proposition~3.6(i).
Let
$$j(t):=\phi\circ\widehat j(t):\h\to\so(2m)=\g_1\subset\so(2m+2)=\g.
$$
{}From condition~1.)~on the~$\widehat j(t)$ and the fact that $\Phi$
and~$\tilde \Phi$ are homomorphisms with $\Phi_*=\tilde\Phi_*=\phi$,
it follows that for each $Z\in\h$ we have
$j_Z(t)\in\Ad_{K_1}(j_Z(0))=\Ad_{\tilde K_1}(j_Z(0))$ for all~$t$.
We interpret the~$j(t)$ as linear maps from~$\h$ to~$\g$
and define $\la(t):={}^Tj(t):\g\to\h$ as the transpose of~$j(t)$
with respect to~$g_0$\,.
{}From~\thetag{10} we conclude, exactly as in Example~3.7\,/\,3.8, that
the~$\la(t)$ pairwise satisfy the conditions of Proposition~3.1 applied
to $(G,g_0)$, resp.~to $(\tilde G,g_0)$.
We thus obtain isospectral families
$$(\SO(2m+2),g_{\la(t)})\quad\text{ and }\quad(\Spin(2m+2),g_{\la(t)})
$$
for all $2m+2\ge2\cdot 3+2=8$, where $g_{\la(t)}$ is the left invariant
metric associated with~$\la(t)$ and~$g_0$ as in Proposition~3.1.
In particular, for $m=3$ we get continuous families of left invariant
isospectral metrics on $\SO(8)$, resp.~on $\Spin(8)$.

Finally note that condition~2.)~on the $\widehat j(t)$ implies that also
$\tr\bigl((j_{Z_1}(t)^2+j_{Z_2}(t)^2)^2\bigr)$ is nonconstant in~$t$.
In fact,
$\tr\bigl((j_{Z_1}(t)^2+j_{Z_2}(t)^2)^2\bigr)=
2\,\tr\bigl((\widehat j_{Z_1}(t)^2+\widehat j_{Z_2}(t)^2)^2\bigr)$, which is 
nonconstant in~$t$ by condition~2.).

Once more, the Ricci tensors of these manifolds have different norms
(Proposition~3.15), and the manifolds are not pairwise isospectral for the
Laplace operator acting on $1$@-forms (Corollary~3.17).

\bigskip

\subheading{Remark 3.11}
In the context of Example~3.10, we obtain an explicit example of the
data~$\h\subset \so(8)$ and $j(t):\h\to\k_1\subset\so(8)$ by using the
specific family of isospectral linear maps from~$\R^2$ to~$\su(3)$ which
was given in Proposition~3.6(ii):
$$j_{Z_1}(t)=\left(\smallmatrix &&&1&0&0&&\\
  &&&0&0&0&&\\&&&0&0&-1&&\\-1&0&0&&&&&\\0&0&0&&&&&\\0&0&1&&&&&\\&&&&&&0&0\\
  &&&&&&0&0\endsmallmatrix\right),\qquad
  j_{Z_2}(t)=\left(\smallmatrix 0&t&f(t)&&&&&\\-t&0&t&&&&&\\
  -f(t)&-t&0&&&&&\\&&&0&t&f(t)&&\\&&&-1&0&t&&\\
  &&&-f(t)&-t&0&&\\&&&&&&0&0\\&&&&&&0&0\endsmallmatrix\right),
$$
where all missing entries are zero, $f(t)=\sqrt{1-2t^2}$,
$t\in[-1/\sqrt2,1/\sqrt2]$, and
$\h=\spann\{Z_1\,,Z_2\}\subset\so(8)$ with
$$Z_1=\frac1{\sqrt3}\left(\smallmatrix &&&-1&0&0&&\\&&&0&-1&0&&\\&&&0&0&-1&&
  \\1&0&0&&&&&\\0&1&0&&&&&\\0&0&1&&&&&\\&&&&&&0&0
  \\&&&&&&0&0\endsmallmatrix\right),\qquad
  Z_2=\left(\smallmatrix 0&0&0&0&0&0&&\\0&0&0&0&0&0&&\\0&0&0&0&0&0&&\\
  0&0&0&0&0&0&&\\0&0&0&0&0&0&&\\0&0&0&0&0&0&&\\&&&&&&0&-1\\
  &&&&&&1&0\endsmallmatrix\right).
$$

\bigskip

\subheading{Remarks 3.12}

(i) In all our examples
(3.3, 3.7, 3.8, 3.9, and 3.10), the key tool were continuous families of
linear maps $j(t):\R^2\to\so(m)$, resp\. $j(t):\R^2\to\su(m)$, satisfying
conditions 1.)~and~2.)~from Example~3.3, resp\. from Proposition~3.6(i).
Note that both conditions are scaling invariant; that is, if the family
$t\mapsto j(t)$ satisfies them, then so does the family $t\mapsto\alpha j(t)$
for each $\alpha>0$. Rescaling the family $j(t)$ in any of our above examples
is equivalent to rescaling $\la(t)={}^Tj(t)$ by the same factor.
Note that for $\alpha\to0$, the isospectral families of metrics~$t\mapsto
g_{\alpha\la(t)}$ collapse
to the trivial family $g_{0\cdot\la(t)}\equiv g_0$\,, where $g_0$ is the
chosen bi-invariant metric. We conclude that continuous families of
isospectral, locally non-isometric left invariant metrics occur in fact
arbitrarily close to any fixed bi-invariant metric.

(ii) If a metric on a semisimple compact Lie group is sufficiently close
to a bi-invariant metric then it is of positive Ricci curvature.
By the argument in~(i) we thus obtain isospectral families of left invariant
metrics of positive Ricci curvature on $\SO(n\ge8)$ and $\SU(n\ge6)$.
These are the first examples of continuous families of isospectral
manifolds of positive Ricci curvature.
(However, in all of these isospectral
families the metrics are of mixed sectional curvature.)

\bigskip

\smallskip

\subhead 3.3 Ricci curvature and $1$@-form heat invariants\endsubhead
\bigskip

\noindent
All examples of families of isospectral left invariant metrics given in
Section~3.1 were applications of Proposition~3.1 (some of them via the more
special Corollary~3.2).
In this section we compute the Ricci curvature of the left invariant metrics
of the type occurring in Proposition~3.1 (see Lemma~3.18) and establish an
algebraic criterion which decides whether for a pair of isospectral left
invariant metrics arising from Proposition~3.1 the associated Ricci
curvatures have different norms (Theorem~3.14). In particular, it turns out
that in all the isospectral families from Section~3.2 the norm of the Ricci
tensor varies during the deformation (Proposition~3.15). This implies
not only that the manifolds are not pairwise locally isometric, but also,
as can be seen by using heat invariants, that they are not isospectral for
the Laplace operator acting on $1$@-forms (Corollary~3.17).

It be should mentioned that the scalar curvature can of course {\it not\/}
be used here to distinguish between the metrics: Since the manifolds are
homogeneous, the associated scalar curvature is constant on each of them;
the fact that volume and total scalar curvature are heat invariants thus
implies that this constant is the same for all metrics in the isospectral
family.

\bigskip

We fix certain objects and notations which we will use throughout
this section.

\bigskip

\subheading{Notation 3.13}
\newline
(i) Let $G$ be a compact Lie group with Lie algebra~$\g$ and a bi-invariant
metric~$g_0$\,. Let $H\subset G$ be a torus in~$G$ with Lie algebra
$\h\subset\g$, and denote by~$\uu$ the $g_0$@-orthogonal complement of the
centralizer~$\zh$ of~$\h$ in~$\g$.

(ii) We consider linear maps $\la:\g\to\g$ whose image is contained
in~$\h$ and which satisfy $\la\restr{\h\oplus\uu}=0$.
For any such~$\la$, we let $g\subla$ be the left invariant
metric on~$G$ which corresponds to the scalar product $(\Id+\la)^*g_0$
on~$\g$. We denote this scalar product on~$\g$, and the corresponding
ones on tensors, by $\scp\subla$\,, and we let $\Ric\supla$ be the
Ricci tensor associated with~$g\subla$\,. For $Z\in\h$ we denote by $\la_Z$
the $1$@-form $\<\la(\,.\,),Z\>_0\in\g^*$.
We define $j:\g\to\g$ as the transpose of~$\la$ with respect to $\scp_0$\,;
note that~$j$ vanishes on the $g_0$@-orthogonal complement
of~$\h$, and its image is contained in $\zh\cap\h^\perp$.

(iii) For any $X\in\g$ we write $\tilde X=(\Id-\la)(X)$.
Note that $\<\tilde X,\tilde Y\>\subla=\<X,Y\>_0$ for all $X,Y\in\g$.
Finally, we choose $g_0$@-orthonormal
bases $\{Z_1\,,\,\ldots,Z_r\}$ of~$\h\subset\g$ and
$\{V_1\,,\,\ldots,V_d\}$ of~$\g$.

\bigskip

\proclaim{Theorem 3.14}
Let $\la,\lap:\g\to\h\subset\g$ be two linear maps as above which
moreover satisfy condition~\thetag{$*3$} of Proposition~{\rm3.1}; i.e.,
for every $Z\in\h$ there exists $a_Z\in G$ such that
  $\lap_Z=\Ad_{a_Z}^*\la_Z$ and $\Ad_{a_Z}\restr{\h}=\Id\restr{\h}$\,.
Then we have, using the above notation:
$$\|\Ric\supla\|\subla^2-\|\Ric\suplap\|\sublap^2=\frac14
  \Bigl(\sum_{i,k=1}^r\tr\bigl((\ad_{j_{Z_i}})^2(\ad_{j_{Z_k}})^2\bigr)-
  \sum_{i,k=1}^r\tr\bigl((\ad_{j'_{Z_i}})^2(\ad_{j'_{Z_k}})^2\bigr)\Bigr).
  \tag{11}
$$
\endproclaim

\smallskip

We postpone the proof of Theorem~3.14 to the end of this section and first
deduce from it that the norm of the Ricci tensor varies indeed in all
the isospectral families of left invariant metrics given in the examples
in Section~3.2.

\bigskip

\proclaim{Proposition 3.15}
Let $g_{\la(t)}$ be any of the isospectral famlies of left invariant metrics
from Example {\rm3.3, 3.7, 3.8, 3.9,} or~{\rm3.10} of Section~{\rm3.2}.
Then $\|\Ric^{\la(t)}\|^2_{\la(t)}$ is nonconstant in~$t$.
\endproclaim

\smallskip

\demo{Proof}
First of all, note that in Theorem~3.14 the right hand side of~\thetag{11}
is zero if and only if the two sums running only
over pairs~$i,k$ with $i\ne k$ are equal.
In fact, $\tr\bigl((\ad_{j_{Z_i}})^4\bigr)
=\tr\bigl((\ad_{j'_{Z_i}})^4\bigr)$ for all~$i$ because $j_{Z_i}$ and
$j'_{Z_i}$ are conjugate by an automorphism of~$\g$ by
condition~\thetag{$*3$}.
Throughout Section~3.2 we worked with $\dimm\,\h=2$; hence we only need to
show that
$$\tr\bigl((\ad_{j_{Z_1}(t)})^2(\ad_{j_{Z_2}(t)})^2\bigr)
  \ne\const\text{ in~$t$}\tag{12}
$$
for each of the families $j(t)={}^T\la(t):\h\to\g$ from
the examples in Section~3.2.

Recall that in some of those examples, $\g$ was equal to a matrix algebra
$\m=\so(n)$ or $\m=\su(n)$ (Examples 3.8, 3.9, 3.10); in the remaining
examples~3.3 and~3.7, $\g$~was the direct sum of such an algebra~$\m$
with an abelian Lie algebra. In each case, the images of the maps
$j(t):\h\to\g$ were contained in~$\m$. In~\thetag{12} we can therefore
interpret ``$\ad$'' as the adjoint representation of~$\m$, and
``$\tr$'' as the trace over~$\m$.

We want to apply the formulas given in Lemma~3.16 below in order to
simplify~\thetag{12}. For this, we first recall that
the families $j(t):\h\to\m$ are isospectral (in the sense
of Definition~1.10 or Definition~3.5, respectively), which means that
for each $Z\in\h$ the $j_Z(t)$ are all conjugate to each other by elements
of $\O(n)$, resp\. $\SU(n)$; hence $\tr(j_Z(t))^2$ is constant in~$t$.
Moreover, $2\,\tr(j_{Z_1}(t)j_{Z_2}(t))$ equals the coefficient at
$su$ of $\tr(j_{sZ_1+uZ_2}(t)^2)$ and is thus constant in~$t$.
Finally,
$2\,\tr(j_{Z_1}(t)j_{Z_2}(t)j_{Z_1}(t)j_{Z_2}(t))+4\,\tr(j_{Z_1}(t)^2
j_{Z_2}(t)^2)$ is also constant in~$t$ since it equals the coefficient
at $s^2u^2$ of $\tr(j_{sZ_1+uZ_2}(t)^4)$.
Lemma 3.16 thus implies in our context:
$$\align\text{If }\m=\so(n)\text{ then }
  &\tr\bigl((\ad_{j_{Z_1}(t)})^2(\ad_{j_{Z_2}(t)})^2\bigr)=
  \const+(n-2)\tr\bigl(j_{Z_1}(t)^2j_{Z_2}(t)^2\bigr);\\
  \text{if }\m=\su(n)\text{ then }
  &\tr\bigl((\ad_{j_{Z_1}(t)})^2(\ad_{j_{Z_2}(t)})^2\bigr)=
  \const+2n\,\tr\bigl(j_{Z_1}(t)^2j_{Z_2}(t)^2\bigr),\endalign
$$
where ``const'' means constant in~$t$.
Note that $n-2\ne0$ since in all examples we needed $n>2$.
Thus in any of the isospectral families $j(t):\h\to\m$ from Section~3.2
we have that~\thetag{12} is equivalent to
$$\tr\bigl(j_{Z_1}(t)^2j_{Z_2}(t)^2\bigr)\ne\const\text{ in }t.\tag{13}
$$
But this was indeed always the case. In fact, in all our examples we
had $\tr\bigl((j_{Z_1}(t)^2+j_{Z_2}(t)^2)^2\bigr)\ne\const$ in~$t$
by condition~2.)~of Example~3.3 (which was assumed in Examples 3.3,
3.8, and shown to hold in Example~3.10), resp\. condition~2.)~of
Proposition~3.6(i) (which was assumed in Examples 3.7, 3.9).
This implies~\thetag{13} since $\tr(j_{Z_1}(t)^4)$ and $\tr(j_{Z_2}(t)^4)$
are constant in~$t$ by the isospectrality assumption.
\qed\enddemo

\bigskip

In the proof of Proposition~3.15 we have used the following formulas for
which we did not find a reference.

\bigskip

\proclaim{Lemma 3.16}
\roster
\item"(i)"
Let $X,Y\in\so(n)$ and $\ad$ be the adjoint representation
of~$\so(n)$ on itself. Then
$$\align &\qquad\tr((\ad_X)^2(\ad_Y)^2)\\&\qquad\qquad
  =(n-6)\tr(X^2Y^2)-2\,\tr(XYXY)+\tr(X^2)\tr(Y^2)
  +2(\tr(XY))^2.\endalign
$$
\item"(ii)"
Let $X,Y\in\su(n)$ and $\ad$ be the adjoint representation
of~$\su(n)$ on itself. Then
$$\tr((\ad_X)^2(\ad_Y)^2)=2n\,\tr(X^2Y^2)+2\,\tr(X^2)\tr(Y^2)
  +4(\tr(XY))^2.
$$
\endroster\endproclaim

\smallskip

\demo{Proof}
(i) Note that the adjoint representation of $\so(n)$ on itself is equivalent
to the canonical representation~$\rho$ of~$\so(n)$ on $\bigwedge^2\R^n$,
given by $\rho_X(y\wedge v)=Xy\wedge v+y\wedge Xv$.
Let $\{e_1\,,\,\ldots,e_n\}$ be the standard basis of~$\R^n$, and define
a scalar product on $\bigwedge^2\R^n$ by $\<y\wedge v,w\wedge z\>=
\<y,w\>\<v,z\>-\<y,z\>\<v,w\>$. Then
$$\tr(\rho_X^2\rho_Y^2)=\frac12\sum_{i,k=1}^n\<\rho_X^2\rho_Y^2(e_i\wedge
  e_k),(e_i\wedge e_k)\>.
$$
That this is indeed equal to the right hand side of the formula in~(i)
follows by straightforward calculation.

(ii) First consider the adjoint representation of the complex Lie algebra
$\gl(n,\C)$ of dimension~$n^2$ on itself, which is equivalent to the
canonical representation~$\rho$ of~$\gl(n,\C)$ on $(\C^n)^*\otimes\C^n$
given by $\rho_X(y^*v)=-({}^tXy)^*v+y^*(Xv)$.
Similarly as in~(i) we obtain by direct computation that
$$\align \tr(\rho_X^2\rho_Y^2)&=2n\,\tr(X^2Y^2)+2\,\tr(X^2)\tr(Y^2)
  +4(\tr(XY))^2\\
  &\qquad-4\,\tr(X)\tr(XY^2)-4\,\tr(Y)\tr(X^2Y)\endalign
$$
for all $X,Y\in\gl(n,\C)$. For $X,Y\in\su(n)$ the last two terms vanish.
Moreover, for $X,Y\in\su(n)$ the trace of $\rho_X^2\rho_Y^2$ on
$\gl(n,\C)$, interpreted now as the {\it real\/} Lie algebra
$\uu(n)\oplus i\uu(n)$ of dimension~$2n^2$, equals two times the right hand
side of the above
formula on the one hand, and two times the trace of $\ad_X^2\ad_Y^2$
on~$\uu(n)$ on the other hand.
Since $\uu(n)$ is the sum of its center (spanned by~$i\Id$) and~$\su(n)$,
the assertion of~(ii) now follows.
\qed\enddemo

\bigskip

\proclaim{Corollary 3.17}
In all examples given in Section~{\rm3.2}, the isospectral manifolds
$(G,g_{\la(t)})$ are not pairwise isospectral for the Laplace operator
acting on $1$@-forms.

More generally,
if $(M,g)$, $(M',g')$ is any pair of homogeneous manifolds which are
isospectral for the Laplace operator on functions and for which the
{\rm(}constant\/{\rm)} functions
$\|Ric^g\|_g^2$ and $\|Ric^{g'}\|_{g'}^2$ are nonequal,
then the associated Laplace operators on $1$@-forms are not isospectral.
\endproclaim

\bigskip

\demo{Proof}
We only need to prove the second statement since
by Proposition~3.15 the norm of the Ricci tensors associated with the
metrics~$g_{\la(t)}$ from Section~3.2 does change nontrivially as~$t$ varies.

For any Riemannian metric~$g$ on a closed Riemannian manifold~$M$ the
heat invariants for the associated Laplace operator on $p$@-forms
$(0\le p\le\dimm\,M)$ are the coefficients $a_i^p(g)$ occurring in the
asymptotic expansion
$$\tr\bigl(\exp(-s\Delta_g^p)\bigr)\;\sim\;(4\pi s)^{-\dimm M/2}
\sum_{i=0}^\infty a_i^p(g)s^i\quad\text{ for }\,s\searrow 0.
$$
By \cite{18}, Theorem 4.8.18 we have
$$\align a_0^0(g)&=\voll(M,g),\quad a_1^0(g)=\tsize\frac16\tsize\int_M
  \scal^g\dvol_g\,,\\
  a_2^0(g)&=\tsize\frac1{360}\tsize\int_M\bigl(5(\scal^g)^2-2\|\Ric^g\|_g^2
  +2\|R^g\|_g^2\bigr)\dvol_g\,,\\
  a_2^1(g)&=a_2^0(g)\cdot\dimm\,M-\tsize\frac1{12}\tsize\int_M\bigl(
  2(\scal^g)^2-6\|\Ric^g\|_g^2+\|R^g\|_g^2\bigr)\dvol_g\,,\endalign
$$
where $\scal^g$, $\Ric^g$, and $R^g$ denote the scalar curvature, Ricci
tensor, and curvature tensor associated with~$g$.
The first two of the above heat invariants imply that if $(M,g)$ and
$(M',g')$ are homogeneous and isospectral on functions, then their
(constant) scalar curvatures are the same; in particular, we then also
have $\int_M(\scal^g)^2\dvol_g=\int_{M'}(\scal^{g'})^2\dvol_{g'}$\,.
By $a_2^0(g)=a_2^0(g')$, the numbers $x:=\int_M\|\Ric^g\|^2_g\dvol_g-
\int_{M'}\|\Ric^{g'}\|_{g'}^2\dvol_{g'}$ and $y:=\int_M\|R^g\|_g^2\dvol_g
-\int_{M'}\|R^{g'}\|_{g'}^2\dvol_{g'}$ satisfy $-2x+2y=0$. If now in
addition the two manifolds were isospectral on $1$@-forms, then
$a_2^1(g)=a_2^1(g')$ and thus $-6x+y=0$. These two equations together imply
$x=y=0$; but $x=0$ contradicts our assumption.
\qed\enddemo

\bigskip

The rest of this section is devoted to the proof of Theorem~3.14.
We continue to use Notation~3.13; recall in particular that we consider
linear maps $\la:\g\to\g$ with image in~$\h$ and $\la\restr{\h\oplus\uu}=0$,
that $j={}^T\la$ with respect to~$g_0$\,, and that $\tilde X=X-\la(X)$ for
$X\in\g$.
First we need a formula for the Ricci tensor $\Ric\supla$ of $(G,g\subla)$.

\bigskip

\proclaim{Lemma 3.18}
For all $X\in\g$ we have
$$\Ric\supla(\tilde X,\tilde X)=\Ric^0((\Id+j)X,(\Id+j)X)-\<\ad_X\,,
  \la\circ\ad_X\>_0-\frac12\|\la\circ\ad_X\|_0^2\,.
$$
\endproclaim

\smallskip

\demo{Proof}
Using the general formula for the Ricci tensor of a homogeneous manifold
given in \cite{5}, Corollary 7.38, and the fact that $\{\tilde V_1\,,
\,\ldots,\tilde V_d\}$ is a $g\subla$@-orthonormal basis of~$\g$, we have
$$\Ric\supla(\tilde X,\tilde X)=-\frac12\|\ad_{\tilde X}\|\subla^2
  -\frac12\tr((\ad_{\tilde X})^2)+\frac14\sum_{i,k=1}^d\<\tilde X,
  [\tilde V_i\,,\tilde V_k]\>\subla^2\,.\tag{14}
$$
We will show that
$$\aligned -\frac12\|\ad_{\tilde X}\|\subla^2
  -\frac12\tr((\ad_{\tilde X})^2)=
  &-\<\ad_X\,,\la\circ\ad_X\>_0-
  \frac12\|\la\circ\ad_X\|_0^2\\
  &+\<\ad_X\,,\ad_X\circ\la\>_0-
  \frac12\|\ad_X\circ\la\|_0^2\endaligned\tag{15}
$$
and
$$\frac14\sum_{i,k=1}^d\<\tilde X,[\tilde V_i\,,\tilde V_k]\>\subla^2
  =\Ric^0((\Id+j)X,(\Id+j)X)
  -\<\ad_X\,,\ad_X\circ\la\>_0+\frac12\|\ad_X\circ\la\|_0^2\,.\tag{16}
$$
These two formulas, together with~\thetag{14}, will clearly imply our
statement.

Note that $\la[X,\la(Y)]=0$ for all $X,Y\in\g$ since $\ad_{\la(Y)}$
annihilates $\zh$ and thus has image in~$\uu$, on which~$\la$ vanishes.
Therefore
$$\align \|\ad_{\tilde X}\|\subla^2&=\sum_{i=1}^d
  \|[\tilde X,\tilde V_i]\|\subla^2
  =\sum_{i=1}^d
  \|[X,V_i]-[\la(X),V_i]-[X,\la(V_i)]+\la([X,V_i])\|_0^2\\
  &=\|\ad_X-\ad_{\la(X)}-\ad_X\circ\la+\la\circ\ad_X\|_0^2\,.
  \endalign
$$
Moreover,
$$\align -\tr((\ad_{\tilde X})^2)=\sum_{i=1}^d
  \|[\tilde X,V_i]\|_0^2
  &=\sum_{i=1}^d\|[X,V_i]-[\la(X),V_i]\|_0^2
  =\|\ad_X-\ad_{\la(X)}\|_0^2\,.\endalign
$$
Thus
$$-\|\ad_{\tilde X}\|\subla^2-\tr((\ad_{\tilde X})^2)=
  2\<\ad_X-\ad_{\la(X)}\,,\ad_X\circ\la-\la\circ\ad_X\>_0
  -\|\ad_X\circ\la-\la\circ\ad_X\|_0^2\,.
$$
Since $\ad_{\la(X)}$ annihilates~$\zh$ and has image in~$\uu$,
it is $g_0$@-orthogonal to both of $\ad_X\circ\la$ (which annihilates~$\uu
=\zh^\perp$) and $\la\circ\ad_X$ (whose image is contained in
$\h\subset\uu^\perp$).
Moreover, $\ad_X\circ\la$ has image in~$\uu$ and is therefore
$g_0$@-orthogonal to $\la\circ\ad_X$\,.
Formula~\thetag{15} now follows; it remains to show formula~\thetag{16}.
We have
$$\allowdisplaybreaks
  \align&\sum_{i,k=1}^d\<\tilde X,[\tilde V_i\,,\tilde V_k]\>\subla^2
  =\sum_{i,k=1}^d \<X,[V_i\,,V_k]-[V_i\,,\la(V_k)]+[V_k\,,\la(V_i)]
  +\la([V_i\,,V_k])\>_0^2\\
  &\qquad=\sum_{i,k=1}^d \bigl(\<X+j(X),[V_i\,,V_k]\>_0^2+\<X,[V_i\,,
  \la(V_k)]-[V_k\,,\la(V_i)]\>_0^2\\
  &\qquad\qquad -4\<X,[V_i\,,V_k]\>_0\<X,[V_i\,,\la(V_k)]\>_0
  -4\<X,\la([V_i\,,V_k])\>_0\<X,[V_i\,,\la(V_k)]\>_0\bigr)\\
  &\qquad=\sum_{k=1}^d \bigl(\|[X+j(X),V_k]\|_0^2+2\|[X,\la(V_k)]\|_0^2
  +2\<X,[V_k\,,\la([X,\la(V_k)])]\>_0^2\\
  &\qquad\qquad -4\<[X,V_k],[X,\la(V_k)]\>_0+4\<X,\la([[X,
  \la(V_k)],V_k])\>_0\bigr)\\
  &\qquad={\vphantom{\sum^d}}4\,\Ric^0
  (X+j(X),X+j(X))+2\|\ad_X\circ\la\|_0^2+0-4\<\ad_X\,,
  \ad_X\circ\la\>_0+0.\endalign
$$
\smallskip\noindent
Here, the third term is zero because of $\la([X,\la(\,.\,)])=0$
(see above); to see that the fifth term is zero, we assume the
$g_0$@-orthonormal basis $\{V_1\,,\,\ldots,V_d\}$ to be adapted to the
$g_0$@-orthogonal decomposition $\g=(\zh\cap\h^\perp)\oplus(\h\oplus\uu)$.
The term $\la([[X,\la(V_k)],V_k])$ vanishes for $V_k\in\h\oplus\uu\subseteq
\kernn\la$; but for $V_k\in\zh\cap\h^\perp$ (which is obviously a Lie
subalgebra of~$\g$), we have $[\uu,V_k]\subseteq\uu$ and thus
$\la([[X,\la(V_k)],V_k])=0$ (recall that $\ad_X\circ\la$ has image in
$\uu\subset\kernn\la$). Formula~\thetag{16} now follows.
\qed\enddemo

\smallskip

\subheading{Proof of Theorem 3.14}
\newline
For all $X,Y\in\g$ we have
$$\align\<\ad_X\,,\la\circ\ad_Y\>_0&=\sum_{k=1}^d \<[X,V_k],\la([Y,V_k])\>_0
  =\sum_{k=1}^d\sum_{i=1}^r \<[X,V_k],Z_i\>_0\<j_{Z_i}\,,[Y,V_k]\>_0\\
  &=\sum_{i=1}^r\<[X,Z_i],[Y,j_{Z_i}]\>_0=-\sum_{i=1}^r\<\ad_{j_{Z_i}}
  \ad_{Z_i}X,Y\>_0\endalign
$$
which is, in particular, symmetric in $X$ and~$Y$ since $j_{Z_i}\in\zh$
commutes with $Z_i\in\h$.
Similarly,
$$\align \<\la\circ\ad_X\,,\la\circ\ad_Y\>_0&=\sum_{k=1}^d\<\la([X,V_k]),
  \la([Y,V_k])\>_0\\ &\hskip-20pt
  =\sum_{k=1}^d\sum_{i=1}^r\<j_{Z_i}\,,[X,V_k]\>_0
  \<j_{Z_i}\,,[Y,V_k]\>_0
  =-\sum_{i=1}^r\<(\ad_{j_{Z_i}})^2X,Y\>_0\,.\endalign
$$
We can thus reformulate Lemma~3.18 as
$$\aligned\Ric\supla(\tilde X,\tilde Y)&=\<\Ric^0(\Id+j)X,(\Id+j)Y\>_0
  +\sum_{i=1}^r\<\ad_{Z_i}\ad_{j_{Z_i}}X,Y\>_0\\
  &\quad+\frac12
  \sum_{i=1}^r\<(\ad_{j_{Z_i}})^2X,Y\>_0\,.\endaligned\tag{17}
$$
For any given pair $Z,W\in\h$ there exist, by our assumption on~$\la$
and~$\lap$,
elements $a_{s,u}\in G$ such that $j'_{sZ+uW}=\Ad_{a_{s,u}}j_{sZ+uW}$
and $\Ad_{a_{s,u}}\restr{\h}=\Id\restr{\h}$ for all $s,u\in\R$.
Note that the endomorphism $\sum_{i=1}^d \ad_{V_i}^2=-4\,\Ric^0$ of~$\g$
restricts to a scalar multiple of the Casimir operator of the adjoint
action on each irreducible
component of~$\g$ and thus commutes with every inner automorphism.
Consequently,
$$\align &\<\Ric^0(\Id+j')(sZ+uW),(\Id+j')(sZ+uW)\>_0\\&\qquad\quad=
  \<\Ad_{a_{s,u}}\Ric^0(\Id+j)(sZ+uW),\Ad_{a_{s,u}}(\Id+j)(sZ+uW)\>_0\\
  &\qquad\quad=\<\Ric^0(\Id+j)(sZ+uW),(\Id+j)(sZ+uW)\>_0\endalign
$$
for all $s,u\in\R$; comparing the coefficients at~$su$, we obtain that
$\<\Ric^0(\Id+j')Z,\allowbreak(\Id+j')W\>_0=\<\Ric^0(\Id+j)Z,(\Id+j)W\>_0$
for all $Z,W\in\h$. Since the $j_{Z_i}$ commute with~$\h$,
we have by~\thetag{17} that
these expressions equal $\Ric\supla(Z,W)$ and $\Ric\suplap(Z,W)$,
respectively. We conclude
$$
\Ric\supla\drestr{\h\times\h}=\Ric\suplap\drestr{\h\times\h}\,.\tag{18}
$$
For $Z\in\h$ and $X\in\h^\perp$ we have $j_X=j'_X=0$ and we obtain, similarly
as above, $\<\Ric^0(\Id+j')Z,X\>_0\allowbreak=\<\Ad_{a_Z}\Ric^0(\Id+j)Z,X\>_0
\allowbreak=\<\Ric^0(\Id+j)Z,\Ad_{a_Z}\inv X\>_0$\,.
Since $\Ad_{a_Z}$ preserves~$\h^\perp$, this implies by~\thetag{17} that
$$\|\Ric\supla\drestr{\h\times\h^{\perp\subla}}\|\subla^2
  =\|\Ric\suplap\drestr{\h\times\h^{\perp\sublap}}\|\sublap^2\,,\tag{19}
$$
where $\h^{\perp\subla}$ denotes the $g\subla$@-orthogonal complement
of~$\h$ in~$\g$. Note also that
$$\align &\sum_{i=1}^r\<\Ric^0,\ad_{Z_i}\ad_{j'_{Z_i}}+\tfrac12
  (\ad_{j'_{Z_i}})^2\>_0\\
  &\qquad\quad=\sum_{i=1}^r\<\Ric^0,\Ad_{a_{Z_i}}(\ad_{Z_i}
  \ad_{j_{Z_i}}
  +\tfrac12(\ad_{j_{Z_i}})^2)\Ad\inv_{a_{Z_i}}\>_0\\
  &\qquad\quad=\sum_{i=1}^r\<\Ric^0,\ad_{Z_i}\ad_{j_{Z_i}}
  +\tfrac12(\ad_{j_{Z_i}})^2\>_0.\endalign
$$
Moreover, both of $\ad_{Z_i}$ and $\ad_{j_{Z_i}}$ vanish on~$\h$ and have
image in~$\h^\perp$ because of $Z_i\,,j_{Z_i}\in\zh$. Therefore, again
by~\thetag{17}:
$$\aligned
  \|\Ric^{g\subla}\drestr{\h^{\perp\subla}\times\h^{\perp\subla}}\|\subla^2
  &=\|\Ric^0\drestr{\h^\perp\times\h^\perp}\|_0^2+2\sum_{i=1}^r
  \<\Ric^0,\ad_{Z_i}\ad_{j_{Z_i}}+\tfrac12(\ad_{j_{Z_i}})^2\>_0\\
  &\qquad +\|\sum_{i=1}^r \bigl(\ad_{Z_i}\ad_{j_{Z_i}}+\tfrac12
  (\ad_{j_{Z_i}})^2\bigr)\|_0^2\,.\endaligned\tag{20}
$$
As we just saw, only the third summand in this expression might differ
from the corresponding one for~$\lap$.
We have
$$\aligned &\|\sum_{i=1}^r \bigl(\ad_{Z_i}\ad_{j_{Z_i}}+\tfrac12
  (\ad_{j_{Z_i}})^2\bigr)\|_0^2\\
  &\qquad\quad=\sum_{i,k=1}^r\tr(\ad_{Z_i}\ad_{j_{Z_i}}
  \ad_{Z_k}\ad_{j_{Z_k}})
  +\sum_{i,k=1}^r\tr(\ad_{Z_i}\ad_{j_{Z_i}}(\ad_{j_{Z_k}})^2)\\
  &\qquad\qquad\quad+\frac14\sum_{i,k=1}^r\tr((\ad_{j_{Z_i}})^2
  (\ad_{j_{Z_k}})^2).\endaligned\tag{21}
$$
Noting again that the $Z_i$ commute with the $j_{Z_k}$, we see that
$\tr(\ad_{Z_i}\ad_{Z_k}\ad_{j'_{Z_i}}\ad_{j'_{Z_k}})$ equals half the
coefficient at~$su$ of
$$\align \tr(\ad_{Z_i}\ad_{Z_k}(\ad_{j'_{sZ_i+uZ_k}})^2)
    &=\tr(\Ad_{a_{s,u}}\ad_{Z_i}\ad_{Z_k}(\ad_{j_{sZ_i+uZ_k}})^2
    \Ad_{a_{s,u}}\inv)\\
    &=\tr(\ad_{Z_i}\ad_{Z_k}(\ad_{j_{sZ_i+uZ_k}})^2).\endalign
$$
Similarly, $\tr(\ad_{Z_i}\ad_{j'_{Z_i}}(\ad_{j'_{Z_k}})^2)$ equals
one third of the coefficient at~$su^2$ of
$$\align \tr(\ad_{Z_i}(\ad_{j'_{sZ_i+uZ_k}})^3)
    &=\tr(\Ad_{a_{s,u}}\ad_{Z_i}(\ad_{j_{sZ_i+uZ_k}})^3
    \Ad_{a_{s,u}}\inv)\\
    &=\tr(\ad_{Z_i}(\ad_{j_{sZ_i+uZ_k}})^3).\endalign
$$
This, together with the formulas~\thetag{18}--\thetag{21},
implies the statement of the theorem.
\qed

\bigskip

\smallskip

\subhead 3.4 Infinitesimal spectral rigidity of bi-invariant
  metrics\endsubhead
\bigskip

\noindent
In Section~3.2 we used Proposition~3.1 to construct many examples
of continuous families $g_{\la(t)}$ of left invariant, isospectral,
non-isometric metrics on compact Lie groups~$G$.
As we saw, these families can occur arbitrarily close to
a bi-invariant metric~$g_0$ on~$G$ (Remark~3.12(i)).
Obviously, however, our construction never yields any isospectral
deformations containing $g_0$ itself. In fact, if $\la=0$ and $\la,\lap$
satisfy condition~\thetag{$*3$} of Proposition~3.1, then also $\lap=0$.

A natural question to ask in this context is whether nontrivial,
continuous isospectral deformations of bi-invariant metrics within
the class of left invariant metrics, even though not available by
our construction, might exist nevertheless. By the following theorem
the answer to this question is no.

\bigskip

\proclaim{Theorem 3.19}
Let $G$ be a compact Lie group and $g_0$ be a bi-invariant metric
on~$G$. Let $\eta>0$ and $\{g(t)\}\restr{t\in(-\eta,\eta)}$ be a
continuous family of left invariant metrics on~$G$ such that
$g(0)=g_0$\,. If the metrics $g(t)$ are pairwise isospectral,
then $g(t)\equiv g_0$ for all~$t$.
\endproclaim

\smallskip

\demo{Proof}
As in the ``algebraic proof'' of Proposition~3.1, denote the right-regular
unitary representation of~$G$ on $L^2(G)$ by~$\rho$.
Let $U\subseteq L^2(G)$ be a linear subspace which is invariant
under~$\rho$ and irreducible.
Since $G$ is compact and the action is unitary, any such~$U$ is
finite dimensional. If $g$ is any left invariant metric on~$G$
then we have by~\thetag{7} that
$\Delta_g=-\sum_{i=1}^d X_i^2=-\sum_{i=1}^d(\rho_*X_i)^2$,
where $\{X_1\,,\,\ldots,X_d\}$ is a left invariant $g$@-orthonormal frame.
Since~$U$ is invariant under~$\rho$ it is also invariant under~$\Delta_g$
and thus spanned by eigenfunctions; in particular $U\subset C^\infty(G)$
because~$U$ is finite dimensional.

Now consider our family $\Delta_{g(t)}$ restricted to~$U$.
For every~$t$, $\spec(\Delta_{g(t)}\restr U)$ is contained
in the discrete set $\spec(G,g(t))$ which is independent of~$t$
by assumption.
Since $g(t)$ and therefore $\spec(\Delta_{g(t)}\restr U)$
depends continuously on~$t$, $\spec(\Delta_{g(t)}\restr U)$
must be independent of~$t$, too.

The metric $g(0)=g_0$ is bi-invariant, hence $\Delta_{g(0)}$ commutes
with right translations (which are $g_0$@-isometries) and therefore
with the representation~$\rho$.
Since $U$ is irreducible, it follows from Schur's lemma
that $\Delta_{g(0)} \restr U$ is a multiple of the identity.
Now $\spec(\Delta_{g(t)}\restr U)=\spec(\Delta_{g(0)}\restr U)$
implies that $\Delta_{g(t)}\restr U$ too is a multiple of the
identity and equals $\Delta_{g(0)}\restr U$\,.

Recall that $L^2(G)$ is a sum of invariant, irreducible subspaces
such as the one we just considered. Since $\Delta_{g(t)}$ and
$\Delta_{g(0)}$ coincide on each of these, they coincide completely
on $C^\infty(G)$.
But if two Riemannian metrics on a manifold have equal Laplacians
then they are theirselves equal.
Thus $g(t)\equiv g(0)$ for all~$t$, as claimed.
\qed\enddemo

\bigskip

\bigskip
\noindent
\gross 4.\hbox to6pt{}Conformally equivalent manifolds which are
  isospectral and\newline
\hbox to17pt{}not locally isometric\rm

\bigskip

\noindent
In the first section of this chapter we present a canonical generalization
of Theorem~1.6\,/\,Proposition~1.8 from Chapter~1; see
Theorem~4.3\,/\,Proposition~4.5. Here the fibers of the torus bundles
under consideration are in general no longer totally geodesic.
In particular, we hereby leave the context of Theorem~1.3 which was
our general starting point in Chapter~1. However, Theorem~4.3 and
Proposition~4.5 can be regarded as special versions of another theorem
which was established recently by C.~Gordon and Z.~Szab\'o in~\cite{23};
see Remark~4.4(ii) below.

In Section~4.2 we use Proposition~4.5 to construct the first pairs of
isospectral manifolds which are conformally equivalent and not locally
isometric. Note that in exactly one instance there have previously
been examples (even continuous families) of isospectral, conformally
equivalent manifolds; namely, those constructed in~1990 by R.~Brooks
and C.~Gordon~\cite{6}. However, the manifolds in these isospectral
families had (as all examples of isospectral manifolds known at that time)
a common Riemannian covering and thus were locally isometric.

For proving that our new examples are not locally isometric we use results
from Chapter~3 to show that the preimages of the maximal scalar curvature
on the two manifolds are not locally isometric because their Ricci
tensors (associated with the induced metrics) have different norms;
see Proposition~4.7.

\bigskip

\smallskip

\subhead 4.1\ \ Isospectral torus bundles whose fibers are not totally
  geodesic\endsubhead
\nopagebreak
\bigskip
\nopagebreak
\subheading{Notation 4.1}
Let $H$ be a torus with Lie algebra $\h=T_eH$, and let $H$ be equipped
with a fixed invariant metric. Let $M$ be a principal $H$@-bundle over
a closed Riemannian manifold $(N,h)$, and let $\phi,\psi\in
C^\infty(N,\R_+)$.
\roster
\item"(i)"
Given a connection form~$\om$ on~$M$
we denote by~$g\suboii$ the unique $H$@-invariant
Riemannian metric on~$M$ which satisfies:\newline
\qquad\qquad 1.)
For each $p\in N$, the induced metric ${g\suboii}\restr{\pi\inv_H(p)}$
on the fiber over~$p$ equals $\phi(p)$ times the given invariant metric
(induced from the metric on~$H$).\newline
\qquad\qquad 2.) The projection $\pi_H:M\to N$ is a Riemannian submersion
with respect to $g\suboii$ on~$M$ and $\psi h$ on~$N$.\newline
\qquad\qquad 3.) The $\om$@-horizontal distribution $\kernn\om$ is
$g\suboii$@-orthogonal to the fibers.
\newline
In particular, note that $g\suboii(X,Z)=\tilde\phi(x)\om_Z(X)$ for
all $X\in T_xM$ and $Z\in\h$, where $\tilde\phi$ denotes
the lift of~$\phi$ to~$M$.
\item"(ii)"
If $M$ is the trivial bundle $N\times H$ and $\la$ is an $\h$@-valued
$1$@-form on~$N$ we write $g\sublii:=g\suboii$\,, where $\om$ is the
connection form on~$N\times H$ defined by $\om(X,Z)=\la(X)+Z$ for
all $(X,Z)\in T(N\times H)\cong TN\times\h$.

\endroster

\bigskip

\subheading{Remark 4.2}
The metric $g\subom$ which was defined in Notation~1.5(iii) is just the same
as $g_{\om,1,1}$\,. In other words, introducing the metric $g\suboii$
on~$M$ can be described as first introducing~$g\subom$ and then stretching
vertical vectors by $\tilde\phi^{1/2}$ and horizontal vectors
by $\tilde\psi^{1/2}$, where $\tilde\phi$ and~$\tilde\psi$ are the
lifts of $\phi$ and~$\psi$ to~$M$.
In the context of~4.1(ii), introducing~$g\sublii$ can be described
analogously, this time by first introducing the metric~$g\subla$ on
$N\times H$ which was defined in Notation~1.5(iv).

Note also that for $\psi:=\phi$, the metric $g_{\om,\phi,\phi}$ equals
$\tilde\phi g\subom$ which is conformally equivalent to~$g\subom$\,.
Similarly we have $\;g_{\la,\phi,\phi}=\tilde\phi g\subla$ in the context
of~4.1(ii).
\bigskip

\proclaim{Theorem 4.3}
Let $(N,h)$ be a closed Riemannian manifold and $H$ be a torus equipped
with an invariant metric. Let $M$ be a principal $H$@-bundle over
$(N,h)$, let $\om,\omp$ be two connection forms on~$M$, and let
$\phi,\phip,\psi,\psip\in C^\infty(N,\R_+)$. Assume:
\roster
\item"($*5$)" For every $Z\in\h$ there exists a bundle automorphism
  $F_Z:M\to M$ which induces an isometry~$\fzb$ on the base
  manifold $(N,h)$ and satisfies\newline 
  $\omp_Z=F_Z^*\om_Z$\,,
  $\;\phip=\fzb^*\phi$, \;and $\;\psip=\fzb^*\psi$.
\endroster
Then $(M,g\suboii)$ and $(M,g\suboiip)$ are isospectral.
\endproclaim

\smallskip

\demo{Proof}
In the following we write $g=g\suboii$ and $g'=g\suboiip$\,.
Let $\HH=L^2(M,g)=L^2(M,g')$.
For any closed connected subgroup~$W\subset H$ of codimension~$1$ we denote,
as in the proof of Theorem~1.3, by $\HH_W$ the space of $W$@-invariant
functions in~$\HH$. Let $\CC_W:=C^\infty(M)\cap\HH_W$\,. Finally we denote
by $\HH_0$ the space of $H$@-invariant functions in~$\HH$ and let
$\CC_0:=C^\infty(M)\cap\HH_0$\,. We claim that
$${\Delta_{g'}\vphantom)}\restr{\CC_W}
  =(F_Z^*\circ\Delta_g\circ{F_Z^*}\inv)\restr{\CC_W}\,,\tag{22}
$$
where $Z\in\h\minzero$ is chosen orthogonal to $T_eW$, and $F_Z$ is as
in~\thetag{$*5$}.
Note that~$F_Z$\,, being a bundle automorphism, leaves the spaces $\CC_W$
and $\CC_0\subset\CC_W$ invariant.
Therefore equation~\thetag{22} implies $\spec(\Delta_g\restr{\HH_0})
=\spec(\Delta_{g'}\restr{\HH_0})$ and $\spec(\Delta_g\restr{\HH_W})
=\spec(\Delta_{g'}\restr{\HH_W})$. Since $W$ was arbitrary, we will then
be done by the decomposition~\thetag{1} from the proof of Theorem~1.3.

It remains to prove~\thetag{22}.
In contrast to the situation of Theorem~1.6, the fibers of the Riemannian
submersions $\pi_H:(M,g)\to(N,\psi h)$ and $\pi_H:(M,g')\to(N,\psip h)$
are in general not totally geodesic now (unless $\phi$ is constant).
The proof of equation~\thetag{22} has to take into account the mean
curvature vector fields~$V$ (on $(M,g)$) and~$V'$ (on $(M,g')$) of the
$H$@-orbits in~$M$. We first compute $V$ and~$V'$.
Let $\nabla$ be the Levi-Civit\`a connection of~$g$, and denote by
$\tilde\phi, \tilde\psi$ the lifts of $\phi,\psi$ to~$M$.
For any $Z\in\h$ the vector field $\nabla_Z Z$ is obviously $g$@-orthogonal
to the $H$@-orbits and thus $\om$@-horizontal. For any $\om$@-horizontal,
$H$@-invariant vector field~$X$ on~$M$ we have, noting that $X$ commutes
with~$Z$:
$$\align\psi h(\pi_{H*}(\nabla_Z Z),\pi_{H*}X)\;&=\;g(\nabla_Z Z,X)
  \;=\;Z(g(Z,X))-g(Z,\nabla_ZX)\\ 
  &=\; Z(\tilde\phi g\subom(Z,X))-g(Z,\nabla_XZ)\\
  &=\; 0-\tfrac12X(\tilde\phi g\subom(Z,Z))
  \;=\;-\tfrac12|Z|^2X(\tilde\phi)\\ 
  &=\;-\tfrac12|Z|^2\psi h(\tfrac1\psi
  \grad_h\phi,\pi_{H*}X).\endalign
$$
Letting $Z$ run through $g$@-orthonormal bases of the tangent spaces
to the $H$@-orbits and summing up, we get
$$V\;=\;\text{$\om$@-horizontal lift of }\;\frac{-\dimm\,\h}{2\phi\psi}
  \,\grad_h\phi.
$$
Analogously, $V'$ is the $\omp$@-horizontal lift of $\;\frac{-\dimm\h}
{2\phip\psip}\,\grad_h\phip$\,.

Now let $x\in M$ and $p=\pi_H(x)\in N$.
Choose a local frame $\{E_1\,,\,\ldots,E_n\}$ on a neighbourhood~$U$
of~$p$ such that $\{E_1(p),\,\ldots,E_n(p)\}$ is a $\psip h$@-orthonormal
basis of $T_pN$, and such that the integral curves of the~$E_i$ through~$p$
are geodesics in $(N,\psip h)$.
Denote the $\omp$@-horizontal lift of~$E_i$ to $\pi_H\inv(U)\subseteq M$
by~$X_i$\,. Since $\pi_H:(M,g')\to(N,\psip h)$ is a Riemannian submersion,
the integral curves of~$X_i$ through~$x$ are geodesics in $(M,g')$.
Thus
$$\Delta_{g'}\restr x\;=\;-\sum_{i=1}^n X_i\restr x X_i +
  \tilde\phip(x)\inv{\Delta_\h}
  \restr x + V'\restr x\;,
$$
where $\Delta_\h:=-\sum_{k=1}^r Z_k^2$ and $\{Z_1\,,\,\ldots,Z_r\}$
is an orthonormal basis of~$\h$.
Now let $y:=F_Z(x)$ and $Y_i:=F_{Z*}(X_i)$. Then the~$Y_i$ are $H$@-invariant
vector fields defined in an $H$@-invariant neighbourhood of~$y$.
Since $\omp_Z=F_Z^*\om_Z$\,, $\;\phip=\bar F_Z^*\phi$,
$\psip=\bar F_Z^*\psi$, and~$\bar F_Z$ is an isometry of $(N,h)$, the vector
field $F_{Z*}V'$ equals~$V$ up to errors tangent to the $W$@-orbits.
Moreover,
each~$Y_i$ is $\om$@-horizontal up to errors tangent to the $W$@-orbits.
We write $Y_i=A_i+U_i$\,, where $A_i$ is $\om$@-horizontal and $U_i$
is tangent to the $W$@-orbits. Note that $A_i$ and $U_i$ are again
$H$@-invariant, and $[A_i\,,U_i]$ is tangent to the $W$@-orbits.
Hence for $f\in\CC_W$ we have
$$V\restr x (f) = V'\restr x(F_Z^*f)\quad\text{and}\quad
  A_i\restr y A_i(f) = Y_i\restr y Y_i (f) = X_i\restr x X_i (F_Z^*f).
$$
Since $F_Z$ induces an isometry from $(N,\psip h)$ to $(N,\psi h)$
the~$A_i\restr y$ are $g$@-orthonormal, and thus
$${\Delta_g}\restr y = -\sum_{i=1}^n A_i\restr y A_i +
  \tilde\phi(y)\inv{\Delta_\h}\restr y + V\restr y\,.
$$
Therefore we have indeed
$$\align(\Delta_g f)(y) &\;=\;\Bigl(-\sum_{i=1}^n X_i^2(F_Z^*f)+\tilde
  \phi^{\prime-1}\Delta_\h(F_Z^*f) + V'(F_Z^*f)\Bigr)(x)\\
  &\;=\;\big(\Delta_{g'}(F_Z^*f)\big)(x)\;=\;\big(({F_Z^*}\inv
  \circ\Delta_{g'}\circ F_Z^*)f\big)(y)\endalign
$$
for each $f\in\CC_W$\,.
\qed\enddemo

\smallskip

\subheading{Remarks 4.4}

(i) In the special case $\phi=\psi=1$, the accordingly simplified proof
constitutes an alternative proof for Theorem~1.6 from Chapter~1
(recall Remark~1.7).

(ii) In turn, there is also an alternative proof of the above Theorem~4.3
along the lines of the proof of Theorem~1.6 given in Chapter~1.
It involves showing that $F_Z$ (for $Z\ne0$ in the orthogonal complement
of~$T_eW$) induces an isometry from $(M/W,g^W\suboiip)$ to $(M/W,g^W\suboii)$
which, moreover, carries the projected $g\suboiip$@-mean curvature vector
field of the $W$@-orbits to the projected $g\suboii$@-mean curvature
vector field of the $W$@-orbits.

{}From this version of the proof one sees immediately that Theorem~4.3
is actually a special case of a theorem which was proven recently by
C.~Gordon and Z.~Szab\'o (Theorem~1.2 in~\cite{23}). Their theorem
reads like Theorem~1.3 by C.~Gordon in Chapter~1, with the following changes:
Condition~(i) is dropped; in condition~(ii), ``isospectral'' is replaced
by ``isometric'', and moreover it is required that there exists an isometry
between the quotient manifolds which intertwines the projected mean
curvature vector fields of the $W$@-orbits.

(iii) In \cite{23}, C.~Gordon and Z.~Szab\'o apply the theorem mentioned
in~(ii), and a version of it for the case of manifolds with boundary, to
construct
a specific class of interesting new examples of isospectral, locally
non-isometric manifolds which arise as torus bundle whose fibers
are not totally
geodesic. The manifolds are diffeomorphic to products of spheres with tori,
resp\. balls (or bounded cylinders) with tori. In the case of manifolds
with boundary, they obtain continuous isospectral families of negatively
curved manifolds, which contrasts with the spectral rigidity result by
C.~Croke and V.~Sharafutdinov for closed negatively curved
manifolds~\cite{14}.

Without going into detail, we mention here that those examples which
Gordon and Szab\'o construct in the case of closed manifolds
can also be viewed as arising from our (more special) Theorem~4.3;
more precisely, from Proposition~4.5 below.

(iv) As we mentioned at the beginning of this chapter, the only examples
of conformally equivalent, isospectral manifolds which were previously
known had been given by R.~Brooks and C.~Gordon in~\cite{6}.
We note here, again without giving details, that those examples too
can be interpreted as an application of Theorem~4.3.

\bigskip

We finish this preparatory section by specializing Theorem~4.3 to
the case of products; the following proposition is related to
Proposition~1.8 in the same way as Theorem~4.3 is to Theorem~1.6.
We use Notation~4.1(ii).

\bigskip

\proclaim{Proposition 4.5}
Let $(N,h)$ be a closed Riemannian manifold and $H$ be a torus equipped
with an invariant metric. Let $\h=T_eH$ and $\la,\lap$ be two $\h$@-valued
$1$@-forms on~$N$. Let $\phi,\psi,\phip,\psip\in C^\infty(N,\R_+)$.
Assume:
\roster
\item"($*6$)" For every $Z\in\h$ there exists an isometry
  $f_Z$ of $(N,h)$ which satisfies\newline
  $\lap_Z=f_Z^*\la_Z$\,, $\;\phip=f_Z^*\phi$, \;and $\;\psip=f_Z^*\psi$.
\endroster
Then $(N\times H,g\sublii)$ and $(N\times H,g\subliip)$ are isospectral.
\endproclaim

\smallskip

\demo{Proof}
The connection forms $\om,\omp$ associated with $\la,\lap$ satisfy
condition~\thetag{$*5$} from Theorem~4.3. In fact, $F_Z:=(f_Z\,,\Id):
N\times H\to N\times H$ has all the properties required there.
\qed\enddemo

\bigskip

\smallskip

\subhead 4.2\ \ Conformally equivalent examples\endsubhead
\bigskip

\noindent
We will now apply Proposition~4.5 to construct the first examples of
isospectral manifolds which are conformally equivalent and not locally
isometric. The idea is  to find, in the context of Proposition~4.5,
a situation where $\la=\lap$ and nevertheless there exists a pair
of functions $\phi=\psi$ and $\phip=\psip$ such that condition~\thetag{$*6$}
is satisfied {\it nontrivially\/}; that is, the isometries~$f_Z$ cannot
be chosen independently of~$Z$.

\bigskip

\subheading{Example 4.6}
Let $K$ be a compact Lie group with Lie algebra $\k=T_eK$, and let~$h$
be a bi-invariant metric on~$K$. Let $H$ be a torus with Lie algebra
$\h=T_eH$, equipped with an invariant metric. Let $\la,\lap:\k\to\h$
be two linear maps which satisfy condition~\thetag{$*4$} from
Corollary~3.2;
i.e., for each $Z\in\h$ there exists $a_Z\in K$ such that $\lap_Z=
\Ad_{a_Z}^*\la_Z$\,. We endow $K\times K$ with the product metric~$\bar h$,
and define
$$\bar\la:\k\oplus\k\ni(X,Y)\mapsto\la(X)+\lap(Y)\in\h.
$$
Choose a class function $\phi\in C^\infty(K,\R_+)$ (i.e., one that is
invariant under inner automorphisms of~$K$),
and define $\phi_1\,,\phi_2:K\times K\to\R_+$ by
$$\phi_1(x,y):=\phi(x),\quad\phi_2(x,y):=\phi(y)
$$
for all $x,y\in K$.
Denote the lifts of $\phi_1\,,\phi_2$ to $K\times K\times H$
by $\tilde\phi_1$ and $\tilde\phi_2$\,, respectively.
We claim that the conformally equivalent manifolds
$$(K\times K\times H,\,\tilde\phi_1\,g_{\bar\la})\text{\ \ and\ \ }
  (K\times K\times H,\,\tilde\phi_2\,g_{\bar\la})
$$
are isospectral by Proposition~4.5, where~$g_{\bar\la}$ is the
left invariant metric on $K\times K\times H$ associated with~$\bar h$
and~$\bar\la$ as in Corollary~3.2.

In fact, for any given $Z\in\h$ choose~$a_Z\in K$ such that $\lap_Z
=\Ad_{a_Z}^*\la_Z$\,, and define
$$f_Z:K\times K\ni(x,y)\mapsto(I_{a_Z}(y),I_{a_Z}\inv(x))\in K\times K,
$$
where $I_{a_Z}$ denotes conjugation by~$a_Z$\,.
Then~$f_Z$ is an isometry by the bi-invariance of~$\bar h$; moreover,
$$\align (f_Z^*\bar\la_Z)(X,Y)&=\la_Z(\Ad_{a_Z}(Y))+\lap_Z(\Ad_{a_Z}\inv(X))
  =\lap_Z(Y)+\la_Z(X)\\&=\bar\la_Z(X,Y)\endalign
$$
for all $(X,Y)\in\k\oplus\k$. Finally, $(f_Z^*\phi_2)(x,y)=
\phi_2(I_{a_Z}(y),I_{a_Z}\inv(x))=\phi(x)=\phi_1(x,y)$ for all $x,y\in K$
since~$\phi$ is a class function. So $f_Z$ satisfies all the conditions
from Proposition~4.5, with~$\bar\la$ playing the role of both~$\la$
and~$\lap$ from the proposition, $\phi_1$~playing the role of $\phi=\psi$,
and~$\phi_2$ the role of $\phip=\psip$.
Recall from Remark~4.2 that $g_{\bar\la,\phi_1,\phi_1}=
\tilde\phi_1\,g_{\bar\la}$ and $g_{\bar\la,\phi_2,\phi_2}=\tilde\phi_2
\,g_{\bar\la}$\,.

\bigskip

The following proposition shows that in many cases the two resulting
isospectral, conformally equivalent manifolds are not locally isometric.

\bigskip

\proclaim{Proposition 4.7}
In the context of Example~{\rm 4.6} assume that
\roster
\item"(i)" $K=\SO(m)$, $m\ge5$, or
\item"(ii)" $K=\SU(m)$, $m\ge3$,
\endroster
and that the torus~$H$ is two-dimensional.
For the isospectral pair $\la,\lap:\k\to\h$ {\rm(}as in the example\/{\rm)}
assume that
$\|\Ric^{g\subla}\|^2\ne\|\Ric^{g\sublap}\|^2$,
where $g\subla$ and $g\sublap$ are the associated left invariant metrics
on $K\times H$. Finally, define a class function~$\phi$ on~$K$ by
$$\phi(x):=e^{2\eps\,\tr(x)/m}\,,
$$
where $\tr$ denotes the trace on $(m\times m)$@-matrices, and
$0<\eps<1/8$.
Then the conformal metrics $\tilde\phi_1\,g_{\bar\la}$ and
$\tilde\phi_2\,g_{\bar\la}$ on $K\times K\times H$, defined as in
Example~{\rm4.6}, are not locally isometric.

More precisely, the preimage in $K\times K\times H$ of the maximal scalar
curvature of~$\tilde\phi_1\,g_{\bar\la}$\,,
resp\. of~$\tilde\phi_2\,g_{\bar\la}$\,, is a submanifold which, when
endowed with the induced metric, is isometric to
$(K\times H,e^{2\eps}g\sublap)$,
respectively $(K\times H,e^{2\eps}g\subla)$,
whose Ricci tensors have different norms by the choice of~$\la$ and~$\lap$.
\endproclaim

\smallskip

Recall that if $K$ and $H$ are of the above type then isospectral pairs
$\la,\lap:\k\to\h$ with $\|\Ric^{g\subla}\|^2\ne
\|\Ric^{g\sublap}\|^2$ exist indeed, even many of them.
See the examples~3.3\,/\,3.7 in connection with Proposition~3.15.

\bigskip

\demo{Proof}
Let $\Id$ denote the neutral element of $K=\SO(m)$, resp\.
$K=\SU(m)$, where $m\ge5$, resp\. $m\ge3$.
We claim that for our choice of~$\phi$ the preimage of the
maximal scalar curvature of
$(K\times K\times H,\,\tilde\phi_1\,g_{\bar\la})$, resp\.
$(K\times K\times H,\,\tilde\phi_2\,g_{\bar\la})$, is precisely
$$\{\Id\}\times K\times H,\text{\quad resp.\quad}
  K\times\{\Id\}\times H.
$$
These submanifolds, endowed with the metric induced by $\tilde\phi_1
\,g_{\bar\la}$\,, resp\. $\tilde\phi_2\,g_{\bar\la}$\,, are isometric
to $(K\times H,e^{2\eps}g\sublap)$, resp\. $(K\times H,e^{2\eps}
g\subla)$, by the definition of~$\bar\la$ and the fact that
$\phi_1(\Id,x)=\phi_2(x,\Id)=\phi(\Id)=e^{2\eps}$ for all $x\in K$.
The statement of the proposition will thus follow.

We now prove our above claim. Since $K$ is irreducible by assumption,
the bi-invariant metric~$h$ is a scalar multiple of the Killing metric
$-B:(X,Y)\mapsto-\tr(\ad_X\ad_Y)$; let $c\in\R$ be such that
$h=c^2\cdot(-B)$. Define $\tau(x):=\eps\,\tr(x)/m$,
$\;\tau_1(x,y):=\tau(x)$
for $(x,y)\in K\times K$, and let $\tilde\tau_1$ be the lift of~$\tau_1$
to $K\times K\times H$; we thus have $\tilde\phi_1=\exp(2\tilde\tau_1)$.
Note that $(\Delta_{g_{\bar\la}}\tilde\tau_1)(x,y,z)=(\Delta_h\tau)(x)
=\mu\,\tau(x)=\mu\,\tilde\tau_1(x,y,z)$
for all $(x,y,z)\in K\times K\times H$, where
$$\mu=\frac1{c^2}\cdot\frac1m\cdot\frac1{k(m)}\cdot\dimm\,K
$$
with $k(m)=m-2$ in case~(i) and $k(m)=2m$ in case~(ii). By \cite{5},
Theorem~1.159 we have
$$\scal^{\tilde\phi_1 g_{\bar\la}}=\frac{\alpha+\beta\tilde\tau_1
  -(N-1)(N-2)\|d\tilde\tau_1\|^2_{g_{\bar\la}}}{\exp(2\tilde\tau_1)}\,,
$$
where $N:=\dimm(K\times K\times H)=2\,\dimm\,K+2$, $\;\alpha:=\scal^{g_{\bar
\la}}$, and $\beta:=2(N-1)\mu$.
Note that $\alpha\le\scal^{\bar h}=2\,\scal^h=\frac2{c^2}\cdot\frac
{\dimm\,K}4$\,, where the inequality follows from Proposition~2.8.

We observe that the function $\R\ni s\mapsto\frac{\alpha+\beta s}
{\exp(2s)}\in\R$ is strictly monotonously increasing in $s\in(-\infty,
\frac12-\frac\alpha\beta]$.
The image of our~$\tilde\tau_1$ is contained in this interval:
In fact, we obviously have $\max\tilde\tau_1=\max\tau=\eps$ and
$$\frac\alpha\beta\le\frac{\dimm\,K\,/\,2c^2}{2\,(2\,\dimm\,K+1)\,\dimm\,K
  \,/\,c^2m\,k(m)}\le\frac{m\,k(m)}{8\,\dimm\,K}\,,
$$
which in case~(i) equals $\frac{2m(m-2)}{8m(m-1)}\le\frac14<\frac12-
\eps$\,, and in case~(ii) equals $\frac{2m^2}{8(m^2-1)}$ which
for $m\ge3$ is also smaller than $\frac12-\eps$ by the choice of~$\eps$.

Since additionally we have $d\tilde\tau_1\restr p=0$ if $\tilde\tau_1(p)$
is maximal, we conclude that $\scal^{\tilde\phi_1g_{\bar\la}}$
attains its maximum precisely in those points where~$\tilde\tau_1$
does so; namely, in $\{\Id\}\times K\times H$, as claimed.
In the same way we show that $\scal^{\tilde\phi_2g_{\bar\la}}$
attains its maximum precisely in $K\times\{\Id\}\times H$.
\qed\enddemo

\enddocument